\documentclass[11pt]{amsart}

\usepackage{amsmath}
\usepackage[all,poly,knot,dvips]{xy}
\usepackage{pstricks,pst-poly,pst-node,pstcol}

\usepackage{amssymb,latexsym}

\usepackage{amsthm,latexsym}
\usepackage{eucal,latexsym}

 \newtheorem{thm}{Theorem}[subsection]
 \newtheorem{cor}[thm]{Corollary}
 \newtheorem{lem}[thm]{Lemma}
 \newtheorem{prop}[thm]{Proposition}
 
 \theoremstyle{definition}
 \newtheorem{defn}[thm]{Definition}
 \theoremstyle{remark}
 \newtheorem{rem}[thm]{Remark}
 \numberwithin{equation}{subsection}

\newcommand{\CS}{{\mathcal{S}}}

\newcommand{\CC}{{\mathcal{C}}}

\newcommand{\BC}{\mathbb{C}}

\newcommand{\CM}{{\mathfrak{M}}}
\newcommand{\s}{\sigma}
\newcommand{\B}{\beta}
\newcommand{\Au}{\text{Aut}}

\begin{document}

 \title[Moves]{On moves between branched coverings of $S^3$\\
the case of Four sheets}

\author{Nikos Apostolakis}

\address{Department of Mathematics\\
 University of California \\
Riverside CA 92521}

\email{nea@math.ucr.edu}

\thanks{A major part of this work was completed
at the Graduate Center of CUNY as part of the authors doctoral dissertation.}

\subjclass{Primary 57M12; Secondary 57M25}

\keywords{branched covering, covering move, colored braid, colored link, 3-manifold}

\begin{abstract}
A combinatorial presentation of closed orientable $3$-manifolds
as bi-tricolored links is given together with two
versions of a calculus via moves to manipulate bi-tricolored
links without changing the represented manifold. That is,
we provide a finite set of moves sufficient to relate  any two manifestations
of the same $3$-manifold as a simple $4$-sheeted branched covering of $S^3$.
\end{abstract}

\maketitle

 \setcounter{section}{-1} \tableofcontents

\section{Introduction}

It is known that every oriented $3$-dimensional manifold
manifests itself as (the total space of) a branched
covering over the $3$-dimensional sphere $S^3$, with simple branching
 behaviour~(\cite{A}). Furthermore it is
known that the degree of the covering can chosen to be any natural
number greater or equal to three(\cite{H},
\cite{M1}). This fact leads to a combinatorial presentation of $3$-manifolds
via colored link diagrams, that is
diagrams of (planar projections of) links whose arcs are
labeled by transpositions of some symmetric group. The
following question then
 arises: \\

\begin{center} \textit{Find a finite set of (preferably local) moves,
 such that any two presentations of the same $3$-manifold as a
colored link can be related by a finite sequence of moves from this set.}
\end{center}

 By a local move we vaguely mean a move that applies to a relatively small
part of the diagram independently
of the structure of the link outside that part. We don't attempt to give
a more precise meaning to the term local
move.\

 In this study we concentrate on coverings of degree $4$ and give two answers:\

Theorem~\ref{theoremX} gives a set of five moves (two local and three non local)
sufficient to relate any two
presentations of the same manifold as a $4$-sheeted simple branched covering of the $3$-sphere.\

Theorem~\ref{stab} proves that after adding a fifth sheet in a standard way the
two local moves of
Theorem~\ref{theoremX} suffice.\

These results and their proofs are modeled on similar results of Piergallini
for coverings of degree $3$
(\cite{P1}, \cite{P2}). It is unknown at this time if similar results are
true for coverings of degree greater or
equal to $5$.\

The rest of this work is organized as follows:\\

 Section~\ref{prem} sets the stage: Basic definitions are
 given and basic results are stated (and sometimes proved) to be
 used in the later sections.\

 Section~\ref{3sheets} reviews, for the sake of completeness, what
 is known for simple branched coverings of degree $3$. We emphasize that
 all results in this section have been discovered by other
 people. In particular subsection~\ref{2sphere3} summarizes the
 relevant results of \cite{BW1} and subsection~\ref{3sphere3}
 the results of \cite{P1} and \cite{P2}.\

 In Section~\ref{2sphere4}  the main technical results of this
 work are proved. We study simple $4$-sheeted coverings of the
 $2$-sphere and the braid group action on them. This section
 culminates with Theorem~\ref{kernel} in which normal generators
 for the kernel of the quotioned lifting homomorphism are given.\

 In Section~\ref{3sphere4}  the above mentioned results about moves are proved.\

 Section~\ref{questions} concludes this work by asking some questions that arise
 in the course of our study.\

 Finally, to fix conventions etc.,  we sketch in four
 Appendices background material that is used throughout this work.

\section{Preliminaries}\label{prem} In what follows we always work in the PL
category. Thus manifold means  PL manifold, submanifold means
locally flat submanifold, homeomorphism means  PL homeomorphism,
etc.
\subsection{Generalities}  Lets begin with the basic
definitions:
 \begin{defn}
A {\em branched covering} $p:E(p)\longrightarrow B(p)$ is a surjective
map between manifolds such that there is a codimension-$2$
submanifold $L$ of $B(p)$ with the property that $p$ restricted to
the preimage of $B(p)\setminus L$ is a finite covering. $p$ is
called the (covering) projection, $E(p)$ the total space, $B(p)$
the base space and $L$
the (downstairs) branching locus.\\
 If the base space $B$ has a basepoint $*\notin L$ then an $m$-sheeted {\em pointed}
 branched
covering of $B$ is an $m$-sheeted branched covering together with
an identification of $p^{-1}(*)$ with $\{1,\dotsc,m\}$.\\
 An {\em equivalence} between  branched coverings is a commutative diagram
 $$
\xymatrix{
 {E(p_1)} \ar[d]_{p_1} \ar[r]^{f} & {E(p_2)} \ar [d]^{p_2}\\
 {B(p_1)}  \ar[r]_{\bar{f}} & {B(p_2)}        }
$$
where $f$, $\bar{f}$ are  homeomorphisms. $\bar{f}$ is said to
cover , or to be over, or to lift $f$.\\
An equivalence between pointed branched coverings is a commutative
diagram as above with $\bar{f}$ basepoint preserving and $f$
commuting with the identifications of the fibers with
$\{1,\dotsc,m\}$. Such an $f$ is called {\em pointed}.\\
 An {\em isomorphism} of (pointed) branched coverings is an equivalence over the
identity.
\end{defn}

\begin{rem} A much more general definition of branched coverings can be found
in~\cite{F}.
\end{rem}

\begin{prop}\label{equiv} Let $M$ be a manifold and $L$ a
codimension-$2$ submanifold of M. Then:\

 Equivalence classes of
pointed $m$-sheeted coverings of $M$ branched over $L$ are in
bijection with $Hom(\pi_1(M\setminus L,*),\CS_m)$, where $\CS_m$
denotes the symmetric group on $m$ symbols.\

 Equivalence
classes of (unpointed) $m$-sheeted coverings of $M$ branched over
$L$ are in bijection with $Hom(\pi_1(M\setminus
L,*),\CS_m)/\CS_m$, where $\CS_m$ acts by conjugation in the
range.\

 The pullback of a (pointed) branched covering by a homeomorphism is again
a (pointed) branched covering branched over the preimage of the branching
locus. The above bijections are equivariant with
respect to pullbacks.
\end{prop}

\begin{proof}Fox proved in \cite{F} that a branched covering is determined, up
to equivalence, by its restriction outside the branching locus and
that any finite covering of $B\setminus L$ can be extended to a
branched covering of $B$. Applying the standard theory of covering
spaces one gets the result.
\end{proof}

   The homomorphism $\rho:\pi_1(M\setminus
L,*)\longrightarrow\CS_m$ corresponding to $p$ is called the {\em
monodromy} of $p$. From now on a branched covering will
systematically be confused with its monodromy in the sense that
the same symbol will be used for both.\

 The following holds:
\begin{lem}\label{lift} Let $p:\widetilde M\longrightarrow M$ be a pointed branched
covering and $h:M\longrightarrow M$ a based homeomorphism that
fixes the branching locus as a set. Then $h$ lifts to a pointed
equivalence of $\tilde h:\widetilde M\rightarrow\widetilde M$ iff
the pullback $h^*(p)$ is isomorphic  to $p$.
\end{lem}
\begin{proof}
Consider the following commutative diagram:
 $$\xymatrix{
 {h^*(\widetilde M)}\ar[d]_{h^*p} \ar[r]_{h^*} \ar@/^1pc/@{.>} [rr]^{\phi}
  & {\widetilde M}\ar[d]_{p} \ar@{.>}[r]_{\tilde h^{-1}} & {\widetilde M}\ar[d]^{p}\\
 {M}\ar[r]^{h} \ar@/_1pc/ [rr]_{\text{id}} & {M} \ar[r]^{h^{-1}} &{M}
}$$
 Obviously one of the doted arrows exists iff  the other does.
\end{proof}
Notice that $\tilde h$ is unique (when it exists).

\begin{defn} An $m$-sheeted branched covering is called simple if the preimage
of any point has cardinality at least $m-1$.
\end{defn}

\subsection{Branched coverings over $D^2$}\label{2D} From now on, unless
explicitly mentioned otherwise, covering will mean a {\em
branched} covering. Although we are mainly interested in coverings
over the sphere $S^2$ it is convenient to consider coverings over
the disc $D^2$.\

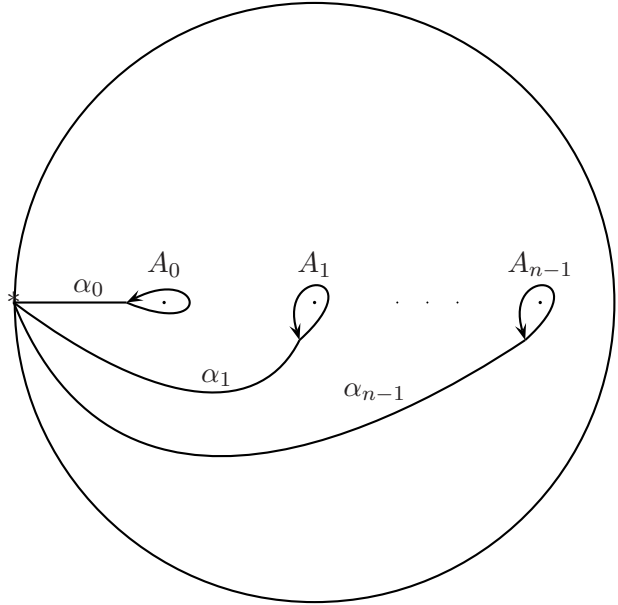
\begin{figure}[htp]
\begin{pspicture}(-1,-1)(7,7)
 \pscircle(3,3){4}
 \rput(-1,3){*}
 \psdots[dotstyle=*,dotscale=.5](1,3)(3,3)(6,3)
 \psline(-1,3)(.5,3)
 \psbezier{->}(.5,3)(1.7,2.5)(1.5,3.5)(.5,3)
 \psbezier
 (-1,3)(1,1.5)(2.3,1.5)(2.8,2.5)
\psbezier{->}(2.8,2.5)(3.8,3.4)(2.5,3.5)(2.8,2.5)
 \psbezier
 (-1,3)(.4,-.5)(4,1.3)(5.8,2.5)
 \psbezier{->}(5.8,2.5)(6.8,3.4)(5.5,3.5)(5.8,2.5)
 \rput(1,3.5){$A_0$}
 \rput(3,3.5){$A_1$}
 \rput(6,3.5){$A_{n-1}$}
 \rput(0,3.2){$\alpha_0$}
 \rput(1.7,2){$\alpha_1$}
 \rput(3.8,1.8){$\alpha_{n-1}$}
 \psdots[dotstyle=*,dotscale=.3](4.1,3)(4.5,3)(4.9,3)
\end{pspicture}
\caption{The free generators of $\pi_1(D^2\setminus
L,*)$}\label{sperma}
\end{figure}

 To set the stage consider $D^2$ as the unit disk in $\BC$ with
basepoint $*=-1$. The branching locus $L$ of a covering will be a
finite set of points, $L=\{A_0,A_1,\dotsc,A_{n-1}\}$, which are
assumed to lie in the real axis in that order. $\pi_1(D^2\setminus
L,*)$ is a free group on the generators $\alpha_0$, $\alpha_1$,
$\dotsc,\alpha_{n-1}$, where $\alpha_i$ is represented by a {\em
lasso} based at $*$ and going around $A_i$ in the counterclockwise
direction as shown in Figure~\ref{sperma}.

 Then using Proposition~\ref{equiv} one can represent an $m$-sheeted
pointed covering branched over $L$ as a sequence of permutations
$\sigma_0$, $\sigma_1$,$\dotsc,\sigma_{n-1}$ of $\CS_m$, where
such a sequence represents the covering
 $$\xymatrix@1{ {\rho:\pi_1(D^2\setminus L,*)}\ar[r] & {\CS_m} }$$
 $$\rho(\alpha_i)=\sigma_i \qquad {\text{for }} i=0,\dotsc,n-1$$
Such a sequence represents a simple covering iff all $\sigma_i$'s
are transpositions.Unless explicitly mentioned otherwise {\em all
coverings are assumed simple and pointed} from now on.\

\noindent {\bf Coverings of $S^2$}\label{2sphere}
The number of boundary components of the total space of $\rho$
equals the number of cycles (including cycles of length $1$) in a
decomposition of
$$\sigma=\prod_{i=0}^{n-1}\sigma_i$$
into disjoint cycles. Indeed $\sigma$ is the monodromy around the
(positively oriented) boundary.\\
 In particular if $\sigma=\text{id}$, $E(\rho)$ has $m$ boundary
components, each mapping homeomorphically onto $\partial D^2$.
Thus one can ``fill them in'' with $m$ discs to get a covering of
$$S^2=D^2\bigcup_{\partial D^2}D^2$$
with each filled in disk mapping homeomorphically onto the second
$D^2$.\\
 Conversely given a covering  of $S^2$ one can cut off a disk not
containing any branch values to get a covering of a disk with
boundary monodromy equal to id.\\
 Thus coverings of $D^2$ with identity boundary monodromy can (and will) be
identified with coverings of $S^2$.\

 If $\rho$ is a covering of $S^2$ the Euler characteristic of $E(\rho)$ can be
calculated by choosing compatible triangulations upstairs and
downstairs. One gets the following version of the Riemann-Hurwitz
formula:
 \begin{equation}\label{char}
 \chi(E(\rho))=2m-n
\end{equation}
where $m$ is the degree and $n$ the number of branch values. As a
corollary the number of branch values is even.  In terms of genus
(assuming that $E(\rho)$ is connected):
\begin{equation}\label{genus}
g(E(\rho))=\frac{n}{2}-m+1\quad .
\end{equation}

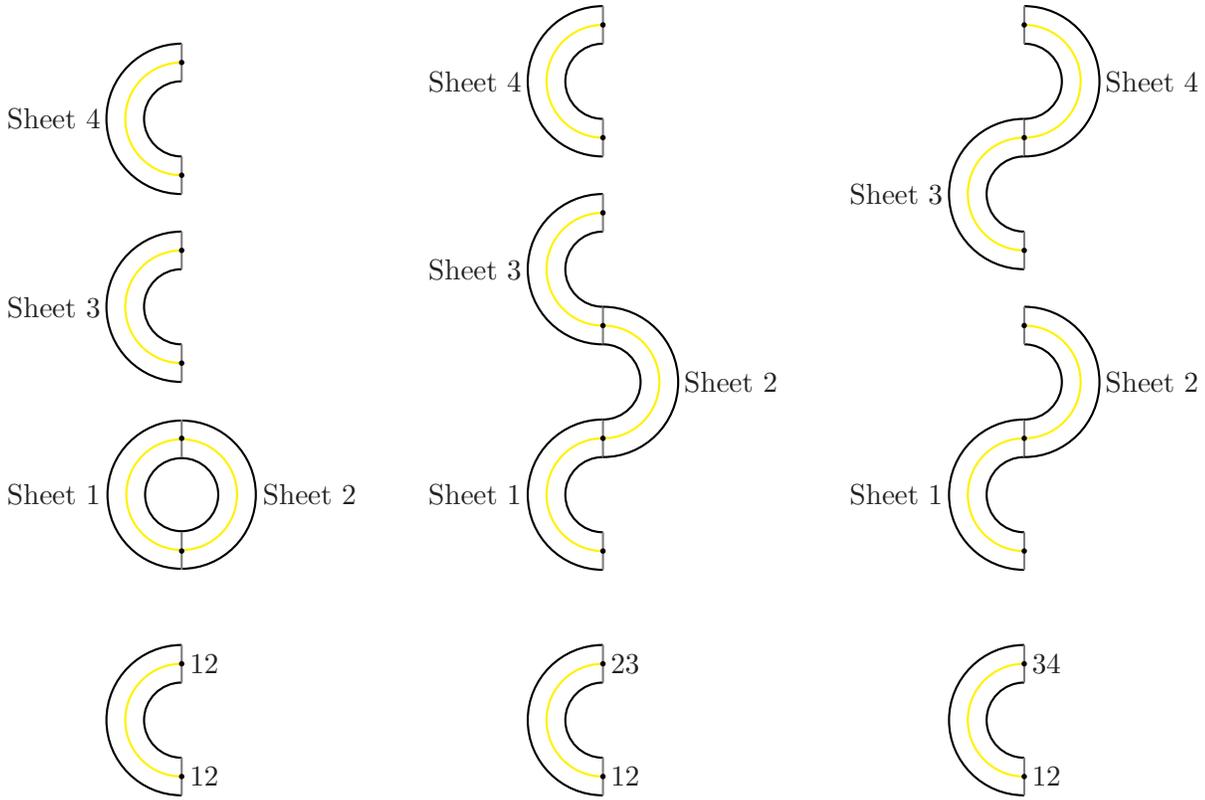
\begin{figure}[htp]
 \[
 \begin{pspicture}(3.5,11)
 \rput(-.7,4){Sheet $1$}
 \rput(2.7,4){Sheet $2$}
 \rput(-.7,6.5){Sheet $3$}
\rput(-.7,9){Sheet $4$}
 \psarc(1,1){1}{90}{270}
 \psarc[linecolor=yellow](1,1){.75}{90}{270}
 \psarc(1,1){.5}{90}{270}
 \psline[linecolor=gray](1,2)(1,1.5)
 \psline[linecolor=gray](1,.5)(1,0)
 \rput(1.3,.25){12}
 \rput(1.3,1.75){12}
    \pscircle(1,4){1}
    \pscircle[linecolor=yellow](1,4){.75}
    \pscircle(1,4){.5}
     \psline[linecolor=gray](1,5)(1,4.5)
     \psline[linecolor=gray](1,3.5)(1,3)
 \psarc(1,6.5){1}{90}{270}
 \psarc[linecolor=yellow](1,6.5){.75}{90}{270}
 \psarc(1,6.5){.5}{90}{270}
 \psline[linecolor=gray](1,7.5)(1,7)
 \psline[linecolor=gray](1,6)(1,5.5)
  \psarc(1,9){1}{90}{270}
  \psarc[linecolor=yellow](1,9){.75}{90}{270}
  \psarc(1,9){.5}{90}{270}
  \psline[linecolor=gray](1,10)(1,9.5)
  \psline[linecolor=gray](1,8.5)(1,8)
     \psdots[dotsize=.1](1,.25)(1,1.75)(1,3.25)(1,4.75)(1,5.75)(1,7.25)(1,8.25)(1,9.75)
 \end{pspicture}
 \hspace{2.1cm}
 \begin{pspicture}(3.5,11)
 \rput(-.7,4){Sheet $1$}
 \rput(2.7,5.5){Sheet $2$}
 \rput(-.7,7){Sheet $3$}
\rput(-.7,9.5){Sheet $4$}
 \psarc(1,1){1}{90}{270}
 \psarc[linecolor=yellow](1,1){.75}{90}{270}
 \psarc(1,1){.5}{90}{270}
 \psline[linecolor=gray](1,2)(1,1.5)
 \psline[linecolor=gray](1,.5)(1,0)
 \rput(1.3,.25){12}
 \rput(1.3,1.75){23}
    \psarc(1,4){1}{90}{270}
    \psarc[linecolor=yellow](1,4){.75}{90}{270}
    \psarc(1,4){.5}{90}{270}
     \psline[linecolor=gray](1,5)(1,4.5)
     \psline[linecolor=gray](1,3.5)(1,3)
 \psarc(1,7){1}{90}{270}
 \psarc[linecolor=yellow](1,7){.75}{90}{270}
 \psarc(1,7){.5}{90}{270}
 \psline[linecolor=gray](1,8)(1,7.5)
 \psline[linecolor=gray](1,6.5)(1,6)
  \psarc(1,5.5){1}{270}{90}
  \psarc[linecolor=yellow](1,5.5){.75}{270}{90}
  \psarc(1,5.5){.5}{270}{90}
\psarc(1,9.5){1}{90}{270}
 \psarc[linecolor=yellow](1,9.5){.75}{90}{270}
 \psarc(1,9.5){.5}{90}{270}
 \psline[linecolor=gray](1,10.5)(1,10)
 \psline[linecolor=gray](1,9)(1,8.5)
 \psdots[dotsize=.1](1,.25)(1,1.75)(1,3.25)(1,4.75)(1,6.25)(1,7.75)(1,8.75)(1,10.25)
 \end{pspicture}
 \hspace{2.1cm}
 \begin{pspicture}(2.5,11)
 \rput(-.7,4){Sheet $1$}
 \rput(2.7,5.5){Sheet $2$}
 \rput(-.7,8){Sheet $3$}
\rput(2.7,9.5){Sheet $4$}
 \psarc(1,1){1}{90}{270}
 \psarc[linecolor=yellow](1,1){.75}{90}{270}
 \psarc(1,1){.5}{90}{270}
 \psline[linecolor=gray](1,2)(1,1.5)
 \psline[linecolor=gray](1,.5)(1,0)
 \rput(1.3,.25){12}
 \rput(1.3,1.75){34}
    \psarc(1,4){1}{90}{270}
    \psarc[linecolor=yellow](1,4){.75}{90}{270}
    \psarc(1,4){.5}{90}{270}
     \psline[linecolor=gray](1,5)(1,4.5)
     \psline[linecolor=gray](1,3.5)(1,3)
 \psarc(1,5.5){1}{270}{90}
 \psarc[linecolor=yellow](1,5.5){.75}{270}{90}
 \psarc(1,5.5){.5}{270}{90}
 \psline[linecolor=gray](1,6.5)(1,6)
     \psarc(1,8){1}{90}{270}
    \psarc[linecolor=yellow](1,8){.75}{90}{270}
    \psarc(1,8){.5}{90}{270}
     \psline[linecolor=gray](1,9)(1,8.5)
     \psline[linecolor=gray](1,7.5)(1,7)
 \psarc(1,9.5){1}{270}{90}
 \psarc[linecolor=yellow](1,9.5){.75}{270}{90}
 \psarc(1,9.5){.5}{270}{90}
 \psline[linecolor=gray](1,10.5)(1,10)
 \psdots[dotsize=.1](1,.25)(1,1.75)(1,3.25)(1,4.75)(1,6.25)(1,7.25)(1,8.75)(1,10.25)
 \end{pspicture}
 \]
\caption{How various intervals lift}\label{fig:intlift}
\end{figure}

\noindent {\bf Constructing branched coverings by cutting and pasting.}
  There is a well known method (going back to Riemann) for
  getting models for coverings over
$D^2$ by ``cutting and pasting''. Let $\rho=\s_0$, $\s_1$,
$\dotsc,\s_{n-1}$ be a, not necessarily simple, $m$-sheeted
covering of $D^2$ branched over $L$. For $i=0,1,\dotsc,n-1$ take
an arc $\gamma_i$ connecting the branching value $A_i$ to a
boundary point (other than  the basepoint $*$), in such a way that
two different $\gamma_i$'s are disjoint. Cut the disc along these
arcs to get a new disk with $2n$ distinguished arcs on its
boundary: each $\gamma_i$ gives two arcs joined at $A_i$, call
these two arcs $\gamma_i^0$ and $\gamma_i^1$. Now take $m$ copies
of the cut disc (the sheets of the covering) numbered $1$ through
$m$ and glue them along their distinguished arcs according to the
$\s_i$'s, that is for $i=0,1,\dotsc,n-1$ and $j=1,\dotsc,m$ glue
the $j$-th copy of $\gamma_i^0$ to the $\s_i(j)$-th copy of
$\gamma_i^1$. The result of these gluings is the total space of
the covering. The projection is the obvious one. See
Figure~\ref{fig:intlift} for an example of this construction.\

Notice that since the cuts described above make $D^2\setminus L$
\label{cutpaste} simply connected (in fact contractible)  this
method can be used to construct a model for {\em any} covering of
$D^2$ branched over $L$. However in order to construct a given
covering one can cut along any $1$-dimensional subcomplex of
$D^2\setminus L$ as long as the fundamental group of the cut disk
is send to identity by the monodromy of the covering. This
observation will be used later when we will construct specific
coverings of the $2$-dimensional sphere.\

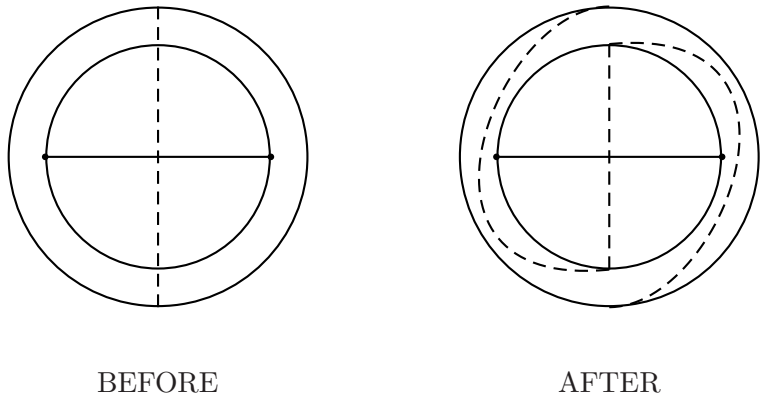
\begin{figure}[htp]
\[
\begin{pspicture}(0,-1)(4,4)
 \pscircle(2,2){2}
 \pscircle(2,2){1.5}
 \psdots(.5,2)(3.5,2)
 \psline(.5,2)(3.5,2)
 \psline[linestyle=dashed](2,4)(2,0)
 \rput(2,-1){BEFORE}
\end{pspicture}
\hspace{2cm}
 \begin{pspicture}(0,-1)(4,4)
 \pscircle(2,2){2}
 \pscircle(2,2){1.5}
 \psdots(.5,2)(3.5,2)
 \psline(.5,2)(3.5,2)
 \psline[linestyle=dashed](2,.5)(2,3.5)
 \psbezier[linestyle=dashed]
 (2,4)(.5,4)(-1,.2)(2,.5)
 \psbezier[linestyle=dashed]
 (2,0)(3.5,0)(5,3.8)(2,3.5)
  \rput(2,-1){AFTER}
\end{pspicture}
\]
\caption{Rotation around an interval.}\label{rot}
\end{figure}

\noindent {\bf The Braid group action.} See
 Appendix~\ref{braids} for the definition of the braid group
 $B_n$.
 It is known that the group $\CM_{0,1,n}$ of isotopy classes of
homeomorphisms of $(D^2,L,\text{rel}\quad \partial D^2)$, is
isomorphic with $B_n$ the braid group on $n$ strands (see for
example \cite{B}). The isomorphism is given by sending the
generator $\B_i$ of $B_n$ to (the isotopy class of) the rotation
around the interval $x_i$ for $i=0,\dotsc,n-1$, where $x_i$ is the
interval between $A_i$ and $A_{i+1}$ on the real axis. In general
\begin{defn} An interval is a path $x$ in the interior of $D^2$
with endpoints in $L$ but otherwise missing $L$.\\
 A rotation
around $x$ is a homeomorphism that fixes pointwise the complement
of a disc neighborhood $U$ of $x$ in $D^2\setminus L$ while
rotating the interior of $U$ by $180^o$ counterclockwise, mapping
$x$ to itself with its endpoints reversed (see Figure~\ref{rot}).
The class in $\CM_{0,1,n}$ of a rotation around an interval $x$
depends only on the isotopy class rel endpoints of $x$ and is
denoted by $\B_x$.
\end{defn}

 A braid will be confused with its image on $\CM_{0,1,n}$ but all
the induced actions of $B_n$ are considered {\em right} actions.
In particular $B_n$ acts on the right on\

\begin{itemize}
\item $\pi_1(D^2\setminus L,*)$.\\
 This action is given by
$$ (\alpha_j)\B_i=\begin{cases} \alpha_j &\text{if $j\neq
i,i+1$},\\
  \alpha_{i+1} & \text{if $j=i$},\\
 \alpha_{i+1}\alpha_i\alpha_{i+1}^{-1} & \text{if $j=i+1$}.
\end{cases}\qquad \text{for}\quad i,j=0,\dotsc,n-1.$$\\

\item Isotopy classes of intervals.\\
This action has the property
$$\B_{(x)\B}=\B^{-1}\B_x\B
\quad \forall \text{ intervals } x,\forall \B\in B_n\quad.$$

\underline{Notation:} For two braids $\B$, $\B'$ denote
$\B^{-1}\B'\B$ by $[\B']\B$, so that the
rotation around the interval $(x)\B$ is denoted by $[\B_x]\B$.\\

\item $Hom(\pi_1(D^2\setminus L,*),\CS_m)$.\\
This action is defined by:
$$((\rho)\B)(\alpha)=\rho(\alpha\B^{-1})$$
and corresponds to pullback of coverings by homeomorphisms (which
act on the left).\\
\end{itemize}

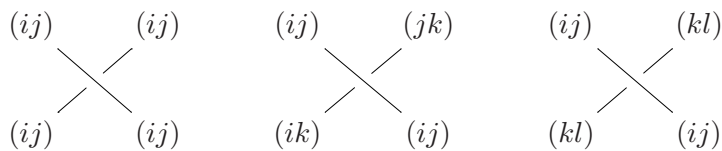
\begin{figure}[htp]
\[
  \xymatrix{
 {(ij)}\ar@{-}[dr] & {(ij)}\ar@{-}[dl]|\hole\\
 {(ij)}            &   {(ij)}  }
\hspace{1cm} \xymatrix{
 {(ij)}\ar@{-}[dr] & {(jk)}\ar@{-}[dl]|\hole\\
 {(ik)}            &   {(ij)}  }
 \hspace{1cm}
\xymatrix{
 {(ij)}\ar@{-}[dr] & {(kl)}\ar@{-}[dl]|\hole\\
 {(kl)}            &   {(ij)}  }
 \]
 \caption{How a generator of $B_n$ acts.}\label{Baction}
 \end{figure}

 This last action of $B_n$ can be expressed combinatorially via
 ``colored braids''. Represent a covering $\rho$ as a ``coloring'' of
the branching locus $L$, i.e. as a labelling of $A_i$'s by
transpositions of the symmetric group $\CS_m$ (referred to as
colors). To see how a braid $\B$ acts on $\rho$ draw a diagram of
$\B$ with its top endpoints coinciding with $L$, and let the
colors of $\rho$ ``flow down'' through the diagram according to
the rule that when a strand passes under another strand its color
gets conjugated by the color of the over strand. That way we get a
new coloring of $L$ at the bottom of the diagram. This bottom
coloring represents $(\rho)\B$.
 Figure~\ref{Baction} shows the action of a generator
in the cases that the top colors coincide, ``intersect'', or are
disjoint, respectively.

\

\noindent {\bf The Category of Colored Braids.}\label{ColBraid}
  The above can be said
  in a categorical language, with
colorings of $L$ being objects and colored braids being morphisms.
Then we have a ``lifting functor'' from the category of colored
braids to the category with objects coverings of $D^2$ branched
over $L$ and morphisms equivalences between them. Indeed it
follows from Lemma~\ref{lift} that (the isotopy class of) the
homeomorphism corresponding to a colored braid lifts to an
equivalence between the coverings corresponding to its ends. In
particular there is a lifting homomorphism
$$\xymatrix@1{ {\lambda(\rho):\Au(\rho)}\ar[r] & {\CM(E(\rho)}  }$$
where $\CM$ denotes the mapping class group and $\Au(\rho)$
denotes the group of automorphisms of $\rho$ as an object in this
category. By Lemma~\ref{lift} $\Au(\rho)$ consists of those braids
that lift to isomorphisms of $\rho$.\

The following lemma describes the ``lifting status'' of rotations
around intervals.
\begin{lem}\label{intlift} Let $x$ be an interval and $\alpha$ an
element of $\pi_1(D^2\setminus L,*)$ represented by a path that
goes around one of the endpoints of $x$ and meets $x$ only once.
Then:
 \begin{itemize}
\item[(a)] If $\rho(\alpha\B_x)=\rho(\alpha)$ then $p^{-1}(x)$
consists of a loop and $m-2$ arcs mapping homeomorphically onto
$x$. A rotation around $x$ will then lift to the composition of a
Dehn twist around the loop and rotations around the arcs, so that
$\B_x$ lifts to a Dhen twist around the loop.
\item[(b)] If $\rho(\alpha\B_x)$ and $\rho(\alpha)$ are
intersecting transpositions then $p^{-1}(x)$ consists of a
``length $3$'' arc and $m-3$ arcs mapping homeomorphically onto
$x$.Then $\B_x$ doesn't lift but $\B_x^3$ lifts to id.
\item[(c)]If $\rho(\alpha\B_x)$ and $\rho(\alpha)$ are
disjoint transpositions then $p^{-1}(x)$ consists of two ``length
$2$'' arcs and $m-4$ arcs mapping homeomorphically onto $x$. Then
$\B_x$ doesn't lift but $\B_x^2$ lifts to id.
 \end{itemize}
\end{lem}
\begin{proof}
Assume without loss of generality, that the subgroup of
$\mathcal{S}_n$ generated by the monodromies of $\alpha$,
$\alpha\B_x$ fixes the complement of $\{1,2,3,4\}$. Take a disk
neighborhood $U$ of $x\cup\alpha$ that contains $\alpha\B_x$ and
misses all branch values except the endpoints of $x$ (clearly such
a neighborhood exists). The restriction of $p$ in U is obviously
equivalent to one of the following three coverings
\begin{itemize}
\item[(a)] $(12),(12)$
\item[(b)] $(12),(23)$
\item[(c)] $(12),(24)$
\end{itemize}
via an equivalence that sends $x$ to $x_0$ and $\alpha$,
$\alpha\B_x$ to $\alpha_0$, $\alpha_1$. The proof then reduces to
verifying the statements for those three coverings. For this see
Figure~\ref{fig:intlift}, which gives explicit models of these
coverings constructed by cutting and pasting. The intervals and
their lifts are shown in yellow while the cuts are shown in grey.
It is to be understood that the ``free'' grey arcs are ``folded in
half''. Using this picture one can easily check by projecting that
the maps lift as claimed.
\end{proof}

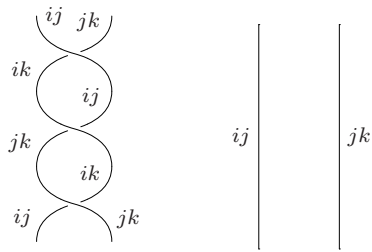
\begin{figure}[htp]
\[
\begin{xy}/r1cm/:
\vtwistneg|(0){ij}<{jk}>{ik}\vcross<{jk}>{ij}\vtwistneg|{jk}>{ij}<{ik}
\end{xy}
\hspace{1cm} \xymatrix{ \ar@{-}[ddd]_{ij}&\ar@{-}[ddd]^{jk}\\
&\\
&\\
&  }
\]
\caption{Move $\mathcal{M}$}\label{moveM}
\end{figure}

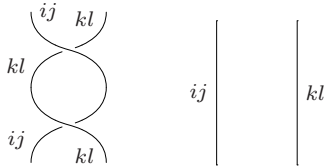
\begin{figure}[htp]
\[
\begin{xy}/r1cm/:
\vtwistneg|(0){ij}<{kl}>{kl}\vcross<{ij}|{kl}
\end{xy}
\hspace{1cm} \xymatrix{ \ar@{-}[dd]_{ij}&\ar@{-}[dd]^{kl}\\
&\\
&  }
\]
\caption{Move $\mathcal{P}$}\label{moveP}
\end{figure}

  Consider the moves between colored braids shown in Figures~\ref{moveM}
  and \ref{moveP}. It should be emphasized that in this section they are
  considered as moves between braid diagrams, that is the strands
  are always {\em vertical}.\

 Let $N(\rho)$ \label{N} be the set of colored braids that can be related to
the identity automorphism of $\rho$, that is the identity braid
colored by the colors of $\rho$, via a finite sequence of these
two moves. Since both moves are local it is clear that $N(\rho)$
is a normal subgroup of $\Au(\rho)$ and furthermore it is natural,
that is if $\rho\B=\rho'$ then $N(\rho')=\B^{-1}N(\rho)\B$.
\begin{defn}
$$\overline{\Au}(\rho)=\frac{\Au(\rho)}{N(\rho)}$$
\end{defn}

 By Lemma~\ref{intlift}, $N\vartriangleleft Ker(\lambda(\rho))$ and
therefore the lifting homomorphism $\lambda(\rho)$ descends to a
homomorphism $\xymatrix@1{
{\bar{\lambda}(\rho):\overline{\Au}(\rho)}\ar[r] & {\CM(E(\rho))}
}$, that is, there is a commutative diagram:
\begin{equation}\label{lambda}
 \xymatrix{
 {\Au(\rho)} \ar[d] \ar[r]^{\lambda(\rho)} & {\CM(E(\rho))}\\
 {\overline{\Au}(\rho)} \ar[ur]_{\bar{\lambda}(\rho)}   }
\end{equation}

\subsection{Coverings over $S^3$}\label{3sphere} Consider an $m$-sheeted covering of $S^3$.
The branching locus will be a codimension-$2$ submanifold of $S^3$
that is a link $L$. According to Proposition~\ref{equiv} the
covering is determined by its monodromy. The total space of the
covering is connected iff the image of the monodromy homomorphism
is a transitive subgroup of $\CS_m$. From now on unless explicitly
mentioned otherwise {\em the total space of a covering will always
be assumed connected}. Recall that coverings are assumed to be
simple which in this context means that the monodromy associated
to each Wirtinger generator is a transposition.

\begin{defn} An $m$-colored link is a link $L$ together with a
simple, transitive representation of its group $\pi_1(S^3\setminus
L)$ into the symmetric group $\CS_m$. A $3$-colored link will be
called tricolored and a $4$-colored link will be called
bi-tricolored.
\end{defn}

Given an isotopy between two links there is an induced isomorphism
between their groups (see Appendix~\ref{diagrams}). Therefore an
isotopy between two links induces a bijection between their
colorings.
\begin{defn}Two colored links are called {\em color-isotopic} (or when
there is no danger of confusion just isotopic) if there exists an
isotopy between them that transforms the coloring of one link to
the coloring of the other link.
\end{defn}

Recall (see Appendix~\ref{diagrams}) that given a diagram of a
link there is an associated Wirtinger presentation of the link
group with a generator for each arc of the diagram and a relation
for each crossing. Therefore an $m$-colored link can be
represented combinatorially by an $m$-colored diagram that is a
link diagram with its arcs labelled with transpositions of $\CS_m$
in a manner compatible with the Wirtinger presentation.
Specifically:

\begin{defn} An $m$-colored link diagram is a link diagram
with its arcs labelled by transpositions of $\CS_m$ in such a way
that the subgroup of $\CS_m$ generated by all the labels is
transitive, and furthermore at each crossing the transposition
labelling one of the under arcs equals the transposition of the
other under arc conjugated by the transposition of the over arc.
\end{defn}

A colored link diagram obviously represents  a colored link. Call
two colored link diagrams equivalent if they represent isotopic
colored links.

\begin{rem}By looking at how Reidemeister moves transform the
presentations of the link group (see Remark~\ref{precoliso}) one
can define colored Reidemeister moves and prove that two colored
link diagrams are equivalent iff they can be related by a finite
sequence of colored Reidemeister moves.
\end{rem}

\begin{defn}An $m$-colored link presentation of a $3$-manifold $M$
is an $m$-colored link diagram $L$, such that $M$ is homeomorphic
to $E(\rho)$. In this case one writes $M=M(L)$ and says that $M$
is represented by $L$.
\end{defn}

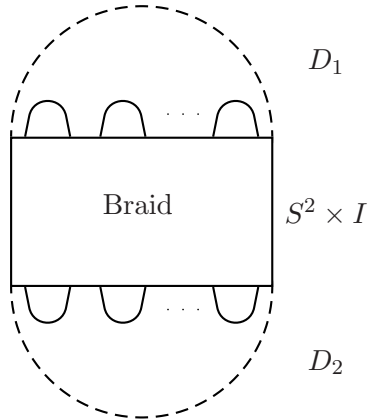
\begin{figure}[htp]
  \begin{pspicture}(5,5)
  \pscircle[linestyle=dashed](2.05,3.5){1.75}
  \pscircle[linestyle=dashed](2.05,1.5){1.75}
   \psframe[linecolor=white,fillstyle=solid,fillcolor=white](.3,3.5)(3.8,1.5)
 \psline[linearc=.25](.5,3.5)(.8,5)(1.1,3.5)
 \psline[linearc=.25](1.5,3.5)(1.8,5)(2.1,3.5)
 \psline[linearc=.25](3,3.5)(3.3,5)(3.6,3.5)
 \psdots[dotsize=.05](2.4,3.8)(2.6,3.8)(2.8,3.8)
 \psframe(.3,3.5)(3.8,1.5)
 \rput(2,2.6){Braid}
  \psline[linearc=.25](.5,1.5)(.8,0)(1.1,1.5)
 \psline[linearc=.25](1.5,1.5)(1.8,0)(2.1,1.5)
 \psline[linearc=.25](3,1.5)(3.3,0)(3.6,1.5)
 \psdots[dotsize=.05](2.4,1.2)(2.6,1.2)(2.8,1.2)
 \rput(4.5,4.5){$D_1$}
 \rput(4.5,.5){$D_2$}
 \rput(4.5,2.5){$S^2\times I$}
  \end{pspicture}
\caption{A plat diagram and the corresponding splitting of
$S^3$}\label{plat}
\end{figure}

\noindent {\bf Plat diagrams and Heegaard splittings.} One way to
understand the manifold represented by a colored link diagram is
the following:
 Start with a plat diagram of the colored link (see for example
 \cite{B} for a proof of the fact that each link has a plat
 diagram),
  that is a
braid diagram closed on its top and bottom by caps and cups and
split the three-sphere $S^3$ as
$$S^3=D_1\bigcup S^2\times I\bigcup D_2$$
where $D_1$ is a $3$-dimensional disc containing all the caps,
$D_2$ is a $3$-dimensional disc containing all the cups, and the
braid is contained in $S^2\times I$ as shown in Figure~\ref{plat}.

Now $D_1$ (respectively $D_2$) is covered by a handlebody $H_1$
(respectively $H_2$) whose genus depends only on the number of
caps (respectively cups), see~\cite{BE}. Since there are as many
caps as cups the two handlebodies have the same genus. On the
other hand $S^2 \times I$ is covered by a ``thick'' surface
``realizing'' a homeomorphism $\xymatrix@1{{f:\partial
H_1}\ar[r]&{\partial H_2}}$. Therefore one gets a Heegaard
splitting of the total space of the covering represented by the
colored link.\

\begin{defn} A \textit{normalized diagram} is a colored link
diagram in a plat form such that its top and bottom are equal to
\[
\rho_{\text{stand}}:=\begin{array}{lr}
\begin{xy}/r.9cm/:
\vcap|{12}
\end{xy}
\hspace{5mm}
\begin{xy}/r.9cm/:
\vcap|{14}
\end{xy}
\hspace{2.5mm} \cdots \hspace{2.5mm}
\begin{xy}/r.9cm/:
\vcap|{1m}
\end{xy}
\hspace{5mm}
\begin{xy}/r.9cm/:
\vcap|{23}
\end{xy}
\hspace{2.5mm} \cdots \hspace{2.5mm}
\begin{xy}/r.9cm/:
\vcap|{23}
\end{xy}
&\\
& \\
\end{array}
\]
\end{defn}

\begin{rem} A variety of different normalizations can be found in
the literature. The one chosen above is most suitable for the
purposes of this work.
\end{rem}

\begin{prop}\label{norm}Every colored link has a normalized diagram.
\end{prop}

The following combinatorial lemma will be used in the proof.
\begin{lem}[\bf Bernstein-Edmmonds]\label{Previous} Let
$\rho=\s_0$, $\s_1$,$\dotsc,\s_{2k-2}$, $\s_{2k}$ be a
sequence of transpositions of $\CS_m$ satisfying the following
four conditions:
 \begin{itemize}
 \item[(i)] $\s_0\s_1\dotsm\s_{2k-1}=\text{id}$
 \item[(ii)] each $\s_i$ is a transposition
 \item[(iii)] the subgroup of $\CS_m$ generated by the $\s_i$'s is transitive
 \item[(iv)] $\s_{2j}=\s_{2j+1}$ for $j=0,\dotsc,k-1$.
 \end{itemize}
Then $\rho$ can be related to any other such sequence via a finite
sequence of the following two moves:
\begin{itemize}
 \item[(a)] interchange the position of any two pairs
 \item[(b)] replace a pair of transpositions by its conjugate by a neighboring pair.
  \end{itemize}
\end{lem}
\begin{proof}Proposition 3.8 of \cite{BE}
\end{proof}

\begin{proof}[Proof of~\ref{norm}] It suffices to show that moves
(a) and (b) of Lemma~\ref{Previous} can be realized by colored
isotopy. This is done in Figure~\ref{NORM}
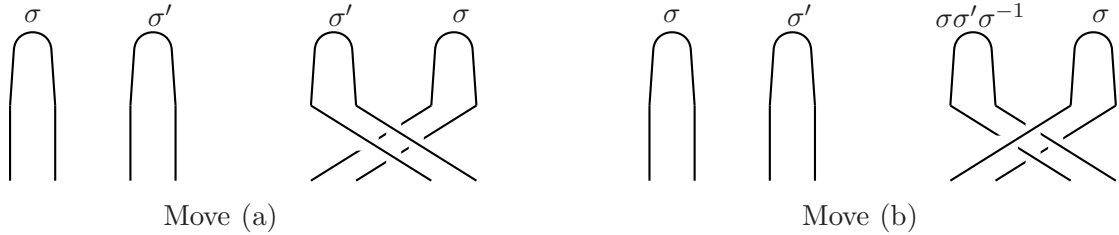
\begin{figure}[htp]
\[
\begin{pspicture}(-5,-1)(1.5,3)
 \psline(-4.8,1)(-4.8,0)  \psline(-4.2,1)(-4.2,0)
 \psline(-3.2,1)(-3.2,0)  \psline(-2.6,1)(-2.6,0)
 \psline[linearc=.25](-4.8,1)(-4.5,5.4)(-4.2,1)
 \psline[linearc=.25](-3.2,1)(-2.9,5.4)(-2.6,1)
 \psline(-.8,0)(.8,1)  \psline(-.2,0)(1.4,1)
 \psframe[linecolor=white,fillstyle=solid,fillcolor=white](-.9,.4)(.1,.6)
 \psframe[linecolor=white,fillstyle=solid,fillcolor=white](.2,.2)(.4,.4)
 \psframe[linecolor=white,fillstyle=solid,fillcolor=white](.2,.6)(.4,.8)
 \psframe[linecolor=white,fillstyle=solid,fillcolor=white](.5,.4)(.7,.6)
 \psline(1.4,0)(-.2,1)  \psline(-.8,1)(.8,0)
  \psline[linearc=.25](-.8,1)(-.5,5.4)(-.2,1)
 \psline[linearc=.25](.8,1)(1.1,5.4)(1.4,1)
   \rput(-2,-.5){Move (a)}
   \rput(-4.5,2.2){$\s$}
   \rput(-2.8,2.2){$\s'$}
   \rput(-.4,2.2){$\s'$}
   \rput(1.2,2.2){$\s$}
\end{pspicture}
\hspace{2cm}
\begin{pspicture}(-5,-1)(1.5,3)
 \psline(-4.8,1)(-4.8,0)  \psline(-4.2,1)(-4.2,0)
 \psline(-3.2,1)(-3.2,0)  \psline(-2.6,1)(-2.6,0)
 \psline[linearc=.25](-4.8,1)(-4.5,5.4)(-4.2,1)
 \psline[linearc=.25](-3.2,1)(-2.9,5.4)(-2.6,1)
   \psline(-.2,0)(1.4,1)
  \psframe[linecolor=white,fillstyle=solid,fillcolor=white](.2,.2)(.4,.4)
  \psframe[linecolor=white,fillstyle=solid,fillcolor=white](.5,.4)(.7,.6)
 \psline(1.4,0)(-.2,1)  \psline(-.8,1)(.8,0)
 \psframe[linecolor=white,fillstyle=solid,fillcolor=white](-.9,.4)(.1,.6)
 \psframe[linecolor=white,fillstyle=solid,fillcolor=white](.2,.6)(.4,.8)
  \psline(-.8,0)(.8,1)
  \psline[linearc=.25](-.8,1)(-.5,5.4)(-.2,1)
 \psline[linearc=.25](.8,1)(1.1,5.4)(1.4,1)
   \rput(-2,-.5){Move (b)}
   \rput(-4.5,2.2){$\s$}
   \rput(-2.8,2.2){$\s'$}
   \rput(-.4,2.2){$\s\s'\s^{-1}$}
   \rput(1.2,2.2){$\s$}
\end{pspicture}
\]
 \caption{Realizing moves (a) and (b)}\label{NORM}
\end{figure}
\end{proof}

\subsection{Dimming the lights}\label{dim} 
\begin{figure}[htp]
\definecolor{lightblue}{rgb}{.5, .65, .9}
\definecolor{lightgreen}{rgb}{.5, .9, .7}
\definecolor{lightred}{rgb}{.9, .5, .5}
\begin{pspicture}(0,0)(5,3.5)
\rput(.85,0){1}
 \rput(2.5,3.2){3}
 \rput(4.15,0){2}
 \rput(2.57,1.1){4}
 \pnode(1,0){A}
 \pnode(2.5,3){B}
 \pnode(4,0){C}
 \pnode(2.75,.9){D}
 \ncline[linecolor=lightgreen]{A}{B}
 \ncline[linecolor=lightred]{A}{C}
 \ncline[linecolor=lightblue]{B}{C}
 \ncline[linestyle=dotted,linecolor=blue]{A}{D}
 \ncline[linestyle=dotted,linecolor=red]{B}{D}
 \ncline[linestyle=dotted,linecolor=green]{C}{D}
 \end{pspicture}
\caption{How to see the map $\kappa$}\label{tetra}
\end{figure}
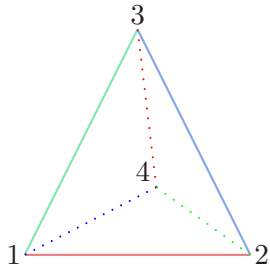

There is an exact sequence of groups:
$$\xymatrix@1{ {1}\ar[r] & {V} \ar[r] & {\CS_4}
\ar[r]^{\kappa} & {\CS_3} \ar[r] & {1} }$$
 where $V=\{\text{id},(12)(34),(14)(23),(13)(24)\}$ is  Klein's
Vierergruppe. For the purposes of this work $\kappa$ is best seen
by identifying the edges of a numbered tetrahedron with the
transpositions of $\CS_4$ and the edges of its front face with the
transpositions of $\CS_3$ (see Figure~\ref{tetra}). Then $\kappa$
fixes the front edges and sends each back edge to its opposite
edge.

 It is customary to identify, see~\cite{P1} for example, the three
transpositions of $\CS_3$ with three colors Red, Green and Blue
via $R=(12)$, $B=(23)$ and $G=(13)$. The above exact sequence then
suggests the identification of the transpositions of $\CS_4$ with
three colors that come in light and dark shades. Denoting by
$\tilde X$ the dark shade of a color $X$ this identification is
given explicitly by:
$$\begin{array}{lll}
R=(12) & B=(23) & G=(13)\\
\tilde R=(24) & \tilde B=(14) &\tilde  G=(24)
\end{array}$$
The homomorphism $\kappa$ can then be described as ``dimming the
lights'' so that the two shades of the same color cannot be
distinguished.\

 There is an induced functor from the category of
``bi-tricolored'' braids to the category of tricolored braids
which will also be referred to as dimming the lights. Dimming the
lights obviously commutes with the lifting functor.

\subsection{Adding a trivial sheet}\label{trivialsheet} Let $\xymatrix@1{
{p:M}\ar[r]&{S^n}}$ be an $m$-sheeted covering over the $n$-sphere
branched over the codimension-$2$ submanifold $L$. One can
construct an $m+1$-sheeted covering of $S^n$ with the same total
space and branched over $L\sqcup K$, where $K$ is an unknotted
$n-2$-sphere contained in a ball which is disjoint from $L$, as
follows:\

Let $D$ be an $(n-1)$-disk disjoint from $L$ with $\partial D=K$.
Since $D$ is simply connected, $p^{-1}(D)$ consists of $m$
disjoint $(n-1)$-disks each mapping homeomorphically onto $D$.
Chose one of these disks, say $\widetilde D$, and cut $M$ open
along $\widetilde D$, to get a manifold $M_{\text{cut}}$  together
with a projection
$\xymatrix@1{{p_{\text{cut}}:M_{\text{cut}}}\ar[r]&{S^n}}$. The
boundary of $M_{\text{cut}}$ is obtained by gluing two copies
$\widetilde D_1$ and $\widetilde D_2$ of $\widetilde D$ along
their boundary and is therefore homeomorphic to $S^{n-1}$.\\
Also cut a copy of $S^{n}$, the $(m+1)^{\text{th}}$ sheet, open
along $K$ to get an $n$-disk $S_{\text{cut}}$ with its boundary
obtained by gluing two copies $D_1$ and $D_2$ of $D$ along their
boundary.\\
Now glue $M_{\text{cut}}$ and $S_{\text{cut}}$ along their
boundary in such a way that $\widetilde D_1$ is glued to $D_2$ and
$\widetilde D_2$ is glued to $D_1$, to get a connected sum of $M$
with $S^n$ which is homeomorphic to $M$. By gluing
$p_{\text{cut}}$ and the gluing map $\xymatrix@1{
S_{\text{cut}}\ar[r]&S^n}$ one gets a projection
$$\xymatrix@1{{p':M}\ar[r]&S^n}.$$

 It is easily checked that $\xymatrix@1{{p':M}\ar[r]&S^n}$ is an
 $(m+1)$-sheeted covering branched over $L\sqcup K$. $p'$ is
 said to be obtained by $p$ by {\em adding a trivial sheet}.
$\pi_1(S^n\setminus L\sqcup K)$ is a free product of
$\pi_1(S^n\setminus L)$
 and an infinite cyclic group
generated by a meridian around $K$ and it is clear that the
monodromy of $p'$ coincides with the monodromy of $p$ in the first
factor and associates a transposition of
 the form $(i,m+1)$ to the meridian around K. Therefore $p'$
is simple whenever $p$ is.

\section{The $3$-sheeted case}\label{3sheets}
In this section we give some results when the degree of the
coverings is $3$. These results will be used either
directly or by analogy in the study of coverings of degree $4$.

\subsection{$3$-sheeted Coverings of $S^2$}\label{2sphere3} Although we are mainly
interested in $3$-sheeted coverings of the $2$-sphere $S^2$ it is
convenient to consider more generally coverings of the $2$-disc
$D^2$. Recall from Section~\ref{2D} how such coverings are
represented as finite sequences of transpositions of $\CS_3$: the
length of the sequence is the number of branch values and the
$i^{\text{th}}$ term of the sequence is the $i^{\text{th}}$
monodromy.\

  Consider the following covering over $D^2$
$$\rho(n):=(12),(12),(23),(23),\dotsc,(23)\quad, $$
where there are $n$ branch values and all non displayed
monodromies are equal to $(23)$. $\rho(n)$ corresponds to a
covering over $S^2$ iff $n$ is even(see Section~\ref{2D} on
page~\pageref{2sphere}). In that case the Riemann-Hurwitz formula
\eqref{genus} gives the following relation between the number of
branch values $n$ and the genus $g$ of the total space:
$$n=2g-4 \quad.$$

Denote by $L(n)$ the subgroup of the braid group $B_n$ that fixes
$\rho(n)$. Recall  from that section~\ref{2D} on
page~\pageref{ColBraid}, that $L(n)$ (denoted $\Au(\rho(n))$
there) coincides with the subgroup of the mapping class group
$\CM_{0,1,n}$ that consists of liftable (isotopy classes of)
homeomorphisms.\

This covering was studied by Birman and Wajnryb in~\cite{BW1}. The
remaining of this section is an exposition of the results of that
paper that are relevant to this study.

\begin{thm}[\bf Birman-Wajnryb]\label{bw1}$L(n)$ is the smallest
subgroup of $B_n$ containing the following elements:

$$\B_0, \B_1^3, \B_2,\B_3,\dotsc,\B_{n-2},\quad \text{and if } n\geq
6,\quad \delta_4$$
 where $\delta_4=[\B_4]\B_3\B_2\B_1^2\B_2\B_3^2\B_2\B_1$ is the
 rotation around the interval $d_4$ shown in Figure~\ref{D4}.

\begin{figure}[htp]
\begin{pspicture}(0,-.9)(10,1.1)
 \psdots[dotstyle=*,dotscale=1](0,0)(1,0)(2,0)(3,0)(4,0)(5,0)(6,0)(7.6,0)
 \pscurve (1,0)(2.7,.6)(4.5,0)(2.7,-.86)(.4,0)(2,.98)(3,1)(4,.8)(5,0)
 \psdots[dotstyle=*,dotscale=.4](6.5,0)(6.8,0)(7.1,0)
 \rput(0,-.26){0}
 \rput(1,-.26){1}
 \rput(2,-.26){2}
 \rput(3,-.26){3}
 \rput(4,-.26){4}
 \rput(5,-.26){5}
 \rput(6,-.26){6}
 \rput(7.6,-.26){$n-1$}
\end{pspicture}
\caption{The interval $d_4$}\label{D4}
\end{figure}
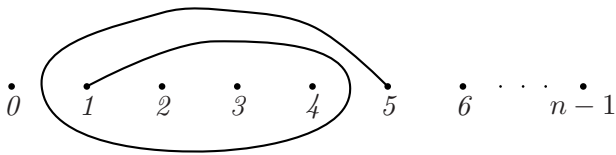
\end{thm}
\begin{proof} Theorem 6.1 of~\cite{BW1}. See also Theorem 3.1 of
the same paper. Note that the same symbol is used there to denote
an interval and the rotation around it.
\end{proof}

 From now on the generators given in the above theorem will be
 referred to as {\em the} generators of $L(n)$. Observe that each
 generator is a power of a rotation around some interval.
 Therefore according to Lemma~\ref{intlift} their lifts can be determined
 by determining how these intervals lift. One gets the following
 theorem, where the notation used for elements of the mapping
 class group $\CM(E(\rho(n))$ is that of~\cite{W} (see Appendix~\ref{Pres}).

\begin{figure}[htp]
\begin{pspicture}(-1,-1)(14,14)

 \psellipse[linecolor=blue](.5,1.5)(.5,.3)
  \psellipse[linecolor=blue](2.5,1.5)(.5,.3)
 \psellipse[linecolor=blue](4.5,1.5)(.5,.3)
 \psellipse[linecolor=blue](6.5,1.5)(.5,.3)
 \psellipse[linecolor=blue](8.5,1.5)(.5,.3)
 \psellipse[linecolor=blue](12.5,1.5)(.5,.3)
 \psbezier(1,1.5)(1.25,1)(1.75,1)(2,1.5)
  \psbezier(3,1.5)(3.25,1)(3.75,1)(4,1.5)
   \psbezier(5,1.5)(5.25,1)(5.75,1)(6,1.5)
    \psbezier(7,1.5)(7.25,1)(7.75,1)(8,1.5)
     \psbezier(9,1.5)(9.25,1)(9.75,1)(10,1.5)
          \psbezier(11,1.5)(11.25,1)(11.75,1)(12,1.5)
     \psbezier
     (0,1.5)(-.5,-1)(13.5,-1)(13,1.5)
      \psbezier[linestyle=dotted,linecolor=red](5,1.5)(5.2,0)(.8,-.8)(.5,1.8)
  \psbezier[linecolor=red](.4,1.2)(0,0)(4,.6)(4.5,1.2)
  \psbezier[linestyle=dotted,linecolor=red](4.5,1.8)(4,.6)(.8,.4)(1,1.5)

\rput(0,7){
 \psellipse[linecolor=blue](.5,1.5)(.5,.3)
 \psframe[linecolor=white,fillstyle=solid,fillcolor=white](-.5,1.5)(1.5,2)
 \psellipse[linestyle=dotted,linecolor=blue](.5,1.5)(.5,.3)
 \psellipse[linecolor=blue](2.5,1.5)(.5,.3)
 \psellipse[linecolor=blue](4.5,1.5)(.5,.3)
 \psellipse[linecolor=blue](6.5,1.5)(.5,.3)
 \psellipse[linecolor=blue](8.5,1.5)(.5,.3)
 \psellipse[linecolor=blue](12.5,1.5)(.5,.3)
 \psbezier(1,1.5)(1.25,1)(1.75,1)(2,1.5)
  \psbezier(3,1.5)(3.25,1)(3.75,1)(4,1.5)
   \psbezier(5,1.5)(5.25,1)(5.75,1)(6,1.5)
    \psbezier(7,1.5)(7.25,1)(7.75,1)(8,1.5)
     \psbezier(9,1.5)(9.25,1)(9.75,1)(10,1.5)
          \psbezier(11,1.5)(11.25,1)(11.75,1)(12,1.5)
     \psbezier
     (0,1.5)(-.5,-1)(13.5,-1)(13,1.5)
  \psframe[linecolor=white,fillstyle=solid,fillcolor=white](9.5,-1)(11.5,2)
  \psdots[dotsize=.1](10.3333,.5)(10.6666,.5)(10.9999,.5)
  \psbezier[linestyle=dotted,linecolor=red](5,1.5)(5.2,0)(.8,-.8)(.5,1.8)
  \psbezier[linecolor=red](.4,1.2)(0,0)(4,.6)(4.5,1.2)
  \psbezier[linestyle=dotted,linecolor=red](4.5,1.8)(4,.6)(.8,.4)(1,1.5)}

   \psbezier
     (1,8.5)(0,11.5)(12.5,8)(13,11.5)
  \psbezier
     (0,8.5)(-.5,12)(1.5,10)(2,11.5)
   \psellipse[linecolor=blue](2.5,11.5)(.5,.3)
   \psellipse[linecolor=blue](4.5,11.5)(.5,.3)
   \psellipse[linecolor=blue](6.5,11.5)(.5,.3)
   \psellipse[linecolor=blue](8.5,11.5)(.5,.3)
   \psellipse[linecolor=blue](12.5,11.5)(.5,.3)
   \psellipse[linecolor=blue](.5,11.5)(.5,.3)
   \psframe[linecolor=white,fillstyle=solid,fillcolor=white](2,11.5)(13,12.5)
   \psellipse[linestyle=dotted,linecolor=blue](2.5,11.5)(.5,.3)
   \psellipse[linestyle=dotted,linecolor=blue](4.5,11.5)(.5,.3)
   \psellipse[linestyle=dotted,linecolor=blue](6.5,11.5)(.5,.3)
   \psellipse[linestyle=dotted,linecolor=blue](8.5,11.5)(.5,.3)
   \psellipse[linestyle=dotted,linecolor=blue](12.5,11.5)(.5,.3)
   \psbezier(3,11.5)(3.25,11)(3.75,11)(4,11.5)
   \psbezier(5,11.5)(5.25,11)(5.75,11)(6,11.5)
    \psbezier(7,11.5)(7.25,11)(7.75,11)(8,11.5)
     \psbezier(9,11.5)(9.25,11)(9.75,11)(10,11.5)
          \psbezier(11,11.5)(11.25,11)(11.75,11)(12,11.5)
\psbezier[linestyle=dotted,linecolor=red](5,11.5)(5.2,9.6)(0,9.6)(-.05,9.8)
  \psbezier[linestyle=dotted,linecolor=red](4.5,11.8)(4.8,10)(0,10.2)(.2,10.5)
    \psbezier[linecolor=red]
    (1,9)(0,11)(3,9)(4.7,11.2)
 \pscurve[linestyle=dotted,linecolor=red](1,9)(.7,8.9)(.5,8.8)
 \psbezier[linecolor=red]
 (.2,10.5)(-.2,8.8)(.7,9.6)(1,8.5)
 \psbezier[linecolor=red]
 (-.05,9.8)(-.1,8.8)(.5,9.6)(.4,8.2)

     \psbezier(1,11.5)(1.25,12)(1.75,12)(2,11.5)
    \psbezier(3,11.5)(3.25,12)(3.75,12)(4,11.5)
   \psbezier(5,11.5)(5.25,12)(5.75,12)(6,11.5)
    \psbezier(7,11.5)(7.25,12)(7.75,12)(8,11.5)
     \psbezier(9,11.5)(9.25,12)(9.75,12)(10,11.5)
          \psbezier(11,11.5)(11.25,12)(11.75,12)(12,11.5)
 \psbezier(0,11.5)(-.5,14)(13.5,14)(13,11.5)
\psbezier[linecolor=red](5,11.5)(5.2,13.2)(.2,12.5)(.5,11.8)
  \psbezier[linecolor=red](4.7,11.2)(4.7,12.4)(.8,12.4)(1,11.5)
  \psbezier[linestyle=dotted,linecolor=red](4.5,11.8)(5.2,13.6)(.2,12.6)(.5,11.2)

 \psframe[linecolor=white,fillstyle=solid,fillcolor=white](9.5,-1)(11.5,14)
  \psdots[dotsize=.1](10.3333,.5)(10.6666,.5)(10.9999,.5)
  \psdots[dotsize=.1](10.3333,7.5)(10.6666,7.5)(10.9999,7.5)
  \psdots[dotsize=.1](10.3333,10.5)(10.6666,10.5)(10.9999,10.5)
  \psdots[dotsize=.1](10.3333,12.5)(10.6666,12.5)(10.9999,12.5)

  \rput(-1,7.5){Sheet 1}
 \rput(-1,10){Sheet 2}
 \rput(-1,12.5){Sheet 3}

 \psline{->}(6.5,6)(6.5,2.5)
\end{pspicture}
\caption{The covering $\rho(n)$ and how $d_4$ lifts}\label{3malako}
\end{figure}
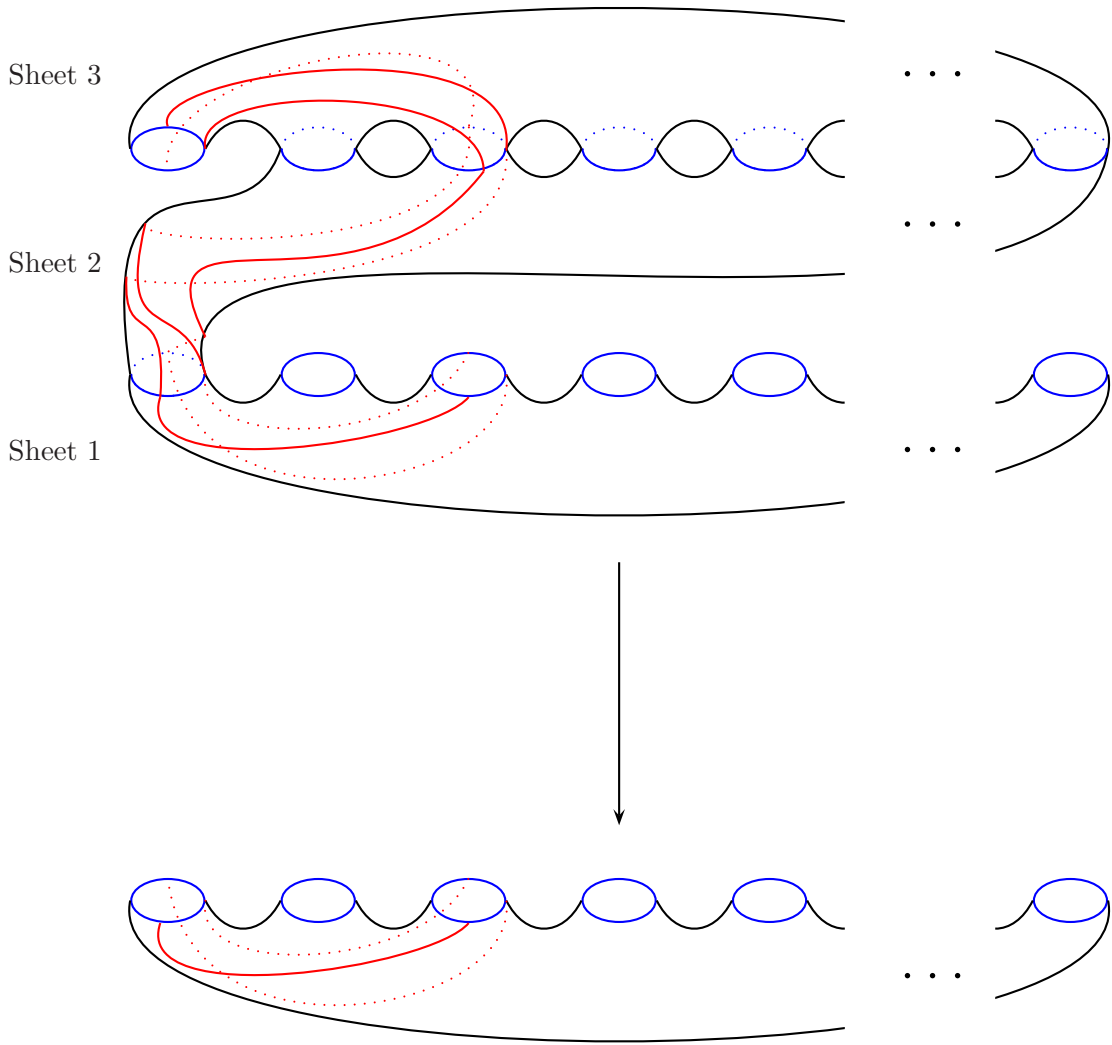

 \begin{thm}[\bf Birman-Wajnryb]\label{bw2} The lifting homomorphism
 $\xymatrix@1{ {\lambda:L(n)}
 \ar[r]&{\CM(E(\rho(n))} }$ is given on the generators of $L(n)$
 by
 \begin{equation}\label{the3map}
\lambda(x)=
\begin{cases} \text{id} & \text{if} \quad x=\B_0\text{ or } \B_1^3\\
                       a_{i} & \text{if $ \quad x=\B_{2i}$ with $i\geq1$}\\
                       b_{i} &\text{if  $\quad x=\B_{2i+1}$ with
                       $i\geq1$}\\
                       d & \text{if $ \quad x=\delta_4$}
                   \end{cases}
  \end{equation}
 \end{thm}
\begin{proof} Figure~\ref{3malako} gives an explicit model of $\rho(n)$.
This model is constructed by cutting and pasting as described in
Section~\ref{2D}, page~\pageref{cutpaste}. The $2$-sphere is cut
open along the intervals $x_0$, $x_2,\dotsc,x_{2g-2}$. Since the
monodromy along a loop that goes around one of those intervals is
id the resulting $g-1$-holed sphere lifts to three disjoint copies
of itself. The cuts are shown in blue and its assumed in
Figure~\ref{3malako} that free blue loops are ``sewn'' to become
intervals so that we get a closed surface covering $S^2$. \
 $d_4$ and how it lifts to an arc and a loop isotopic to $d$ ia
 shown in red. The check for the other generators is rather
 trivial.
\end{proof}

As an immediate corollary one gets the following:
\begin{thm}\label{3surj} For $3$-sheeted coverings of the $2$-sphere the
lifting homomorphism is surjective.
\end{thm}
\begin{proof} All Wajnryb generators are in the image of the
lifting homomorphism.
\end{proof}

Furthermore using Wajnryb's presentation the authors of~\cite{BW1}
gave a set of normal generators of the kernel of the lifting
homomorphism:
\begin{thm}[\bf Birman-Wajnryb]\label{Ker3} The kernel of the lifting homomorphism\\
 $\xymatrix@1{ {\lambda:L(n)}
 \ar[r]&{\CM(E(\rho(n))} }$ is the smallest normal subgroup of
 $L(n)$ containing the following elements (recall that $[\B']\B$ for braids $\B'$, $\B$
 means $\B^{-1}\B'\B$):\\
$\B_0$, $\B_1^3$, $B$ and $D$, where
$$B=(\B_2\B_3\B_4)^4([\delta_4^{-1}]\B_4^{-1}\B_3^{-1}
\B_2^{-2}\B_3^{-1}\B_4^{-1}\B_5^{-1})\delta_4^{-1},$$
$$D=\B_{2g+2}\chi\B_{2g+2}^{-1}\chi^{-1},$$
Where:
$$\chi=\B_{2g+1}\B_{2g}\dotsm\B_3\B_2^2\B_3\dotsm\B_{2g}\B_{2g+1},$$
 \end{thm}
 \begin{proof} See the proof of Theorems 5.1 and 6.1
 of~\cite{BW1}. See also the Errata~\cite{BW2}.
 \end{proof}
 \begin{rem} The normal generators  in~\cite{BW1} are slightly
 different than the ones given above. In particular  the element
  we denote by $B$ is denoted by $B'$ there while $B$ denotes an
  element which is equal to $B'$ modulo the other generators. The
  generators we use are the same as those in~\cite{P1}.
\end{rem}

\subsection{$3$-sheeted Coverings of $S^3$}\label{3sphere3}

Recall from Section~\ref{3sphere} that coverings over $S^3$ are
represented by colored link diagrams. In particular for
$3$-sheeted coverings we have {\em tricolored} link diagrams that
is link diagrams whose arcs are colored with three colors
(corresponding to the three transpositions of $\CS_3$) in such a
way that
\begin{itemize}
\item[a)] at least two colors are used (to ensure that the total
space of the covering is connected), and
\item[b)] at each crossing there are either three colors
or only one (to ensure that the Wirtinger relations are
preserved).
\end{itemize}

The following theorem was proved independently by Montesinos
 in~\cite{M1} and Hilden in~\cite{H}. The proof we provide is essentially that
 of Hilden.

\begin{thm}[\bf Hilden, Montesinos]\label{trinity} Every orientable closed $3$-manifold
is a simple $3$-sheeted covering of the $3$-sphere.
\end{thm}
\begin{proof} Since by Theorem~\ref{3surj} the lifting homomorphism is surjective all
Heegaard splittings can be realized (see Section~\ref{3sphere}) .
\end{proof}

Thus any orientable closed $3$-manifold can be represented by a
tricolored link diagram. Of course the same manifold manifests
itself as a $3$-sheeted covering of $S^3$ in many different ways
and therefore admits several presentations as a tricolored link.
Naturally then one asks the question whether there is a finite set
of moves such that two tricolored link diagrams represent the same
$3$-manifold iff they can be related via a finite sequence of
moves from this set.\
 For some time it was conjectured that the ``Montesinos move''
 $\mathcal{M}$ shown in Figure~\ref{moveM} is the answer. However
 Montesinos gave a counterexample in \cite{M2}. Then in the $90$'s
 Piergallini proved the following theorem in \cite{P1}.

\begin{figure}[htp]
\[
\begin{pspicture}(-.5,0)(5,5)
 \psdots(0,4)(.5,4)(1,4)(1.5,4)(2,4)(2.5,4)(4,4)(4.5,4)
        (0,1)(.5,1)(1,1)(1.5,1)(2,1)(2.5,1)(4,1)(4.5,1)
   \psline(1,4)(1,1)
   \psline(1.5,4)(1.5,1)
   \psline(2,4)(2,1)
   \psline(2.5,4)(2.5,1)
   \psline(4,4)(4,1)
   \psline(4.5,4)(4.5,1)
 \psline(0,4)(0,2.75)
 \psline(0,2.25)(0,1)
 \psline(0.5,4)(0.5,2.75)
 \psline(0.5,2.25)(0.5,1)
  \psline(0,2.75)(.5,2.25)
  \psline(.5,2.75)(.3,2.55)
  \psline(.2,2.44)(0,2.25)
           \psframe(-.3,0)(4.8,1)
           \psframe(-.3,4)(4.8,5)
    \rput(0,4.3){12}
 \rput(0.5,4.3){12}
 \rput(1,4.3){23}
 \rput(1.5,4.3){23}
 \rput(2,4.3){23}
 \rput(2.5,4.3){23}
 \rput(4,4.3){23}
 \rput(4.5,4.3){23}
                   \rput(3.2,4.3){23's}
 \rput(0,.7){12}
 \rput(0.5,.7){12}
 \rput(1,.7){23}
 \rput(1.5,.7){23}
 \rput(2,.7){23}
 \rput(2.5,.7){23}
 \rput(4,.7){23}
 \rput(4.5,.7){23}
                   \rput(3.2,.7){23's}
  \psdots[dotsize=0.1](3,4)(3.3,4)(3.6,4)(3,1)(3.3,1)(3.6,1)
                      (3,2.5)(3.3,2.5)(3.6,2.5)
     \rput(2.5,4.8){$\tau$}
     \rput(2.5,.2){$\tau'$}
\end{pspicture}
\hspace{2cm}
\begin{pspicture}(-.5,0)(5,5)
 \psdots(0,4)(.5,4)(1,4)(1.5,4)(2,4)(2.5,4)(4,4)(4.5,4)
        (0,1)(.5,1)(1,1)(1.5,1)(2,1)(2.5,1)(4,1)(4.5,1)
   \psline(1,4)(1,1)
   \psline(1.5,4)(1.5,1)
   \psline(2,4)(2,1)
   \psline(2.5,4)(2.5,1)
   \psline(4,4)(4,1)
   \psline(4.5,4)(4.5,1)
 \psline(0,4)(0,1)
 \psline(0.5,4)(0.5,1)
             \psframe(-.3,0)(4.8,1)
             \psframe(-.3,4)(4.8,5)
 \rput(0,4.3){12}
 \rput(0.5,4.3){12}
 \rput(1,4.3){23}
 \rput(1.5,4.3){23}
 \rput(2,4.3){23}
 \rput(2.5,4.3){23}
 \rput(4,4.3){23}
 \rput(4.5,4.3){23}
                   \rput(3.2,4.3){23's}
 \rput(0,.7){12}
 \rput(0.5,.7){12}
 \rput(1,.7){23}
 \rput(1.5,.7){23}
 \rput(2,.7){23}
 \rput(2.5,.7){23}
 \rput(4,.7){23}
 \rput(4.5,.7){23}
                   \rput(3.2,.7){23's}
  \psdots[dotsize=0.1](3,4)(3.3,4)(3.6,4)(3,1)(3.3,1)(3.6,1)
                      (3,2.5)(3.3,2.5)(3.6,2.5)
     \rput(2.5,4.8){$\tau$}
     \rput(2.5,.2){$\tau'$}
\end{pspicture}
\]
 \caption{Move $P_{II}$}\label{3II}
\end{figure}
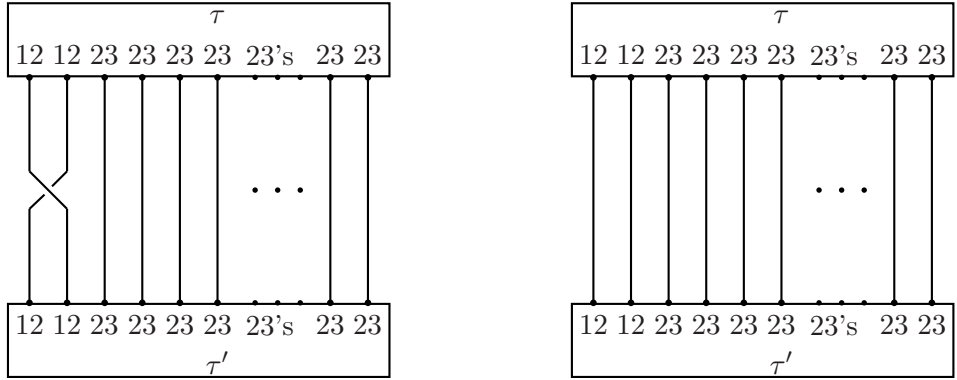

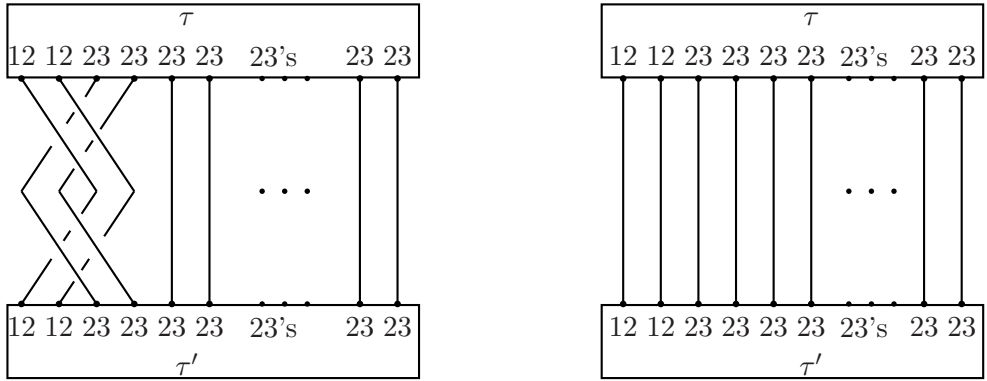
\begin{figure}[htp]
\[
\begin{pspicture}(.5,0)(6.5,5)
 \psdots(1,4)(1.5,4)(2,4)(2.5,4)(3,4)(3.5,4)(5.5,4)(6,4)
        (1,1)(1.5,1)(2,1)(2.5,1)(3,1)(3.5,1)(5.5,1)(6,1)
     \psline(3,4)(3,1)
   \psline(3.5,4)(3.5,1)
   \psline(5.5,4)(5.5,1)
   \psline(6,4)(6,1)
  \psline(1,4)(2,2.5)
  \psline(1.5,4)(2.5,2.5)
  \psline(1,2.5)(2,1)
  \psline(1.5,2.5)(2.5,1)
   \psline(1,2.5)(1.4,3.1)
   \psline(1.6,3.4)(1.7,3.55)
   \psline(1.85,3.755)(2,4)
    \psline(1.5,2.5)(1.7,2.8)
    \psline(1.85,3.025)(1.95,3.175)
    \psline(2.05,3.325)(2.5,4)
    \psline(2,2.5)(1.8,2.2)
    \psline(1.65,1.975)(1.56,1.84)
    \psline(1.4,1.6)(1,1)
    \psline(2.5,2.5)(2.1,1.9)
    \psline(1.9,1.6)(1.8,1.45)
    \psline(1.7,1.3)(1.5,1)
           \psframe(.8,0)(6.3,1)
           \psframe(.8,4)(6.3,5)
  \rput(1,4.3){12}
 \rput(1.5,4.3){12}
 \rput(2,4.3){23}
 \rput(2.5,4.3){23}
 \rput(3,4.3){23}
 \rput(3.5,4.3){23}
 \rput(5.5,4.3){23}
 \rput(6,4.3){23}
                   \rput(4.35,4.3){23's}
  \rput(1,.7){12}
 \rput(1.5,.7){12}
 \rput(2,.7){23}
 \rput(2.5,.7){23}
 \rput(3,.7){23}
 \rput(3.5,.7){23}
 \rput(5.5,.7){23}
 \rput(6,.7){23}
                   \rput(4.35,.7){23's}
  \psdots[dotsize=0.1](4.2,4)(4.5,4)(4.8,4)(4.2,1)(4.5,1)(4.8,1)
                      (4.2,2.5)(4.5,2.5)(4.8,2.5)
     \rput(3.2,4.8){$\tau$}
     \rput(3.2,.2){$\tau'$}
\end{pspicture}
\hspace{2cm}
\begin{pspicture}(-.5,0)(5,5)
 \psdots(0,4)(.5,4)(1,4)(1.5,4)(2,4)(2.5,4)(4,4)(4.5,4)
        (0,1)(.5,1)(1,1)(1.5,1)(2,1)(2.5,1)(4,1)(4.5,1)
   \psline(1,4)(1,1)
   \psline(1.5,4)(1.5,1)
   \psline(2,4)(2,1)
   \psline(2.5,4)(2.5,1)
   \psline(4,4)(4,1)
   \psline(4.5,4)(4.5,1)
 \psline(0,4)(0,1)
 \psline(0.5,4)(0.5,1)
             \psframe(-.3,0)(4.8,1)
             \psframe(-.3,4)(4.8,5)
 \rput(0,4.3){12}
 \rput(0.5,4.3){12}
 \rput(1,4.3){23}
 \rput(1.5,4.3){23}
 \rput(2,4.3){23}
 \rput(2.5,4.3){23}
 \rput(4,4.3){23}
 \rput(4.5,4.3){23}
                   \rput(3.2,4.3){23's}
 \rput(0,.7){12}
 \rput(0.5,.7){12}
 \rput(1,.7){23}
 \rput(1.5,.7){23}
 \rput(2,.7){23}
 \rput(2.5,.7){23}
 \rput(4,.7){23}
 \rput(4.5,.7){23}
                   \rput(3.2,.7){23's}
  \psdots[dotsize=0.1](3,4)(3.3,4)(3.6,4)(3,1)(3.3,1)(3.6,1)
                      (3,2.5)(3.3,2.5)(3.6,2.5)
     \rput(2.5,4.8){$\tau$}
     \rput(2.5,.2){$\tau'$}
\end{pspicture}
\]
 \caption{Move $P_{III}$}\label{3III}
\end{figure}

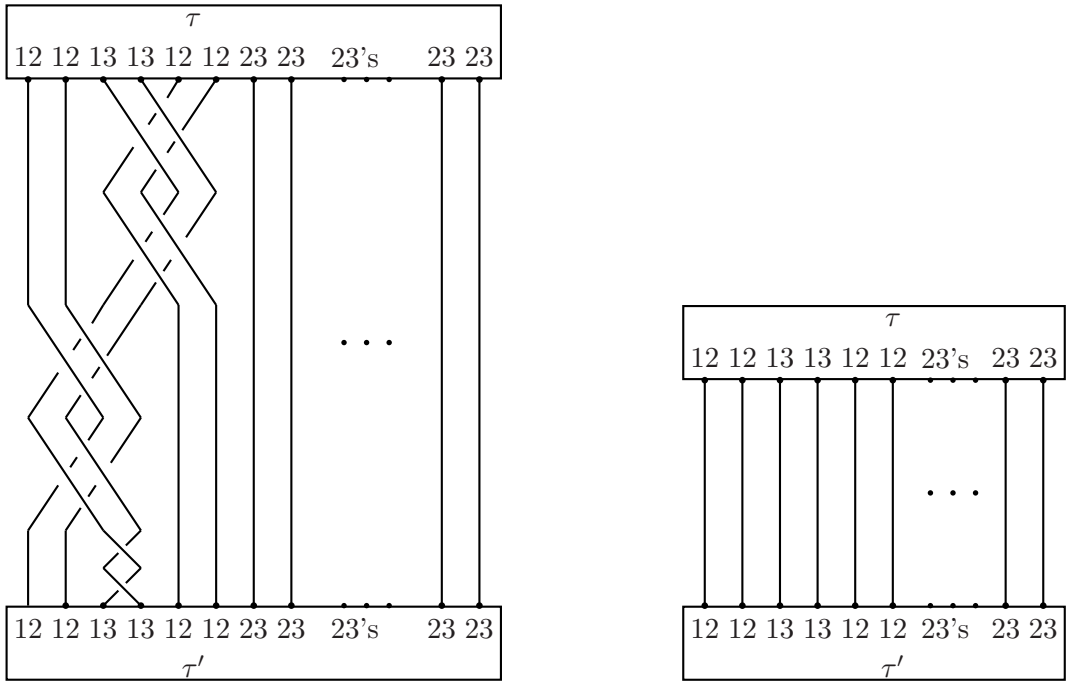
\begin{figure}[htp]
\[
\begin{pspicture}(-.5,-4)(6.5,5)
 \psdots(0,4)(.5,4)(1,4)(1.5,4)(2,4)(2.5,4)(3,4)(3.5,4)(5.5,4)(6,4)
        (.5,-3)(1,-3)(1.5,-3)(2,-3)(2.5,-3)(3,-3)(3.5,-3)(5.5,-3)(6,-3)
   \psline(0,4)(0,1)
   \psline(0.5,4)(0.5,1)
   \psline(3,4)(3,-3)
   \psline(3.5,4)(3.5,-3)
   \psline(5.5,4)(5.5,-3)
   \psline(6,4)(6,-3)
  \psline(1,4)(2,2.5)
  \psline(1.5,4)(2.5,2.5)
  \psline(1,2.5)(2,1)
  \psline(1.5,2.5)(2.5,1)
   \psline(1,2.5)(1.4,3.1)
   \psline(1.6,3.4)(1.7,3.55)
   \psline(1.85,3.755)(2,4)
    \psline(1.5,2.5)(1.7,2.8)
    \psline(1.85,3.025)(1.95,3.175)
    \psline(2.05,3.325)(2.5,4)
    \psline(2,2.5)(1.8,2.2)
    \psline(1.65,1.975)(1.56,1.84)
    \psline(1.4,1.6)(1,1)
    \psline(2.5,2.5)(2.1,1.9)
    \psline(1.9,1.6)(1.8,1.45)
    \psline(1.7,1.3)(1.5,1)
           \psframe(-.3,-4)(6.3,-3)
           \psframe(-.3,4)(6.3,5)
  \rput(0,4.3){12}
 \rput(.5,4.3){12}
 \rput(1,4.3){13}
 \rput(1.5,4.3){13}
 \rput(2,4.3){12}
 \rput(2.5,4.3){12}
 \rput(3,4.3){23}
 \rput(3.5,4.3){23}
 \rput(5.5,4.3){23}
 \rput(6,4.3){23}
                   \rput(4.35,4.3){23's}
  \rput(0,-3.3){12}
 \rput(.5,-3.3){12}
 \rput(1,-3.3){13}
 \rput(1.5,-3.3){13}
 \rput(2,-3.3){12}
 \rput(2.5,-3.3){12}
 \rput(3,-3.3){23}
 \rput(3.5,-3.3){23}
 \rput(5.5,-3.3){23}
 \rput(6,-3.3){23}
                   \rput(4.35,-3.3){23's}
  \psdots[dotsize=0.1](4.2,4)(4.5,4)(4.8,4)(4.2,-3)(4.5,-3)(4.8,-3)
                      (4.2,.5)(4.5,.5)(4.8,.5)
     \rput(2.2,4.8){$\tau$}
     \rput(2.2,-3.8){$\tau'$}
   \psline(0,1)(1,-.5)
  \psline(.5,1)(1.5,-.5)
  \psline(0,-.5)(1,-2)
  \psline(.5,-.5)(1.5,-2)
   \psline(0,-.5)(.4,.1)
   \psline(.6,.4)(.7,.55)
   \psline(.85,.755)(1,1)
    \psline(.5,-.5)(.7,-.2)
    \psline(.85,.025)(.95,.175)
    \psline(1.05,.325)(1.5,1)
    \psline(1,-.5)(.8,-.8)
    \psline(.65,-1.025)(.56,-1.16)
    \psline(.4,-1.4)(0,-2)
    \psline(1.5,-.5)(1.1,-1.1)
    \psline(.9,-1.4)(.8,-1.55)
    \psline(.7,-1.7)(.5,-2)
 \psline(1,-2)(1.5,-2.5)
 \psline(1.5,-2)(1.3,-2.2)
 \psline(1.2,-2.3)(1,-2.5)
   \psline(1,-2.5)(1.5,-3)
   \psline(1.5,-2.5)(1.3,-2.7)
   \psline(1.2,-2.8)(1,-3)
  \psline(0,-2)(0,-3)
  \psline(.5,-2)(.5,-3)
  \psline(2,1)(2,-3)
  \psline(2.5,1)(2.5,-3)
 \end{pspicture}
\hspace{2cm}
\begin{pspicture}(-.5,0)(5,5)
 \psdots(0,4)(.5,4)(1,4)(1.5,4)(2,4)(2.5,4)(4,4)(4.5,4)
        (0,1)(.5,1)(1,1)(1.5,1)(2,1)(2.5,1)(4,1)(4.5,1)
   \psline(1,4)(1,1)
   \psline(1.5,4)(1.5,1)
   \psline(2,4)(2,1)
   \psline(2.5,4)(2.5,1)
   \psline(4,4)(4,1)
   \psline(4.5,4)(4.5,1)
 \psline(0,4)(0,1)
 \psline(0.5,4)(0.5,1)
             \psframe(-.3,0)(4.8,1)
             \psframe(-.3,4)(4.8,5)
 \rput(0,4.3){12}
 \rput(0.5,4.3){12}
 \rput(1,4.3){13}
 \rput(1.5,4.3){13}
 \rput(2,4.3){12}
 \rput(2.5,4.3){12}
 \rput(4,4.3){23}
 \rput(4.5,4.3){23}
                   \rput(3.2,4.3){23's}
 \rput(0,.7){12}
 \rput(0.5,.7){12}
 \rput(1,.7){13}
 \rput(1.5,.7){13}
 \rput(2,.7){12}
 \rput(2.5,.7){12}
 \rput(4,.7){23}
 \rput(4.5,.7){23}
                   \rput(3.2,.7){23's}
  \psdots[dotsize=0.1](3,4)(3.3,4)(3.6,4)(3,1)(3.3,1)(3.6,1)
                      (3,2.5)(3.3,2.5)(3.6,2.5)
     \rput(2.5,4.8){$\tau$}
     \rput(2.5,.2){$\tau'$}
\end{pspicture}
\]
\caption{Move $P_{IV}$}\label{3IV}
\end{figure}

\begin{thm}[\bf Piergallini]\label{pier1} Two tricolored link diagrams represent the same
manifold iff they can be related (up to colored Reidemeister
moves) by a finite sequence of moves of the types $\mathcal{M}$,
$P_{II}$, $P_{III}$, $P_{IV}$ described in Figures~\ref{moveM},
\ref{3II}, \ref{3III},\ref{3IV} respectively.
\end{thm}

\newpage

\begin{proof} Move $\mathcal{M}$ does not change the represented
manifold since the braid on its left side lifts to identity. That
the moves II, III, and IV do not change the represented manifold
follows for example from the proof of Theorem~\ref{pier2} below
where it is shown that they can be realized using moves
$\mathcal{M}$, $\mathcal{P}$ and addition/deletion of a trivial
sheet. (There is no vicious circle here since this part of the
theorem is not used in the proof of \ref{pier2}).\

For the reverse implication we just give a sketch of the contents
of~\cite{P1}. Let $D$ and $D'$ be two tricolored links
representing the same $3$-manifold.

Convert $D$ and $D'$ to normalized diagrams (still denoted $D$ and
$D'$), see Section~\ref{3sphere}. This gives us two Heegaard
diagrams of the same manifold which according to
Appendix~\ref{Heegaard} are stably equivalent. The proof then
proceeds as follows:

\begin{description}
\item[Step $1$] Stabilization of Heegaard splittings can be
realized at the level of normalized diagrams via isotopy.
\end{description}

 Indeed performing the modification shown in Figure~\ref{Heegstab}
 the corresponding Heegaard splitting changes by a connected sum
 with the genus-$1$ splitting of $S^3$.

  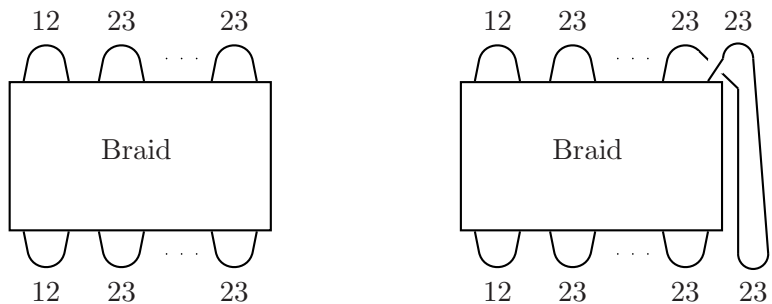
\begin{figure}[htp]
 \[
 \begin{pspicture}(4,5)
 \psline[linearc=.25](.5,3.5)(.8,5)(1.1,3.5)
 \psline[linearc=.25](1.5,3.5)(1.8,5)(2.1,3.5)
 \psline[linearc=.25](3,3.5)(3.3,5)(3.6,3.5)
 \psdots[dotsize=.05](2.4,3.8)(2.6,3.8)(2.8,3.8)
 \psframe(.3,3.5)(3.8,1.5)
 \rput(2,2.6){Braid}
 \rput(.8,4.3){12}
 \rput(1.8,4.3){23}
 \rput(3.3,4.3){23}
 \psline[linearc=.25](.5,1.5)(.8,0)(1.1,1.5)
 \psline[linearc=.25](1.5,1.5)(1.8,0)(2.1,1.5)
 \psline[linearc=.25](3,1.5)(3.3,0)(3.6,1.5)
 \psdots[dotsize=.05](2.4,1.2)(2.6,1.2)(2.8,1.2)
 \rput(.8,.7){12}
 \rput(1.8,.7){23}
 \rput(3.3,.7){23}
 \end{pspicture}
 \hspace{2cm}
  \begin{pspicture}(5,5)
 \psline[linearc=.25](.5,3.5)(.8,5)(1.1,3.5)
 \psline[linearc=.25](1.5,3.5)(1.8,5)(2.1,3.5)
 \psline[linearc=.25](3,3.5)(3.3,5)(3.6,3.5)
 \psdots[dotsize=.05](2.4,3.8)(2.6,3.8)(2.8,3.8)
 \psframe[linecolor=white,fillstyle=solid,fillcolor=white](3.4,3.5)(3.7,3.9)
 \psline(3.48,3.9)(4,3.4)
 \psframe[linecolor=white,fillstyle=solid,fillcolor=white](3.6,3.6)(3.8,3.8)
 \psline(3.6,3.5)(3.8,3.8)
\psline(4,3.4)(4,1.2)
 \psarc(4.2,1.2){.2}{180}{0}
 \psline(4.4,1.2)(4.2,3.8)
 \psarc(4,3.8){.2}{0}{180}
 \psframe(.3,3.5)(3.8,1.5)
 \rput(2,2.6){Braid}
 \rput(.8,4.3){12}
 \rput(1.8,4.3){23}
 \rput(3.3,4.3){23}
 \rput(4,4.3){23}
 \psline[linearc=.25](.5,1.5)(.8,0)(1.1,1.5)
 \psline[linearc=.25](1.5,1.5)(1.8,0)(2.1,1.5)
 \psline[linearc=.25](3,1.5)(3.3,0)(3.6,1.5)
 \psdots[dotsize=.05](2.4,1.2)(2.6,1.2)(2.8,1.2)
 \rput(.8,.7){12}
 \rput(1.8,.7){23}
 \rput(3.3,.7){23}
 \rput(4.2,.7){23}
 \end{pspicture}
 \]
 \caption{Heegaard Stabilization}\label{Heegstab}
 \end{figure}

Therefore by repeatedly applying the modification in
Figure~\ref{Heegstab} one can change the diagrams $D$ and $D'$ so
that they represent equivalent Heegaard splittings i.e. so that we
have a commutative diagram:
$$
\xymatrix{ {H_1}\ar[d]& {\supset}&{\partial
H_1}\ar[d]^{f_1}\ar[r]^{\phi}&{\partial H_2}
\ar[d]^{f_2}& {\subset}&{H_2}\ar[d]\\
{H_1}& {\supset}&{\partial H_1}\ar[r]^{\phi'}&{\partial H_2}&
{\subset} &{H_2} }
$$

where $M=H_1\underset{\phi}{\bigcup}H_2$ and
$M=H_1\underset{\phi'}{\bigcup}H_2$ are the Heegaard splittings
corresponding to the diagrams $D$ and $D'$ respectively, and the
vertical maps are homeomorphisms. In particular the homeomorphisms
$f_1$ and $f_2$ are homeomorphisms of surfaces that extend to
homeomorphisms of the handlebodies.

\begin{description}
\item[Step $2$] For any homeomorphism $\xymatrix@1{
{\phi:\partial H_1}\ar[r]&{\partial H_1} }$ that extends to a
homeomorphism $\xymatrix@1{ {H_1}\ar[r]&{H_1} }$ there is a
tricolored braid that lifts to $\phi$ and can be added at the top
(or bottom) of a normalized diagram via isotopy and move
$\mathcal{M}$.
\end{description}

 This is done using the results of Suzuki in \cite{S} where a finite set
  of generators for the group of isotopy classes
of homeomorphisms of a surface that extend to homeomorphisms of
the handlebody is given. Piergallini in Section 2 of \cite{P1}
shows that each of these generators can be added at the top (and
hence at the bottom) of a normalized diagram using isotopy and
move $\mathcal{M}$.\

 Therefore, referring to the commutative diagram above, one can add a colored braid
lifting to $f_1$ at the top, and a colored braid lifting to
$f_2^{-1}$ at the bottom of $D'$. The resulting diagram (still
denoted $D'$) induces the same Heegaard diagram as $D$ that is,
the braids of $D$ and $D'$ lift to the same homeomorphism of
$\partial H$. So if $\B$ (resp. $\B'$) is the braid of $D$ (resp.
$D'$) then
$$\exists \B''\in Ker(\lambda)\quad \B=\B'\B''$$
where $\lambda$ is the lifting homomorphism.

\begin{description}
\item[Step 3] Any braid in the kernel of the lifting
homeomorphism can be added at the bottom of a normalized diagram
using the moves $\mathcal{M}$, II, III and IV.
\end{description}

 Piergallini shows that in Section 3 of \cite{P1} by showing that
 each of the normal generators in Theorem~\ref{Ker3} can be added
 at the bottom of a normalized diagram.\

 Therefore the braid $\B''$ of the above equation can be added at
 the bottom of $D'$ to convert it to $D$.

\end{proof}

Notice the following difference between move $\mathcal{M}$ and the
other three moves: in order to check that moves $P_{II}$-$P_{IV}$
can be applied one needs to know the structure of the whole link
while for move $\mathcal{M}$ one only needs to know a small part
of the link. In this sense $\mathcal{M}$ is a ``local'' move while
the other three are not.\

\begin{figure}[htp]
\[
\begin{pspicture}(3,3)
\psframe(3,3)
 \rput(1.3,1.7){tricolored}
  \rput (1.5,1.4){D}
  \end{pspicture}
  \hspace{2.5cm}
\begin{pspicture}(5,3)
\psframe(3,3)
 \rput(1.3,1.7){tricolored}
  \rput (1.5,1.4){D}
 \pscircle(4,1.5){.5}
 \rput(4.9,1.5){(i4)}
  \end{pspicture}
\]
\caption{Adding a trivial sheet}\label{next}
\end{figure}
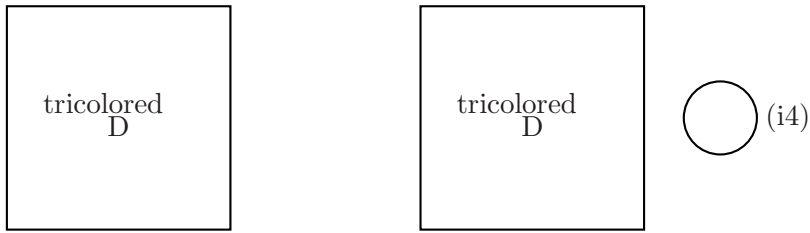

 Recall from Section~\ref{trivialsheet}
the procedure of adding a trivial sheet to coverings over spheres.
In the case of a $3$-sheeted covering over $S^3$ it produces a
$4$-sheeted covering over $S^3$ with the same total space. In
terms of colored diagrams we have the following move
(Figure~\ref{next}): to a tricolored link diagram $D$ add an
unknotted unlinked component colored by a dark color (see
Section~\ref{dim}) to get a bi-tricolored link diagram. This move
will be referred to as addition of a trivial sheet and its reverse
move as deletion of a trivial sheet.

After proving in \cite{M2} that $\mathcal{M}$ is not enough
Montesinos asked in the same paper whether  moves $\mathcal{M}$,
$\mathcal{P}$ (see Figure~\ref{moveP}) together with
addition/deletion of a fourth trivial sheet suffice to relate any
two tricolored link presentations of the same manifold.
Piergallini gave a positive answer in \cite{P2}:

\begin{thm}[\bf Piergallini]\label{pier2} Two tricolored link diagrams represent
the same manifold iff they can be related using moves
$\mathcal{M}$ , $\mathcal{P}$ and addition/deletion of a fourth
trivial sheet (and of course colored Reidemeister moves).
\end{thm}

\begin{proof} Piergallini proves in Theorem A of \cite{P2} that the moves
II, III and IV can be realized using moves $\mathcal{M}$ and
$\mathcal{P}$ after one adds a trivial sheet.
\end{proof}

\section{Four sheets I: Coverings over $S^2$}\label{2sphere4} In this
section we establish for $4$-sheeted coverings of $S^2$ results
similar to the results of Section~\ref{2sphere3}.  Again it is
more convenient to work with coverings over the disc and then
restrict attention to coverings with boundary monodromy id, see
Section~\ref{2D}, page~\pageref{2sphere}.  Specifically we
consider the following $4$-sheeted covering of the $2$-disc:
$$\rho_{23}^n:=(12),(12),(14),(14),(23),\dotsc,(23)$$
 where there are $n$ branch values and all the missing monodromies
 are equal to $(23)$. (This notation is part of a more general
 notation, see Definition~\ref{notation} below.) After dimming
 the lights (see Section~\ref{dim}) we get the $3$-sheeted
 covering $\rho(n)$ of Section~\ref{2sphere3} and since dimming
 the lights commutes with lifting we get that

  $$\Au(\rho_{23}^n)\subset L(n)$$
where $L(n)$ denotes  $\Au(\rho(n)$ as in Section~\ref{2sphere3}.
Furthermore if $\B\in L(n)$ then dimming the lights for
$(\rho_{23}^n)\B$ will again give $\rho(n)$. Thus we have an
action of $L(n)$ on the set $\kappa^{-1}(\rho(n))$ of $4$-sheeted
coverings that give $\rho(n)$ after dimming the lights, and the
isotropy group of $\rho_{23}^n$ under this action is
$\Au(\rho_{23}^n)$.\

 Using this action in Subsection~\ref{thecomplex} we construct a
 $2$-dimensional cell complex whose fundamental group is the
 quotiented automorphism group
 $\overline{\Au}(\rho_{23}^n)$. The complex is then used in
 Subsection~\ref{thegenerators} to produce generators for this
 group. Finally in Subsection~\ref{thekernel} we use these
 generators to find normal generators for the kernel of the
 (quotiented) lifting homeomorphism.

\subsection{The Complex}\label{thecomplex} Recall that $\kappa$
denotes dimming the lights, $L=\{A_0,A_1,\dotsc,A_{n-1}\}$ denotes the set
 of branched values and that $\alpha_i$
for $i=0,1,\dotsc,n-1$ denotes the generator of $\pi_1(D^2\setminus L)$
 that is represented by a lasso going
around $A_i$.

\begin{defn}\label{notation} Let
$$\CC^n:=\{\rho\in\kappa^{-1}(\rho(n))\mid
\rho(\alpha_0)=\rho(\alpha_1)=(12)\}$$ For $I\subseteq
\{2,3,\dotsc,n-1\}$ define the $4$-sheeted coverings $\rho_I^n$,
$\tilde\rho_I^n$ of $D^2$ via:
\begin{align*} \rho_I^n(\alpha_i):=\begin{cases} (12) & \text{if }
i=0,1\\ (14) & \text{if } i\in I\\ (23) & \text{if } i\notin I
\end{cases} &\quad,\qquad \tilde\rho_I^n(\alpha_i)=\begin{cases} (12) & \text{if }
i=0,1\\ (23) & \text{if } i\in I\\ (14) & \text{if } i\notin I
\end{cases}
\end{align*}

\end{defn}

Recall the generators of $L(n)$ given in Theorem~\ref{bw1}. The following lemma
describes how $\delta_4$ acts on
$$\CC^6=\{\rho_{\emptyset}^6,\tilde\rho_{\emptyset}^6,\rho_{ij}^6,
\tilde\rho_{ij}^6,\rho_i^6,\tilde\rho_i^6 \mid i,j=2,3,4,5\}.$$
Since $(\alpha_i)\delta_4=\alpha_i$ for $i\geq 6$ this determines
 how $\delta_4$ acts on elements of $\CC^n$.
 The action of the remaining generators is clear.

\begin{lem}\label{d-action} $\delta_4$ acts on $\CC^6$
 as follows: it fixes $\rho_{\emptyset}$,
 $\tilde\rho_{\emptyset}$, $\rho_{ij}$ and $\tilde\rho_{ij}$,
 while $(\rho_i)\delta_4=\tilde\rho_i$.
\end{lem}
 \begin{proof} Since $\delta_4\in L(n)$ it fixes $\rho_{\emptyset}$,
 $\tilde\rho_{\emptyset}$. For the remaining cases see Figure \ref{fig:daction}.

\begin{figure}[htp]
\[
\begin{pspicture}(1,1.5)(5,9.5)
 \psline(4.5,8)(1.5,7)
 \psline(1.5,4)(4.5,3)
 \psline(2,9)(2,7.25)
 \psline(2,7.05)(2,3.95)
 \psline(2,3.7)(2,1.833)
 \psline(3,9)(3,7.6)
 \psline(3,7.4)(3,3.6)
 \psline(3,3.4)(3,1.833)
 \psline(4,9)(4,7.95)
 \psline(4,7.7)(4,3.3)
 \psline(4,3.05)(4,1.833)
 \psline(1,9)(1.9,8.746)
 \psline(2.1,8.68)(2.9,8.45)
 \psline(3.1,8.39)(3.9,8.16)
 \psline(4.1,8.1)(4.5,8)
 \psline(1.5,7)(1.9,6.867)
 \psline(2.1,6.8)(2.9,6.533)
 \psline(3.1,6.467)(3.9,6.2)
 \psline(4.1,6.133)(4.5,6)
 \psline(4.5,6)(5,5.5)
 \psline(5,5.5)(5,1.833)
 \psline(5,9)(5,6.2)
 \psline(5,6.2)(4.75,5.9)
 \psline(4.62,5.744)(4.1,5.12)
 \psline(3.9,4.96)(3.1,4.64)
 \psline(2.9,4.56)(2.1,4.24)
 \psline(1.9,4.16)(1.5,4)
 \psline(4.5,3)(4.1,2.867)
 \psline(3.9,2.8)(3.1,2.533)
 \psline(2.9,2.467)(2.1,2.2)
 \psline(1.9,2.133)(1.5,2)
 \psline(1.5,2)(1,1.833)
 \psdots(1,9)(2,9)(3,9)(4,9)(5,9)
        (1,1.833)(2,1.833)(3,1.833)(4,1.833)(5,1.833)
 \rput(1,9.3){12}
 \rput(2,9.3){23}
 \rput(3,9.3){23}
 \rput(4,9.3){14}
 \rput(5,9.3){14}
 \rput(1,1.533){12}
 \rput(2,1.533){23}
 \rput(3,1.533){23}
 \rput(4,1.533){14}
 \rput(5,1.533){14}
 \rput(1.2,7){24}
 \rput(1.25,4){24}
 \rput(3.4,5){24}
 \rput(3.6,2.35){12}
 \rput(4.3,7){12}
 \rput(1.7,5.8){34}
 \rput(2.7,5.8){34}
\end{pspicture}
\hspace{3cm}
\begin{pspicture}(1,1.5)(5,9.5)
 \psline(4.5,8)(1.5,7)
 \psline(1.5,4)(4.5,3)
 \psline(2,9)(2,7.25)
 \psline(2,7.05)(2,3.95)
 \psline(2,3.7)(2,1.833)
 \psline(3,9)(3,7.6)
 \psline(3,7.4)(3,3.6)
 \psline(3,3.4)(3,1.833)
 \psline(4,9)(4,7.95)
 \psline(4,7.7)(4,3.3)
 \psline(4,3.05)(4,1.833)
 \psline(1,9)(1.9,8.746)
 \psline(2.1,8.68)(2.9,8.45)
 \psline(3.1,8.39)(3.9,8.16)
 \psline(4.1,8.1)(4.5,8)
 \psline(1.5,7)(1.9,6.867)
 \psline(2.1,6.8)(2.9,6.533)
 \psline(3.1,6.467)(3.9,6.2)
 \psline(4.1,6.133)(4.5,6)
 \psline(4.5,6)(5,5.5)
 \psline(5,5.5)(5,1.833)
 \psline(5,9)(5,6.2)
 \psline(5,6.2)(4.75,5.9)
 \psline(4.62,5.744)(4.1,5.12)
 \psline(3.9,4.96)(3.1,4.64)
 \psline(2.9,4.56)(2.1,4.24)
 \psline(1.9,4.16)(1.5,4)
 \psline(4.5,3)(4.1,2.867)
 \psline(3.9,2.8)(3.1,2.533)
 \psline(2.9,2.467)(2.1,2.2)
 \psline(1.9,2.133)(1.5,2)
 \psline(1.5,2)(1,1.833)
 \psdots(1,9)(2,9)(3,9)(4,9)(5,9)
        (1,1.833)(2,1.833)(3,1.833)(4,1.833)(5,1.833)
 \rput(1,9.3){12}
 \rput(2,9.3){23}
 \rput(3,9.3){23}
 \rput(4,9.3){23}
 \rput(5,9.3){14}
 \rput(1,1.533){12}
 \rput(2,1.533){14}
 \rput(3,1.533){14}
 \rput(4,1.533){14}
 \rput(5,1.533){23}
 \rput(1.2,7){13}
 \rput(1.25,4){24}
 \rput(3.4,5){24}
 \rput(3.6,2.35){12}
 \rput(4.3,7){12}
 \rput(1.7,5.8){12}
 \rput(2.7,5.8){12}
\end{pspicture}
\]
\caption{How $\delta_4$ acts}\label{fig:daction}
\end{figure}
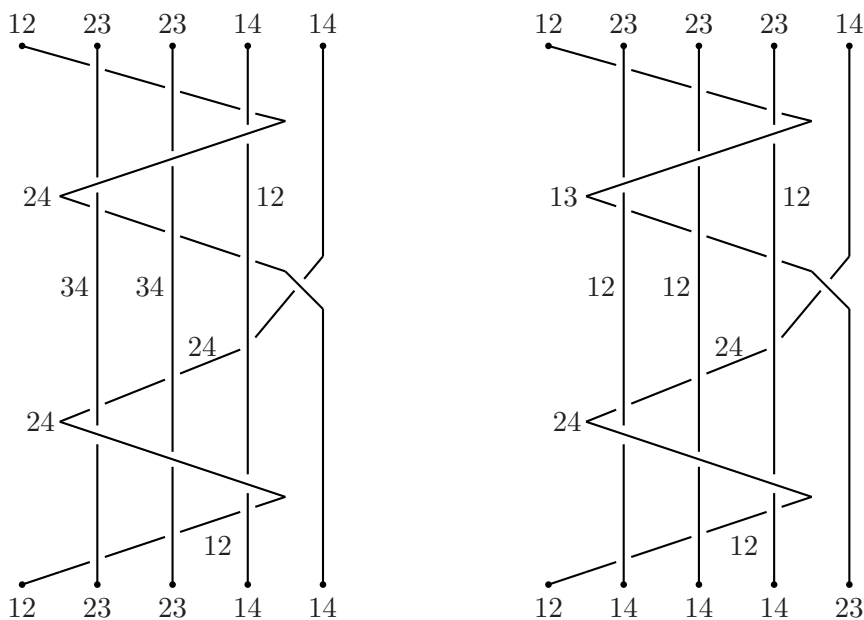

 Notice that these two cases are enough since $\delta_4$ commutes
 with $\B_2$, $\B_3$ and, of course, with conjugations in the
 range. Furthermore by Lemma \ref{Com} it commutes with $\B_4$
  modulo the move $\mathcal{M}$.
\end{proof}

\begin{cor} The action of $L(n)$ on $\kappa^{-1}(\rho(n)$
restricts to an action on $\CC^n$.
\end{cor}

\begin{proof} We have to check that each of the generators of
$L(n)$ preserves $\CC(n)$, that is if the first two monodromies of
a covering $\rho$ are equal to $(12)$ then the first two
monodromies of $(\rho)\B$ are also equal to $(12)$, for $\B\in
L(n)$. Lemma~\ref{d-action} proves this for $\delta_4$, for the
other generators it is obvious.
\end{proof}

 Recall from page~\pageref{N} that for a covering $\rho$, $N(\rho)$ denotes
 the normal subgroup of $\Au(\rho)$ ``generated'' by the moves
 $\mathcal{M}$ and $\mathcal{P}$ and that the quotient of
 $\Au(\rho)$ by that subgroup is denoted by
 $\overline{\Au}(\rho)$.\\
 A standard technique for understanding a group is to realize it
 as the fundamental group of some $2$-complex. We construct such a
 complex for the quotiented automorphism group:

\begin{defn} Given $\rho\in\kappa^{-1}(\rho(n))$ define a $2$-complex
 $\CC(\rho)$ as follows:
\begin{itemize}
\item The $0$-cells of $\CC(\rho)$ are the points of the
$L(n)$-orbit of $\rho$.
\item The $1$-cells of $\CC(\rho)$ are labelled by the
generators of $L(n)$. There is one $1$-cell labelled by $\B$ for
each {\em unordered} pair
$\{\rho',\rho''\}\subseteq\kappa^{-1}(\rho(n))$ with the property
$\rho'\B=\rho''$.
\item A path on the $1$-skeleton gives a word on the generators of
$L(n)$. There is a $2$-cell attached along each closed path whose
word lies in $N(\rho')$ for some $\rho'$ in the path.
\end{itemize}
\end{defn}

Notice that each generator of $L(n)$ acts as an involution and
therefore we don't need to orient the $1$-dimensional cells.
 Clearly $\pi_1(\CC(\rho),\rho)=\overline{\Au}(\rho)$. Also notice
that one can get a complex whose fundamental group is $\Au(\rho)$
by attaching only the $2$-cells that kill the words that are equal
to identity in $L(n)$.\

The complex $\CC(\rho)$ will be confused with its $0$-skeleton
when the confusion is easily resolved.\\

\underline{Terminology:} A vertex is a $0$-cell. An edge is a
$1$-cell with its boundary points attached to distinct vertices. A
loop is a $1$-cell with its boundary points attached to the same
vertex.

\begin{lem}\label{this} At the level of $0$-skeletons:
\begin{itemize}
\item[(a)] $\CC(\rho_{23}^n)$ consists of those coverings which:
\begin{itemize}
\item start with $(12), (12)$,
\item have at least one monodromy equal to $(14)$, and
at least one monodromy equal to $(23)$,
\item have an even number of monodromies equal to
$(14)$.
\end{itemize}
\item[(b)]  $\CC(\tilde\rho_{23}^n)$ consists of those coverings which:
\begin{itemize}
\item start with $(12), (12)$,
\item have at least one monodromy equal to $(14)$, and
at least one monodromy equal to $(23)$,
\item have an even number of monodromies equal to
$(23)$.
\end{itemize}
\item[(c)] For n even, $\CC(\rho_{2}^n)$ consists of those coverings which:
\begin{itemize}
\item start with $(12), (12)$,
\item have at least one monodromy equal to $(14)$, and
at least one monodromy equal to $(23)$,
\item have an odd number of monodromies equal to
$(14)$.
\end{itemize}
\end{itemize}
\end{lem}

\begin{proof} Let $\CC_{\text{id}}^n$ (resp. $\CC_{(23)}^n$)
 be the set of coverings with the properties described in (a) and
 $n$ even (resp. odd).
  Clearly $\rho_{23}^n\in\CC_{\text{id}}^n$ (resp. $\CC_{(23)}^n$)  and
$L(n)$ preserves $\CC_{\text{id}}^n$ (resp. $\CC_{(23)}^n$). This
is obvious for $\B_0$, $\B_1^3$, $\B_2,\cdots, \B_{n-2}$ and
\ref{d-action} proves it for $\delta_4$. \\
To see that every covering in $\CC_{\text{id}}^n$ (resp.
$\CC_{(23)}^n$) is in the $L(n)$-orbit of $\rho_{23}^n$, observe
that $\langle\B_2,\cdots,\B_{n-2}\rangle$ acts transitively on the
subset of $\CC_{\text{id}}^n$ (resp. $\CC_{(23)}^n$) consisting of
coverings with any fixed even number of monodromies equal to
$(14)$, while $\delta_4$, acting on suitably chosen coverings,
increases
or decreases the number of monodromies equal to $(14)$ by $2$.\\
(b) and (c) are proved similarly by (defining and) looking at
$\CC_{(14)}^n$ and $\CC_{(14)(23)}^n$.
\end{proof}
\begin{prop} $\CC_{\bullet}^n$ (defined in the proof of Lemma~\ref{this})
 consists of those coverings of
$\CC^n$ that have connected total space and boundary monodromy
equal to $\bullet$.
\end{prop}
\begin{proof} Obvious. \end{proof}

 Conjugation by $(12)(34)$ (i.e interchanging the dark and light
shades of Blue) defines an involution $\quad\widetilde{}\quad$ on
$\CC^n$. Clearly $\quad\widetilde{}\quad$ interchanges
$\CC_{(23)}^n$ and $\CC_{(14)}^n$ and preserves
$\CC_{\text{id}}^n$ and $\CC_{(14)(23)}^n$. In particular
$\CC_{(23)}^n$ and $\CC_{(14)}^n$ are isomorphic.\

 Now given a covering with $n$ branch values one can ``add''
one more point with prescribed monodromy at the end to get a
covering with $n+1$ branch values. More formally, for each $n$
there are two embeddings:
$$ i_{(23)}^n:\CC^n\longrightarrow\CC^{n+1},$$
and
$$i_{(14)}^n:\CC^n\longrightarrow\CC^{n+1}$$
which are obviously $L(n)$-equivariant if $B_n$ is (as usual)
considered to be the subgroup of $B_{n+1}$ that ``doesn't touch''
the last string.

 \begin{prop}\label{dec}
 For even $n$:
 \begin{itemize}
 \item[a)] At the level of $0$ skeletons:
 $$\CC_{\text{id}}^n=i_{(23)}^{n-1}(\CC_{(23)}^{n-1})\bigsqcup
 i_{(14)}^{n-1}(\CC_{(14)}^{n-1}),$$
and each vertex of $i_{(23)}i_{(14)}(\CC_{(14)(23)}^{n-2})$ is
joined by an edge labelled by $\B_{n-2}$ to the corresponding
vertex of $i_{(14)}i_{(23)}(\CC_{(14)(23)}^{n-2})$. These are the
only new edges.

\item[b)] At the level of $0$-skeletons:
 $$\CC_{(14)(23)}^n=i_{(23)}^{n-1}(\CC_{(14)}^{n-1})\bigsqcup
 i_{(14)}^{n-1}(\CC_{(23)}^{n-1})$$
and each vertex of $i_{(23)}i_{(14)}(\CC_{\text{id}}^{n-2}
\sqcup\{\rho_\emptyset^{n-2}\}\sqcup\{\tilde\rho_\emptyset^{n-2}\})$
is joined by an edge labelled by $\B_{n-2}$ to the corresponding
vertex of $i_{(14)}i_{(23)}(\CC_{\text{id}}^{n-2}
\sqcup\{\rho_\emptyset^{n-2}\}\sqcup\{\tilde\rho_\emptyset^{n-2}\})$.
These are the only new edges.

\end{itemize}
 For odd $n$:
\begin{itemize}
\item[] At the level of $0$ skeletons:
 $$\CC_{(23)}^n=i_{(23)}^{n-1}(\CC_{\text{id}}^{n-1}\sqcup \{\tilde\rho_\emptyset^{n-1}\})
 \bigsqcup i_{(14)}^{n-1}(\CC_{(14)(23)}^{n-1}),$$
and each vertex of $i_{(23)}i_{(14)}(\CC_{(14)}^{n-2}
\sqcup\{\tilde\rho_\emptyset^{n-2}\})$ is joined by an edge
labelled by $\B_{n-2}$ to the corresponding vertex of
$i_{(14)}i_{(23)}(\CC_{\text{id}}^{n-2}
\sqcup\{\tilde\rho_\emptyset^{n-2}\})$. These are the only new
edges.

\end{itemize}
 \end{prop}
\begin{proof} Erase the last two points.
 \end{proof}

 \begin{lem}\label{Com} $\B_4$ commutes with $\delta_4$ modulo the Montesinos
move.
\end{lem}

\begin{proof}
It suffices to show that $\delta_4\B_4\delta_4^{-1}\B_4^{-1}\equiv \text{id}
 \pmod{\mathcal{M}}$. This is done in
the following pictures. The top of the braid is assumed colored by $\rho(n)$.
The colors of the arcs are not shown
as they can be easily deduced.

\newpage
Start with $\delta_4\B_4\delta_4^{-1}\B_4^{-1}$ and perform $\mathcal{M}$ inside
the dotted circles:\\

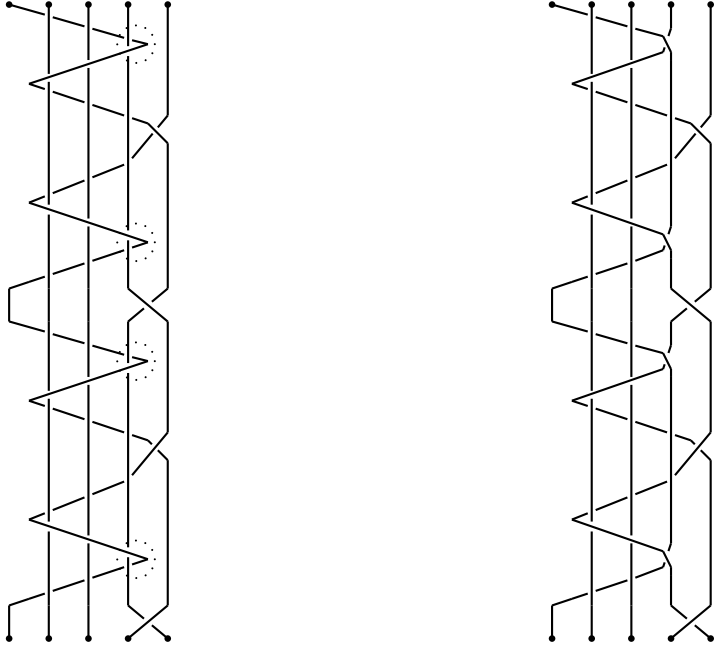
\[
\psset{unit=15pt}
\begin{pspicture}(0,-7)(8,9)
\psdots(1,9)(2,9)(3,9)(4,9)(5,9)               %
       (5,-7)(1,-7)(2,-7)(3,-7)(4,-7)               %
 \psline(4.5,8)(1.5,7)
 \psline(1.5,4)(4.5,3)
 \psline(2,9)(2,7.25)
 \psline(2,7.05)(2,3.95)
 \psline(2,3.7)(2,1.833)
 \psline(3,9)(3,7.6)
 \psline(3,7.4)(3,3.6)
 \psline(3,3.4)(3,1.833)
 \psline(4,9)(4,7.95)
 \psline(4,7.7)(4,3.3)
 \psline(4,3.05)(4,1.833)
 \psline(1,9)(1.9,8.746)
 \psline(2.1,8.68)(2.9,8.45)
 \psline(3.1,8.39)(3.9,8.16)
 \psline(4.1,8.1)(4.5,8)
 \psline(1.5,7)(1.9,6.867)
 \psline(2.1,6.8)(2.9,6.533)
 \psline(3.1,6.467)(3.9,6.2)
 \psline(4.1,6.133)(4.5,6)
 \psline(4.5,6)(5,5.5)
 \psline(5,5.5)(5,1.833)
 \psline(5,9)(5,6.2)
 \psline(5,6.2)(4.75,5.9)
 \psline(4.62,5.744)(4.1,5.12)
 \psline(3.9,4.96)(3.1,4.64)
 \psline(2.9,4.56)(2.1,4.24)
 \psline(1.9,4.16)(1.5,4)
 \psline(4.5,3)(4.1,2.867)
 \psline(3.9,2.8)(3.1,2.533)
 \psline(2.9,2.467)(2.1,2.2)
 \psline(1.9,2.133)(1.5,2)
 \psline(1.5,2)(1,1.833)
      \psline(1,1.833)(1,1)
      \psline(2,1.833)(2,1)
      \psline(3,1.833)(3,1)
     \psline(4,1.833)(5,1)
     \psline(5,1.833)(4.6,1.4998)
     \psline(4.4,1.3332)(4,1)
   \rput(.12,.5) {%
  \begin{pspicture}(0,17)
 \psline(4.5,8)(1.5,7)
 \psline(1.5,4)(4.5,3)
 \psline(2,9)(2,7.25)
 \psline(2,7.05)(2,3.95)
 \psline(2,3.7)(2,1.833)
 \psline(3,9)(3,7.6)
 \psline(3,7.4)(3,3.6)
 \psline(3,3.4)(3,1.833)
 \psline(4,9)(4,7.95)
 \psline(4,7.7)(4,3.3)
 \psline(4,3.05)(4,1.833)
 \psline(1,9)(1.9,8.746)
 \psline(2.1,8.68)(2.9,8.45)
 \psline(3.1,8.39)(3.9,8.16)
 \psline(4.1,8.1)(4.5,8)
 \psline(1.5,7)(1.9,6.867)
 \psline(2.1,6.8)(2.9,6.533)
 \psline(3.1,6.467)(3.9,6.2)
 \psline(4.1,6.133)(4.5,6)
 \psline(4.5,6)(4.6,5.9)
 \psline(4.75,5.75)(5,5.5)
 \psline(5,5.5)(5,1.833)
 \psline(5,9)(5,6.2)
 \psline(5,6.2)(4.1,5.12)
 \psline(3.9,4.96)(3.1,4.64)
 \psline(2.9,4.56)(2.1,4.24)
 \psline(1.9,4.16)(1.5,4)
 \psline(4.5,3)(4.1,2.867)
 \psline(3.9,2.8)(3.1,2.533)
 \psline(2.9,2.467)(2.1,2.2)
 \psline(1.9,2.133)(1.5,2)
 \psline(1.5,2)(1,1.833)
\end{pspicture}     }
 \rput(0,-8){%
 \psline(1,1.833)(1,1)
      \psline(2,1.833)(2,1)
      \psline(3,1.833)(3,1)
     \psline(5,1.833)(4,1)
     \psline(4,1.833)(4.4,1.4998)
     \psline(4.6,1.3332)(5,1) }
  \pscircle[linestyle=dotted](4.2,8){.5}
  \pscircle[linestyle=dotted](4.2,3){.5}
  \pscircle[linestyle=dotted](4.2,0){.5}
  \pscircle[linestyle=dotted](4.2,-5){.5}
\end{pspicture}
\hspace{3cm}
\begin{pspicture}(0,-7)(8,9)
\psdots(1,9)(2,9)(3,9)(4,9)(5,9)               %
       (5,-7)(1,-7)(2,-7)(3,-7)(4,-7)               %
 \psline(3.8,7.8)(1.5,7)
 \psline(1.5,4)(3.8,3.2)
 \psline(2,9)(2,7.25)
 \psline(2,7.05)(2,3.95)
 \psline(2,3.7)(2,1.833)
 \psline(3,9)(3,7.6)
 \psline(3,7.4)(3,3.6)
 \psline(3,3.4)(3,1.833)
 \psline(4,9)(4,8.4)           
 \psline(4,7.6)(4,3.4)         
 \psline(4,2.6)(4,1.833)
 \psline(1,9)(1.9,8.746)
 \psline(2.1,8.68)(2.9,8.45)
 \psline(3.1,8.39)(3.8,8.2)
 \psline(1.5,7)(1.9,6.867)
 \psline(2.1,6.8)(2.9,6.533)
 \psline(3.1,6.467)(3.9,6.2)
 \psline(4.1,6.133)(4.5,6)
 \psline(4.5,6)(5,5.5)
 \psline(5,5.5)(5,1.833)
 \psline(5,9)(5,6.2)
 \psline(5,6.2)(4.75,5.9)
 \psline(4.62,5.744)(4.1,5.12)
 \psline(3.9,4.96)(3.1,4.64)
 \psline(2.9,4.56)(2.1,4.24)
 \psline(1.9,4.16)(1.5,4)
 \psline(3.8,2.8)(3.1,2.533)
 \psline(2.9,2.467)(2.1,2.2)
 \psline(1.9,2.133)(1.5,2)
 \psline(1.5,2)(1,1.833)
      \psline(1,1.833)(1,1)
      \psline(2,1.833)(2,1)
      \psline(3,1.833)(3,1)
     \psline(4,1.833)(5,1)
     \psline(5,1.833)(4.6,1.4998)
     \psline(4.4,1.3332)(4,1)
 \psline(4,7.8)(4,7.6)
 \psline(3.8,8.2)(4,7.8)
 \psline(3.8,7.8)(3.84,7.92)
 \psline(3.935,8.205)(4,8.4)
  \psline(4,2.8)(4,2.6)
  \psline(3.8,3.2)(4,2.8)
  \psline(3.8,2.8)(3.84,2.92)
  \psline(3.935,3.205)(4,3.4)
   \rput(.12,.5) {%
  \begin{pspicture}(0,17)
 \psline(3.8,7.8)(1.5,7)
 \psline(1.5,4)(3.8,3.2)
 \psline(2,9)(2,7.25)
 \psline(2,7.05)(2,3.95)
 \psline(2,3.7)(2,1.833)
 \psline(3,9)(3,7.6)
 \psline(3,7.4)(3,3.6)
 \psline(3,3.4)(3,1.833)
\psline(4,9)(4,8.4)           
 \psline(4,7.6)(4,3.4)         
 \psline(4,2.6)(4,1.833)
 \psline(1,9)(1.9,8.746)
 \psline(2.1,8.68)(2.9,8.45)
 \psline(3.1,8.39)(3.8,8.2)
 \psline(1.5,7)(1.9,6.867)
 \psline(2.1,6.8)(2.9,6.533)
 \psline(3.1,6.467)(3.9,6.2)
 \psline(4.1,6.133)(4.5,6)
 \psline(4.5,6)(4.6,5.9)
 \psline(4.75,5.75)(5,5.5)
 \psline(5,5.5)(5,1.833)
 \psline(5,9)(5,6.2)
 \psline(5,6.2)(4.1,5.12)
 \psline(3.9,4.96)(3.1,4.64)
 \psline(2.9,4.56)(2.1,4.24)
 \psline(1.9,4.16)(1.5,4)
 \psline(3.8,2.8)(3.1,2.533)
 \psline(2.9,2.467)(2.1,2.2)
 \psline(1.9,2.133)(1.5,2)
 \psline(1.5,2)(1,1.833)
  \psline(4,7.8)(4,7.6)
 \psline(3.8,8.2)(4,7.8)
 \psline(3.8,7.8)(3.84,7.92)
 \psline(3.935,8.205)(4,8.4)
  \psline(4,2.8)(4,2.6)
  \psline(3.8,3.2)(4,2.8)
  \psline(3.8,2.8)(3.84,2.92)
  \psline(3.935,3.205)(4,3.4)
\end{pspicture}     }
 \rput(0,-8){%
 \psline(1,1.833)(1,1)
      \psline(2,1.833)(2,1)
      \psline(3,1.833)(3,1)
     \psline(5,1.833)(4,1)
     \psline(4,1.833)(4.4,1.4998)
     \psline(4.6,1.3332)(5,1) }
\end{pspicture}
\]

then isotope and perform $\mathcal{M}$ moves inside the dotted circles:\\

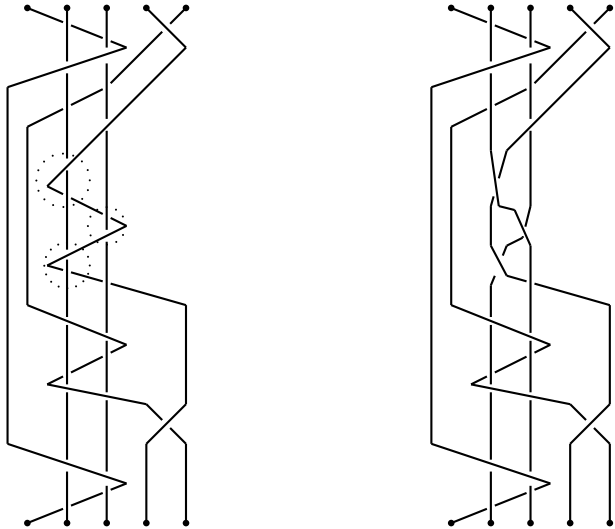
\[
\psset{unit=15pt}
\begin{pspicture}(-1,-2)(4,12)
 \psdots(0,11)(1,11)(2,11)(3,11)(4,11)               %
        (0,-2)(1,-2)(2,-2)(3,-2)(4,-2)               %
\psline(0,-2)(.9,-1.64) \psline(1.1,-1.56)(1.9,-1.24)%
\psline(2.1,-1.16)(2.5,-1)                           %
  \psline(4,3.5)(4,1)                                
 \psline(2.5,10)(-.5,9)                              
  \psline(3,11)(4,10)                                
     \psline(4,11)(3.6,10.6)                         
      \psline(4,1)(3,0)                              
      \psline(3,0)(3,-2)                             
       \psline(4,0)(4,-2)                            
\psline(3,1)(3.4,.6) \psline(3.6,.4)(4,0)            
\psline(-.5,0)(2.5,-1) \psline(-.5,9)(-.5,0)         %
\psline(2,-1.06)(2,-2)  \psline(1,-.7)(1,-2)         %
\psline(2,1.1)(2,-.7)                                %
\psline(0,11)(.9,10.64) \psline(1.1,10.56)(1.9,10.24)%
\psline(2.1,10.16)(2.5,10)                           %
\psline(1.9,8.95)(1.1,8.55) \psline(.9,8.45)(0,8)
   \psline(1,11)(1,9.65)  \psline(3.4,10.4)(2.1,9.1)
   \psline(2,11)(2,10)   \psline(4,10)(.5,6.5)
                \psline(1,9.35)(1,7.2)
                \psline(2,9.6)(2,8.2)
    \psline(2.5,5.5)(.5,4.5) \psline(1,6.9)(1,4.9) \psline(2,7.9)(2,5.4)
 \psline(0,3.5)(2.5,2.5)   \psline(0,8)(0,3.5)
 \psline(1,4.65)(1,3.2)   \psline(2,5.1)(2,2.8)
      \psline(.5,1.5)(3,1)
   \psline(2,2.6)(2,1.3)
   \psline(1,3)(1,1.5) \psline(1,1.3)(1,-.4)
   \psline(.5,1.5)(.9,1.7) \psline(1.1,1.8)(1.9,2.2)
   \psline(2.1,2.3)(2.5,2.5)
   \psline(.5,4.5)(.9,4.38571) \psline(1.1,4.32857)(1.9,4.1)
   \psline(2.1,4.04286)(4,3.5)
   \psline(.5,6.5)(.9,6.3) \psline(1.1,6.2)(1.9,5.8)
   \psline(2.1,5.7)(2.5,5.5)
 \pscircle[linestyle=dotted](.9,6.65){.7}
 \pscircle[linestyle=dotted](2,5.5){.5}
 \pscircle[linestyle=dotted](1,4.5){.6}
\end{pspicture}
\hspace{3cm}
\begin{pspicture}(-1,-2)(4,12)
 \psdots(0,11)(1,11)(2,11)(3,11)(4,11)               %
        (0,-2)(1,-2)(2,-2)(3,-2)(4,-2)               %
\psline(0,-2)(.9,-1.64) \psline(1.1,-1.56)(1.9,-1.24)%
\psline(2.1,-1.16)(2.5,-1)                           %
  \psline(4,3.5)(4,1)                                
 \psline(2.5,10)(-.5,9)                              
  \psline(3,11)(4,10)                                
     \psline(4,11)(3.6,10.6)                         
      \psline(4,1)(3,0)                              
      \psline(3,0)(3,-2)                             
       \psline(4,0)(4,-2)                            
\psline(3,1)(3.4,.6) \psline(3.6,.4)(4,0)            
\psline(-.5,0)(2.5,-1) \psline(-.5,9)(-.5,0)         %
\psline(2,-1.06)(2,-2)  \psline(1,-.7)(1,-2)         %
\psline(2,1.1)(2,-.7)                                %
\psline(0,11)(.9,10.64) \psline(1.1,10.56)(1.9,10.24)%
\psline(2.1,10.16)(2.5,10)                           %
\psline(1.9,8.95)(1.1,8.55) \psline(.9,8.45)(0,8)
   \psline(1,11)(1,9.65)  \psline(3.4,10.4)(2.1,9.1)
   \psline(2,11)(2,10)   \psline(4,10)(1.4,7.4)
                \psline(1,9.35)(1,7.4)
                \psline(2,9.6)(2,8.2)
    \psline(1.8,5.2)(1.4,5) \psline(1,6)(1,5) \psline(2,7.9)(2,6)
 \psline(0,3.5)(2.5,2.5)   \psline(0,8)(0,3.5)
 \psline(1,4)(1,3.2)   \psline(2,5)(2,2.8)
      \psline(.5,1.5)(3,1)
   \psline(2,2.6)(2,1.3)
   \psline(1,3)(1,1.5) \psline(1,1.3)(1,-.4)
   \psline(.5,1.5)(.9,1.7) \psline(1.1,1.8)(1.9,2.2)
   \psline(2.1,2.3)(2.5,2.5)
  \psline(1.4,4.24286)(1.9,4.1)
   \psline(2.1,4.04286)(4,3.5)
    \psline(1.2,6)(1.6,5.9)
  \psline(1,7.4)(1.2,6)  \psline(1.4,7.4)(1.2,6.7)
  \psline(1,6)(1.07,6.245)
  \psline(1.6,5.9)(2,5) \psline(1.81,5.24)(1.8,5.2)
  \psline(1.87,5.48)(2,6)
  \psline(1,5)(1.4,4.24286) \psline(1.4,5)(1.3,4.75)
  \psline(1.1,4.24)(1,4)
\end{pspicture}
\]
\newpage

Isotope and perform $\mathcal{M}$ moves:\\
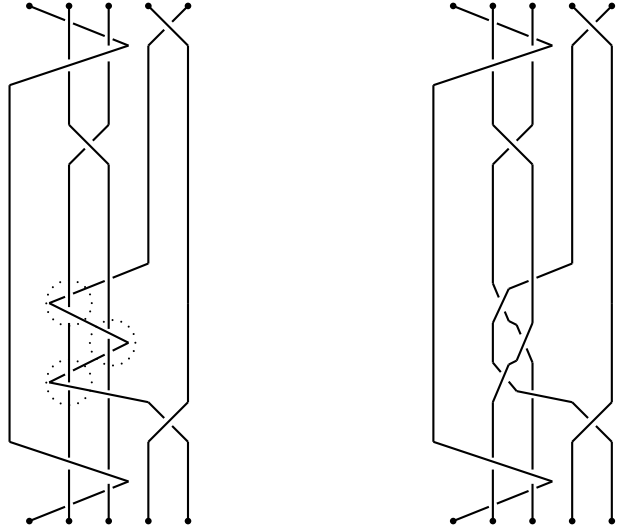
\[
\psset{unit=15pt}
\begin{pspicture}(-1,-2)(4,12)
 \psdots(0,11)(1,11)(2,11)(3,11)(4,11)               %
        (0,-2)(1,-2)(2,-2)(3,-2)(4,-2)               %
\psline(0,-2)(.9,-1.64) \psline(1.1,-1.56)(1.9,-1.24)%
\psline(2.1,-1.16)(2.5,-1)                           %
  \psline(4,3.5)(4,1)                                
 \psline(2.5,10)(-.5,9)                              
  \psline(3,11)(4,10)                                
     \psline(4,11)(3.6,10.6)                         
      \psline(4,1)(3,0)                              
      \psline(3,0)(3,-2)                             
       \psline(4,0)(4,-2)                            
\psline(3,1)(3.4,.6) \psline(3.6,.4)(4,0)            
\psline(-.5,0)(2.5,-1) \psline(-.5,9)(-.5,0)         %
\psline(2,-1.06)(2,-2)  \psline(1,-.7)(1,-2)         %
\psline(2,1.1)(2,-.7)                                %
\psline(0,11)(.9,10.64) \psline(1.1,10.56)(1.9,10.24)%
\psline(2.1,10.16)(2.5,10)                           %
 \psline(4,10)(4,3.5) \psline(3,10)(3,4.5)
    \psline(1,11)(1,9.65)  \psline(3.4,10.4)(3,10)
    \psline(2,11)(2,10)
                \psline(1,9.35)(1,8)
                \psline(2,9.6)(2,8)
 \psline(.5,3.5)(2.5,2.5)
 \psline(1,7)(1,3.4)   \psline(2,7)(2,2.85)
      \psline(.5,1.5)(3,1)
   \psline(2,2.6)(2,1.3)
   \psline(1,3)(1,1.5) \psline(1,1.3)(1,-.4)
   \psline(.5,1.5)(.9,1.7) \psline(1.1,1.8)(1.9,2.2)
   \psline(2.1,2.3)(2.5,2.5)
\psline(.5,3.5)(.9,3.66) \psline(1.1,3.74)(1.9,4.06)
\psline(2.1,4.14)(3,4.5)
  \psline(1,8)(2,7) \psline(1,7)(1.4,7.4) \psline(1.6,7.6)(2,8)
 \pscircle[linestyle=dotted](1,3.5){.6}
 \pscircle[linestyle=dotted](2.1,2.5){.6}
 \pscircle[linestyle=dotted](1,1.5){.6}
\end{pspicture}
\hspace{3cm}
\begin{pspicture}(-1,-2)(4,12)
 \psdots(0,11)(1,11)(2,11)(3,11)(4,11)               %
        (0,-2)(1,-2)(2,-2)(3,-2)(4,-2)               %
\psline(0,-2)(.9,-1.64) \psline(1.1,-1.56)(1.9,-1.24)%
\psline(2.1,-1.16)(2.5,-1)                           %
  \psline(4,3.5)(4,1)                                
 \psline(2.5,10)(-.5,9)                              
  \psline(3,11)(4,10)                                
     \psline(4,11)(3.6,10.6)                         
      \psline(4,1)(3,0)                              
      \psline(3,0)(3,-2)                             
       \psline(4,0)(4,-2)                            
\psline(3,1)(3.4,.6) \psline(3.6,.4)(4,0)            
\psline(-.5,0)(2.5,-1) \psline(-.5,9)(-.5,0)         %
\psline(2,-1.06)(2,-2)  \psline(1,-.7)(1,-2)         %
\psline(2,1.1)(2,-.7)                                %
\psline(0,11)(.9,10.64) \psline(1.1,10.56)(1.9,10.24)%
\psline(2.1,10.16)(2.5,10)                           %
 \psline(4,10)(4,3.5) \psline(3,10)(3,4.5)
    \psline(1,11)(1,9.65)  \psline(3.4,10.4)(3,10)
    \psline(2,11)(2,10)
                \psline(1,9.35)(1,8)
                \psline(2,9.6)(2,8)
 \psline(1.4,3.05)(1.6,2.95)
         \psline(1.4,3.05)(1.3,3.2875)
         \psline(1.15,3.64375)(1,4)
     \psline(1.6,2.95)(1.7,2.7125)
     \psline(1.85,2.35625)(2,2)
 \psline(1,7)(1,4)   \psline(2,7)(2,3)
      \psline(1.6,1.28)(3,1)
  \psline(1.6,1.28)(1.4,1.52)
  \psline(1.2,1.76)(1,2)
   \psline(2,2)(2,1.3)
   \psline(1,3)(1,2) \psline(1,1)(1,-.4)
    \psline(1.4,1.95)(1.6,2.05)
   \psline(1.6,2.05)(2,3)
     \psline(1.4,3.86)(1,3)
   \psline(1.4,1.95)(1,1)
 \psline(1.4,3.86)(1.9,4.06)
\psline(2.1,4.14)(3,4.5)
  \psline(1,8)(2,7) \psline(1,7)(1.4,7.4) \psline(1.6,7.6)(2,8)
\end{pspicture}
\]
\

And finally isotope to get id:\\
\[
\psset{unit=15pt}
\begin{pspicture}(-1,-2)(4,12)
 \psdots(0,11)(1,11)(2,11)(3,11)(4,11)               %
        (0,-2)(1,-2)(2,-2)(3,-2)(4,-2)               %
\psline(0,11)(0,-2) \psline(1,11)(1,-2) \psline(2,11)(2,-2)
\psline(3,11)(3,-2) \psline(4,11)(4,-2)
\end{pspicture}
\]
\end{proof}

\newpage

\begin{prop}\label{sq} Every square of the form
$$\xymatrix{ {\rho_1}\ar@{-}[d]_{x}\ar@{-}[r]^{y} &
{\rho_2}\ar@{-}[d]^{x}\\
{\rho_3}\ar@{-}[r]_{y} & {\rho_4} } $$
 in $\CC(\rho)$ bounds a $2$-cell.
\end{prop}

\begin{proof}Assume first that $x,y\neq \delta_4$. Then it suffices
to prove that $x$ and $y$ cannot be adjacent elementary braids.
Suppose then, without loss of generality, that $x=\B_i$ and
$y=\B_{i+1}$ for some $i\geq 2$. Since $\rho_1$ is moved by $\B_i$
its $i^{\text{th}}$ and $(i+1)^{\text{th}}$ monodromies are
different, and similarly its $(i+1)^{\text{th}}$ and
$(i+2)^{\text{th}}$ monodromies are different. Thus without loss
of generality
$$\rho_1=\dotsc,\underset{i}{(14)},\underset{i+1}{(23)},\underset{i+2}{(14)}
\dotsc$$ but then
\begin{align*} \rho_1\B_i\B_{i+1} &=\dotsc,\underset{i}{(23)},\underset{i+1}
{(14)},\underset{i+2}{(14)}\dotsc\\
&\neq\dotsc,\underset{i}{(14)},\underset{i+1}{(14)},\underset{i+2}{(23)}\dotsc\\
&=\rho_1\B_{i+1}\B_i
\end{align*}
a contradiction.\

If one of the $x$, $y$ is equal to $\delta_4$, say $x=\delta_4$,
one only needs to check the cases $y=\B_4$ and $y=\B_5$ since
$d_4$ doesn't meet any other ``generating'' interval.\\
If $y=\B_4$, apply Lemma \ref{Com}.\\
 Since $\rho_1$ is moved by $\delta_4$ it has either three
 monodromies equal to $(14)$ and one equal to $(23)$ or three
 monodromies equal to $(23)$ and one equal to $(14)$, in the spots
 $2$, $3$, $4$, $5$ and if $y=\B_5$ then the $5^{\text{th}}$ and $6^{\text{th}}$
monodromies of $\rho_1$ are different. But then applying $\B_5$ to
$\rho_1$ one sees that $\rho_2$ has an even number of monodromies
equal to $(14)$ in the spots $2$, $3$, $4$, $5$ and it is
therefore fixed by $\delta_4$, a contradiction.
\end{proof}

\subsection{Generators}\label{thegenerators} We use the same symbol to denote
an element of $\Au(\rho)$  and
its image in $\overline{\Au}(\rho)$.
 This subsection is devoted to the proof of the
following theorem:
\begin{thm}\label{main}
 For all even $n\geq 6$, $\overline{\Au}(\rho_{23}^n)$ is
generated by $\B_0$, $\B_2$, $\B_4$, $\B_5$,$\dotsc, \B_{n-2}$,
$\delta_4$ and, if $n\geq 8$, $\delta_6$, where,
$$\delta_6=[\delta_4]\B_5^{-1}\B_4^{-1}\B_3^{-1}\B_2^{-1}
\B_6^{-1}\B_5^{-1}\B_4^{-1}\B_3^{-1}.$$
\end{thm}

Actually the following stronger statement will be proved:

\begin{thm}\label{that}
\begin{itemize}\item[(a)] For each $n\geq 6$, $\overline{\Au}(\rho_{23}^n)$ is
generated by $\B_0$, $\B_2$, $\B_4$, $\B_5$,$\dotsc, \B_{n-2}$,
$\delta_4$ and, if $n\geq 8$, $\delta_6$.

\item[(b)] For each even $n\geq 6$, $\overline{\Au}(\rho_{2,3,n-1}^n)$ is
generated by $\B_0$, $\B_2$, $\B_4$,$\dots, \B_{n-3}$,
$[\B_{n-2}]\B_{n-3}^{-1}\dotsm\B_3^{-1}$ and if $n\geq8$,
$\delta_4$ and, if $n\geq 10$, $\delta_6$.

\end{itemize}
\end{thm}

The proof will be by induction on $n$.
  In the $n=6$ case is not difficult to actually give generators for the
 automorphism group:

 \begin{prop}\label{n=6}\begin{itemize}\item[(a)] $\Au(\rho_{23}^{6})$ is
 the smallest subgroup of $B_6$
 containing the following elements:
 $\B_0$, $\B_1^3$, $\B_2$, $\B_3^2$, $\B_4$, $\delta_4$,
 $[\B_1^3]\B_2^{-1}\B_3^{-1}$, $[\delta_4]\B_4^{-1}\B_3^{-1}$,
 $[\B_2^2]\B_3^{-1}$, $[\B_4^2]\B_3^{-1}$,
 $[\B_3^2]\B_2^{-1}\B_4^{-1}\B_3^{-1}$

 \item[(b)] $\Au(\tilde\rho_4^{6})$ is the smallest subgroup of $B_6$
 containing the following elements:
$\B_0$, $\B_1^3$, $\B_2$, $\B_3^2$, $\B_4^2$, $\delta_4^2$,
$[\B_2^2]\B_3^{-1}$,  $[\delta_4^2]\B_4^{-1}$,
$\B_4\delta_4\B_4\delta_4^{-1}$, $[\B_1^3]\B_2^{-1}\B_3^{-1}$,
$\delta_4\B_4\delta_4^{-1}\B_4^{-1}$,
$[\delta_4\B_4\delta_4^{-1}]\B_3^{-1}$,
 $[\B_3]\B_4^{-1}$, $[\B_4]\B_3^{-1}$.
 \end{itemize}
\end{prop}
\begin{proof}
\begin{itemize} \item[(a)] The graph in Figure \ref{six} describes
(the $1$-skeleton of) $\CC(\rho_{23}^{6})$, together with a
Schreier transversal drawn with wavy lines. Loops are not depicted
as, for example, the presence of a loop labelled by $\B_1^3$ at
the vertex $\rho_{23}^6$ can be deduced by the absence of an edge
labelled by $\B_1^3$ at that vertex.

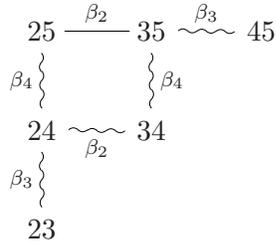
\begin{figure}[htp]
$$
\xymatrix{
 {25}\ar@{~}[d]_{\beta_4}\ar@{-}[r]^{\beta_2}
 &{35}\ar@{~}[d]^{\beta_4}\ar@{~}[r]^{\beta_3}
 &{45}\\
 {24}\ar@{~}[d]_{\beta_3}\ar@{~}[r]_{\beta_2}
 &{34}\\
 {23}         }
$$
\caption{$\CC(\rho_{23}^n)$ with a maximal tree}\label{six}
\end{figure}

According to the Reidemeister-Schreier method (see for example \cite{ZVC})
$\Au(\rho_{23}^n)$ is generated by
$g_is_j\overline{g_is_j}$ where $g_i$ runs through the transversal, $s_j$
runs through the generators and $\bar x$
denotes the element of the transversal which ends at the endpoint of $x$.
Then one gets the following generators:

\begin{itemize}
\item[$g_i=$id]

$$\B_0,\ \B_1^3,\ \B_2,\ \text{id},\ \B_4,\ \delta_4$$

\item[$g_i=\B_3$]

$$\begin{array}{ll}
\B_0,\ \B_1^3,&\\
\B_3\B_2(\overline{\B_3\B_2})^{-1}&=\text{id},\\
\B_3\B_3(\overline{\B_3\B_3})^{-1}&=\B_3^2,\\
\B_3\B_4(\overline{\B_3\B_4})^{-1}&=\text{id},\\
\B_3\delta_4(\overline{\B_3\delta_4})^{-1}&=\B_3\delta_4\B_3^{-1}\\
&=\delta_4.
\end{array}$$

\item[$g_i=\B_3\B_4$]

$$\begin{array}{ll}
\B_0,\ \B_1^3,&\\
\B_3\B_4\B_2(\overline{\B_3\B_4\B_2})^{-1}
&=\B_3\B_4\B_2\B_4^{-1}\B_2^{-1}\B_3^{-1}\\
&=\text{id},\\
\B_3\B_4\B_3(\overline{\B_3\B_4\B_3})^{-1}
&=\B_3\B_4\B_3\B_4^{-1}\B_3^{-1}\\
&=\B_4,\\
\B_3\B_4\B_4(\overline{\B_3\B_4\B_4})^{-1}&=\B_3\B_4^2\B_3^{-1},\\
\B_3\B_4\delta_4(\overline{\B_3\B_4\delta_4})^{-1}
&=\B_3\B_4\delta_4\B_4^{-1}\B_3^{-1}.
\end{array}$$

\item[$g_i=\B_3\B_2$]

$$\begin{array}{ll}
\B_0,&\\
\B_3\B_2\B_1^3(\overline{\B_3\B_2\B_1^3})^{-1}
&=\B_3\B_2\B_1^3\B_2^{-1}\B_3^{-1},\\
\B_3\B_2\B_2(\overline{\B_3\B_2\B_2})^{-1}
&=\B_3\B_2^2\B_3^{-1},\\
\B_3\B_2\B_3(\overline{\B_3\B_2\B_3})^{-1}&=\B_2,\\
\B_3\B_2\B_4(\overline{\B_3\B_2\B_4})^{-1}&=\text{id},\\
\B_3\B_2\delta_4(\overline{\B_3\B_2\delta_4})^{-1}&=\delta_4.
\end{array}$$

\item[$g_i=\B_3\B_2\B_4$]

$$\begin{array}{ll}
\B_0,&\\
\B_3\B_2\B_4\B_1^3(\overline{\B_3\B_2\B_4\B_1^3})^{-1}
&=\B_3\B_2\B_4\B_1^3\B_4^{-1}\B_2^{-1}\B_3^{-1}\\
&=\B_3\B_2\B_1^3\B_2^{-1}\B_3^{-1},\\
\B_3\B_2\B_4\B_2(\overline{\B_3\B_2\B_4\B_2})^{-1}
&=\B_3\B_2\B_4\B_2\B_4^{-1}\B_3^{-1}\\
&=\B_3\B_2^2\B_3^{-1},\\
\B_3\B_2\B_4\B_3(\overline{\B_3\B_2\B_4\B_3})^{-1}&=\text{id},\\
\B_3\B_2\B_4\B_4(\overline{\B_3\B_2\B_4\B_4})^{-1}
&=\B_3\B_2\B_4^2\B_2^{-1}\B_3^{-1}\\
&=\B_3\B_4^2\B_3^{-1},\\
\B_3\B_2\B_4\delta_4(\overline{\B_3\B_2\B_4\delta_4})^{-1}
&=\B_3\B_2\B_4\delta_4\B_4^{-1}\B_2^{-1}\B_3^{-1}\\
&=\B_3\B_4\delta_4\B_4^{-1}\B_3^{-1}.
\end{array}$$

\item[$g_i=\B_3\B_2\B_4\B_3$]

$$\begin{array}{ll}
\B_0,&\\
\B_3\B_2\B_4\B_3\B_1^3(\overline{\B_3\B_2\B_4\B_3\B_1^3})^{-1}
&=\B_3\B_2\B_4\B_3\B_1^3\B_3^{-1}\B_4^{-1}\B_2^{-1}\B_3^{-1}\\
&=\B_3\B_2\B_1^3\B_2^{-1}\B_3^{-1},\\
\B_3\B_2\B_4\B_3\B_2(\overline{\B_3\B_2\B_4\B_3\B_2})^{-1}
&=\B_3\B_2\B_4\B_3\B_2\B_3^{-1}\B_4^{-1}\B_2^{-1}\B_3^{-1}\\
&=\B_3\B_4\B_2\B_3\B_2\B_3^{-1}\B_2^{-1}\B_4^{-1}\B_3^{-1}\\
&=\B_3\B_4\B_3\B_4^{-1}\B_3^{-1}\\
&=\B_4,\\
\B_3\B_2\B_4\B_3\B_3(\overline{\B_3\B_2\B_4\B_3\B_3})^{-1}
&=\B_3\B_2\B_4\B_3^2\B_4^{-1}\B_2^{-1}\B_3^{-1},\\
\B_3\B_2\B_4\B_3\B_4(\overline{\B_3\B_2\B_4\B_3\B_4})^{-1}
&=\B_3\B_2\B_4\B_3\B_4\B_3^{-1}\B_4^{-1}\B_2^{-1}\B_3^{-1}\\
&=\B_3\B_2\B_3\B_2^{-1}\B_3^{-1}\\
&=\B_2,\\
\B_3\B_2\B_4\B_3\delta_4(\overline{\B_3\B_2\B_4\B_3\delta_4})^{-1}
&=\B_3\B_2\B_4\B_3\delta_4\B_3^{-1}\B_4^{-1}\B_2^{-1}\B_3^{-1}\\
&=\B_3\B_4\delta_4\B_4^{-1}\B_3^{-1}.
\end{array}$$
\end{itemize}

\item[(b)] Similar calculations using the graph in Figure~\ref{sixtwo} which
represents the $1$-skeleton of $\CC(\tilde\rho_4^{6})$ together with a Schreier
transversal give the result.

\begin{figure}[htp]
$$
\xymatrix{
 {\tilde 2}\ar@{~}[d]_{\delta_4}\ar@{~}[r]^{\beta_2}
 &{\tilde 3}\ar@{~}[d]_{\delta_4}\ar@{~}[r]^{\beta_3}
 &{\tilde 4}\ar@{~}[d]_{\delta_4}\ar@{~}[r]^{\beta_4}
 &{\tilde 5}\ar@{~}[d]^{\delta_4}\\
 {2}\ar@{-}[r]_{\beta_2}
 &{3}\ar@{-}[r]_{\beta_3}
 &{4}\ar@{-}[r]_{\beta_4}
 &{5}            }
$$
\caption{$\CC(\tilde\rho_4^{6})$ with a maximal
tree}\label{sixtwo}
\end{figure}
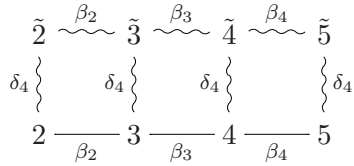
\end{itemize}
\end{proof}

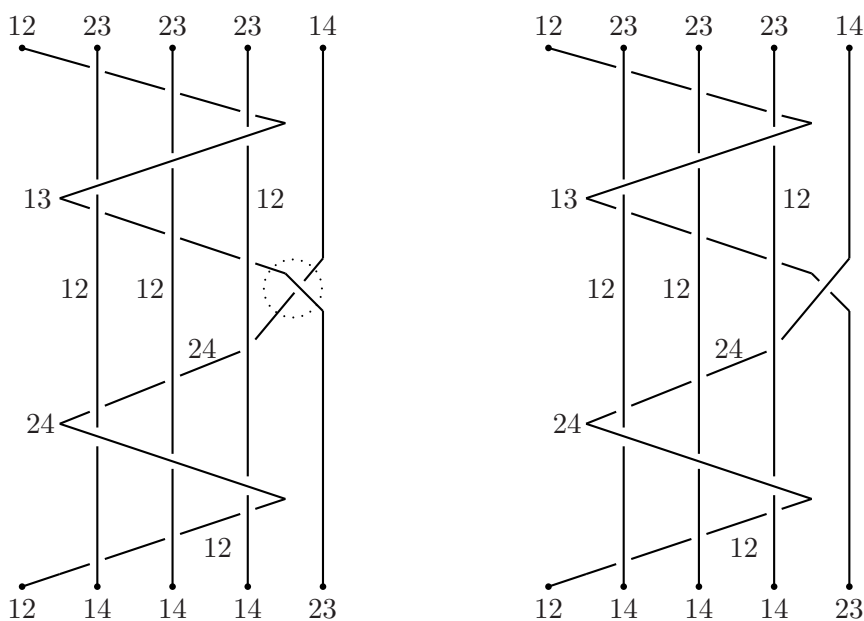
\begin{figure}[htp]
\[
\begin{pspicture}(1,1.5)(5,9.5)
 \psline(4.5,8)(1.5,7)
 \psline(1.5,4)(4.5,3)
 \psline(2,9)(2,7.25)
 \psline(2,7.05)(2,3.95)
 \psline(2,3.7)(2,1.833)
 \psline(3,9)(3,7.6)
 \psline(3,7.4)(3,3.6)
 \psline(3,3.4)(3,1.833)
 \psline(4,9)(4,7.95)
 \psline(4,7.7)(4,3.3)
 \psline(4,3.05)(4,1.833)
 \psline(1,9)(1.9,8.746)
 \psline(2.1,8.68)(2.9,8.45)
 \psline(3.1,8.39)(3.9,8.16)
 \psline(4.1,8.1)(4.5,8)
 \psline(1.5,7)(1.9,6.867)
 \psline(2.1,6.8)(2.9,6.533)
 \psline(3.1,6.467)(3.9,6.2)
 \psline(4.1,6.133)(4.5,6)
 \psline(4.5,6)(5,5.5)
 \psline(5,5.5)(5,1.833)
 \psline(5,9)(5,6.2)
 \psline(5,6.2)(4.75,5.9)
 \psline(4.62,5.744)(4.1,5.12)
 \psline(3.9,4.96)(3.1,4.64)
 \psline(2.9,4.56)(2.1,4.24)
 \psline(1.9,4.16)(1.5,4)
 \psline(4.5,3)(4.1,2.867)
 \psline(3.9,2.8)(3.1,2.533)
 \psline(2.9,2.467)(2.1,2.2)
 \psline(1.9,2.133)(1.5,2)
 \psline(1.5,2)(1,1.833)
 \psdots(1,9)(2,9)(3,9)(4,9)(5,9)
        (1,1.833)(2,1.833)(3,1.833)(4,1.833)(5,1.833)
 \rput(1,9.3){12}
 \rput(2,9.3){23}
 \rput(3,9.3){23}
 \rput(4,9.3){23}
 \rput(5,9.3){14}
 \rput(1,1.533){12}
 \rput(2,1.533){14}
 \rput(3,1.533){14}
 \rput(4,1.533){14}
 \rput(5,1.533){23}
 \rput(1.2,7){13}
 \rput(1.25,4){24}
 \rput(3.4,5){24}
 \rput(3.6,2.35){12}
 \rput(4.3,7){12}
 \rput(1.7,5.8){12}
 \rput(2.7,5.8){12}
 \pscircle[linestyle=dotted](4.6,5.8){.4}
\end{pspicture}
\hspace{3cm}
\begin{pspicture}(1,1.5)(5,9.5)
 \psline(4.5,8)(1.5,7)
 \psline(1.5,4)(4.5,3)
 \psline(2,9)(2,7.25)
 \psline(2,7.05)(2,3.95)
 \psline(2,3.7)(2,1.833)
 \psline(3,9)(3,7.6)
 \psline(3,7.4)(3,3.6)
 \psline(3,3.4)(3,1.833)
 \psline(4,9)(4,7.95)
 \psline(4,7.7)(4,3.3)
 \psline(4,3.05)(4,1.833)
 \psline(1,9)(1.9,8.746)
 \psline(2.1,8.68)(2.9,8.45)
 \psline(3.1,8.39)(3.9,8.16)
 \psline(4.1,8.1)(4.5,8)
 \psline(1.5,7)(1.9,6.867)
 \psline(2.1,6.8)(2.9,6.533)
 \psline(3.1,6.467)(3.9,6.2)
 \psline(4.1,6.133)(4.5,6)
 \psline(4.5,6)(4.6,5.9)
 \psline(4.75,5.75)(5,5.5)
 \psline(5,5.5)(5,1.833)
 \psline(5,9)(5,6.2)
 \psline(5,6.2)(4.1,5.12)
 \psline(3.9,4.96)(3.1,4.64)
 \psline(2.9,4.56)(2.1,4.24)
 \psline(1.9,4.16)(1.5,4)
 \psline(4.5,3)(4.1,2.867)
 \psline(3.9,2.8)(3.1,2.533)
 \psline(2.9,2.467)(2.1,2.2)
 \psline(1.9,2.133)(1.5,2)
 \psline(1.5,2)(1,1.833)
 \psdots(1,9)(2,9)(3,9)(4,9)(5,9)
        (1,1.833)(2,1.833)(3,1.833)(4,1.833)(5,1.833)
 \rput(1,9.3){12}
 \rput(2,9.3){23}
 \rput(3,9.3){23}
 \rput(4,9.3){23}
 \rput(5,9.3){14}
 \rput(1,1.533){12}
 \rput(2,1.533){14}
 \rput(3,1.533){14}
 \rput(4,1.533){14}
 \rput(5,1.533){23}
 \rput(1.2,7){13}
 \rput(1.25,4){24}
 \rput(3.4,5){24}
 \rput(3.6,2.35){12}
 \rput(4.3,7){12}
 \rput(1.7,5.8){12}
 \rput(2.7,5.8){12}
\end{pspicture}
\]
\caption{A $\mathcal{P}$ move (inside the dotted circle) turns
$\delta_4$ into $\delta_4^{-1}$ }\label{label}
\end{figure}

\begin{cor}  \begin{itemize} \item[(a)]
$\overline{\Au}(\rho_{23}^{6})$ is generated by $\B_0$, $\B_2$,
$\B_4$ and $\delta_4$.

\item[(b)]$\overline{\Au}(\rho_2^{6})$ is generated by $\B_0$, $\B_3$ and
$[\B_4]\B_3^{-1}$.

\end{itemize}
\end{cor}
\begin{proof} One just needs to check that the generators given in Proposition \ref{n=6}
 are equivalent $mod N$ to the generators
given above. Given Lemma \ref{Com} the only slightly non trivial
checks is to check that for any element of $\CC(\tilde\rho_5^{6})$
we have $\delta_4^2\equiv\text{id}\pmod{N}$ and
$[\B_4]\B_3^{-1}\equiv[\B_3]\B_4^{-1}\pmod{N}$. For the first look
at Figure \ref{label}, where it is shown that
$\delta_4\equiv\delta_4^{-1}\pmod{N}$ for $\rho_5^6$. For the
second look at Figure~\ref{Nolabel} where a more general statement
is shown.
\end{proof}

Next we do the $n=7$ case. The $1$-skeleton of $\CC(\rho_{23}^7)$
is shown in the Figure~\ref{seven}. Again loops are not shown.

\begin{figure}[htp]
$$
\xymatrix{
 & {\tilde 2}\ar@{-}[dl]_{\delta_4}\ar@{-}[r]^{\beta_2}
 & {\tilde 3}\ar@{-}[dl]_{\delta_4}\ar@{-}[r]^{\beta_3}
 & {\tilde 4}\ar@{-}[dl]_{\delta_4}\ar@{-}[r]^{\beta_4}
 & {\tilde 5}\ar@{-}[dl]_{\delta_4}\ar@{-}[r]^{\beta_5}
 & {\tilde 6}\\
 {26}\ar@{-}[d]_{\beta_5}\ar@{-}[r]^{\beta_2}
 &{36}\ar@{-}[d]_{\beta_5}\ar@{-}[r]^{\beta_3}
 &{46}\ar@{-}[d]_{\beta_5}\ar@{-}[r]^{\beta_4}
 &{56}\\
 {25}\ar@{-}[d]_{\beta_4}\ar@{-}[r]^{\beta_2}
 &{35}\ar@{-}[d]_{\beta_4}\ar@{-}[r]^{\beta_3}
 &{45}\\
 {24}\ar@{-}[d]_{\beta_3}\ar@{-}[r]^{\beta_2}
 &{34}\\
 {23}                }
 $$
 \caption{$\CC(\rho_{23}^7)$}\label{seven}
 \end{figure}

\begin{prop}\label{n=7} $\overline{\Au}(\rho_{23}^7)$ is generated by $\B_0$,
$\B_2$, $\B_4$, $\B_5$, and $\delta_4$.
\end{prop}

\begin{proof} First notice that $\overline{\Au}
(\rho_{23}^7)=\pi_1(\CC(\rho_{23}^7),\rho_{23}^7)$ is generated by
{\em lassoes}, i.e. elements of the form $wxw^{-1}$ where, $w$ is
a path not containing loops (that is, $w$ is a subtree) starting
at $\rho_{23}^7$ and $x$ is a loop at the endpoint of $w$.
Furthermore the tail $w$ is not important: two lassoes with the
same head (at the same vertex with the same label) represent the
same element.\

  This is clear from Figure \ref{seven} and Proposition \ref{sq}. \\

 By Proposition \ref{n=6} one only needs to check  lassoes
 with
 \begin{itemize}
 \item[(a)] heads at $\tilde\rho_4^7$\\
The heads are labelled by $\B_0$ or $\B_2$. Using the path
$w=\B_3\B_4\B_5\B_2\B_3\delta_4$ one gets:
$$w\B_0w^{-1}=\B_0$$
\begin{align*} w\B_2w^{-1} &=
[\B_2]\delta_4^{-1}\B_3^{-1}\B_2^{_1}\B_5^{-1}\B_4^{-1}\B_3^{-1}\\
 &=[\B_2]\B_3^{-1}\B_2^{_1}\B_5^{-1}\B_4^{-1}\B_3^{-1}\\
 &=[\B_3]B_5^{-1}\B_4^{-1}\B_3^{-1}\\
 &=[\B_3]\B_4^{-1}\B_3^{-1}\\
 &=\B_4
\end{align*}

\item[(b)] The lasso with head at $\tilde\rho_3^6$ labelled
by $\B_4$ (corresponding to $[\B_4]\B_3^{-1}$.\\

Using the path $w=\B_3\B_4\B_5\B_2\delta_4$ one gets:

\begin{align*} w\B_4w^{-1} &=
[\B_4]\delta_4^{-1}\B_2^{_1}\B_5^{-1}\B_4^{-1}\B_3^{-1}\\
 &=[\B_4]\B_5^{-1}\B_4^{-1}\B_3^{-1}\\
 &=[\B_5]\B_3^{-1}\\
  &=\B_5
\end{align*}

 \item[(c)] heads labelled by $\B_5$\\
In this case the head is either in $i_{(23)}(\CC(\rho_{23}^6))$
 or in $i_{(14)}(\CC(\rho_2^6))$.\\
 In the first case its obvious that both of them are equal to
 $\B_5$ since one can choose the path $w$ to lie entirely inside
 $i_{(23)}(\CC(\rho_{23}^6))$ and therefore all it's edges to
 commute with $\B_5$.\\
In the second case one only needs to check those with head at
$\rho_5^7$ and $\tilde\rho_2^7$.\\

For the first use the path $w=\B_3\B_2\B_4\B_3\B_5\B_4$ to get:
\begin{align*} w\B_5w^{-1} &=
[\B_5]\B_4^{-1}\B_5^{-1}\B_3^{-1}\B_4^{-1}\B_2^{-1}\B_3^{-1}\\
 & =[\B_4]\B_3^{-1}\B_4^{-1}\B_2^{-1}\B_3^{-1}\\
   & =[\B_3]\B_2^{-1}\B_3^{-1}\\
   & =\B2.
   \end{align*}

For the second use the path $w=\B_3\B_4\B_5\delta_4$ to get:

\begin{align*} w\B_5w^{-1} &=
[\B_5]\delta_4^{-1}\B_5^{-1}\B_4^{-1}\B_3^{-1}\\
&=[\delta_4]\B_4^{-1}\B_3^{-1}\\
&=\delta_4
\end{align*}

 \item[(d)] heads at $\tilde\rho_6^7$.\\

 Those have heads labelled by $\B_2$, or $\B_4$, or $\delta_4$.\\

 For the first one use any path to see that it is equal to a lasso with head at
 $\tilde\rho_5^7$, labelled by $\B_2$. The later has been checked
 in case a).\\

 For the second one use a path that ends with $\B_4\B_5$ say
 $w=w'\B_4\B_5$ to get:
 $$w'\B_4\B_5\B_4\B_5^{-1}\B_4^{-1}(w')^{-1}=w'\B_5(w')^{-1}$$
 which has head at $\tilde\rho_5^7$ and therefore has been checked
 in b). (It turns out to be $\delta_4$).\\

 Finally for the third one use a path of the form $w=w^{\prime}\delta_4\B_5$ to get:
 $$w'\delta_4\B_5\delta_4\B_5^{-1}\delta_4^{-1}
 (w')^{-1}=w'\B_5(w')^{-1}$$
 which has head at $\rho_5^7$ and therefore has been checked
 in b). (It turns out to be $\B_2$).
\end{itemize}
\end{proof}

\begin{lem}\label{vanK} The fundamental group of $\CC_\bullet^n$ is generated by
 lassoes. Furthermore the tail of a lasso is not important: two
lassoes based at the same vertex and having the same head (at the
same vertex labelled by the same element) are equal.

\end{lem}

\begin{proof}  Attach $2$-cells to kill all loops. It suffices to prove that
the resulting complex $\bar\CC_{\bullet}^n$ is simply connected.
By induction
 on the number of branch values $n$: It is true for $n=6,7$. Assume it is true for
 all smaller values of $n$. Then by Lemmas \ref{dec} and \ref{sq} $\bar\CC_{\text{id}}^n$  is
 the union of two simply  connected complexes glued along a connected subcomplex and
 therefore by Van Kampen's  theorem it is simply connected. The second statement follows
  from the fact that the intersection of the two  pieces is actually simply connected.\

The proof for $\CC_{(23)}^n$ and $\CC_{(14)(23)}^n$ is the same
after observing that the extra vertices are glued along isolated
edges.

\end{proof}

\begin{figure}[htp]
\[
\begin{pspicture}(0,-.9)(10,1.1)
 \psdots[dotstyle=*,dotscale=1](0,0)(1,0)(2,0)(3,0)(5.6,0)(6.6,0)(7.6,0)(8.6,0)(10.2,0)
 \pscurve
 (1,0)(3.5,.6)(6.5,.5)(7.8,0)(6.6,-.4)(5.2,0)(3.5,.3)(1.5,0)(1,-.4)(.5,0)(3.5,.8)(7.5,.6)(8.6,0)
 \psdots[dotstyle=*,dotscale=.4](9.1,0)(9.3,0)(9.5,0)
 \psdots[dotstyle=*,dotscale=.4](3.6,0)(3.8,0)(4.1,0)
 \rput(0,-.8){0}
 \rput(1,-.8){1}
 \rput(2,-.8){2}
 \rput(3,-.8){3}
 \rput(5.6,-.8){$k-2$}
 \rput(6.6,-.8){$k-1$}
 \rput(7.6,-.8){$k$}
 \rput(8.6,-.8){$k+1$}
 \rput(10.2,-.8){$n-1$}
\end{pspicture}
\]
\caption{The interval $d_k$}\label{DK}
\end{figure}
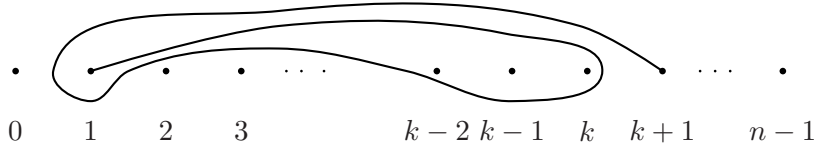

\begin{defn}  For even k, $d_k$ is the interval shown in
Figure~\ref{DK}
 and $\delta_k$ is the rotation around $d_k$.
\end{defn}

\begin{lem}\label{deltakey} Let $G$ be the subgroup of $B_n$ generated by the elements
listed in Theorem \ref{that}. Then $\delta_k$ is in $G$ for all
$k$.
\end{lem}
\begin{proof} It is clear that $\delta_k$
 is in $L(n)$. Now notice that for $k\geq 6$, $\delta_k$ fixes
pointwise a small disc containing $A_2$ and $A_3$ but no other
 branch value. Connect this disc to $\partial D^2$ by an arc that
 misses the $A_i$'s, the $\alpha_i$'s and the
$x_i$'s. Remove the disc and the arc to get a {\em three-sheeted}
covering of a disk which is equivalent to $\rho(n-2)$ via an
equivalence that sends $d_6$ to $d_4$. Apply Theorem~\ref{bw1} to
this covering and then glue the disk and the arc back in to get an
expression of $\delta_k$ as a word in the generators of $G$.
\end{proof}

\begin{proof}[Proof of \ref{that}] The cases $n=6$, $n=7$  have
been proved already. Let $n\geq 8$ and assume the theorem proven
for all smaller values of $n$.\

By \ref{vanK} suffices to prove that all lassoes are in $G$ where
where $G$ is the subgroup of $\overline{\Au}$
generated by the listed elements.\\

 If $n$ is even then (at the level of $0$-skeletons)

 \begin{align*}
\CC(\rho_{23}^n) & = \CC_{\text{id}}^n\\ & =
i_{(23)}(\CC_{(23)}^{n-1}) \bigsqcup i_{(14)}(\CC_{(14)}^{n-1})
\end{align*}
 and the lassoes fall in to two cases:
\begin{description}
\item[Case a] Lassoes with heads at
$i_{(23)}(\CC_{(23)}^{n-1})$.
\end{description}

 In this case by the induction hypothesis, one needs only to check lassoes with head
labelled by $\B_{n-2}$ . After writing
$$\CC_{(23)}^{n-1}=i_{(23)}(\CC_{\text{id}}^{n-2}
\sqcup\{\tilde\rho_\emptyset^{n-2}\})\bigsqcup
i_{(14)}(\CC_{(14)(23)}^{n-2})$$
 and observing that for lassoes in
 $i_{(23)}i_{(23)}(\CC_{\text{id}}^{n-2})$ the tail can be chosen
 to commute with $\B_{n-2}$, one is left with checking the
 lasso with head at $\tilde\rho_{n-2,n-1}^n$ labelled by
 $\B_{n-2}$.\\

 If $n=8$ notice that the path $w=\B_5\delta_4\B_6\B_5$ connects
 $\rho_{56}^8$ to $\tilde\rho_{56}^8$ and

 \begin{align*}
w\B_6w^{-1}&=\B_5\delta_4\B_6\B_5\B_6\B_5^{-1}\B_6^{-1}\delta_4^{-1}\B_5^{-1}\\
 &=\B_5\delta_4\B_5\delta_4^{-1}\B_5^{-1}\\
 &=\delta_4
\end{align*}

 so that the lasso with head at $\tilde\rho_{56}^8$
 labelled by $\B_6$ is equal to the lasso with head at $\rho_{56}^8$
 labelled by $\delta_4$. The later is by definition equal to
 $\delta_6$, since the path $\B_3\B_4\B_5\B_6\B_2\B_3\B_4\B_5$ connects $\rho_{23}^8$
 to $\rho_{56}^8$.\\

If $n\geq 10$ notice that the path
$w=\delta_{n-2}\B_{n-4}\B_{n-5}\B_{n-3}\B_{n-4}\B_{n-2}\B_{n-3}$
connects $\tilde\rho_{n-5,n-3,n-2,n-1}^n$ to
$\tilde\rho_{n-2,n-1}$ and
\begin{align*}
w\B_{n-2}w^{-1}&=
  [\B_{n-2}]\B_{n-3}^{-1}\B_{n-2}^{-1}\B_{n-4}^{-1}\B_{n-3}^{-1}\B_{n-5}^{-1}
    \B_{n-4}^{-1}\delta_{n-2}^{-1}\\
 &=[\B_{n-3}]\B_{n-4}^{-1}\B_{n-3}^{-1}\B_{n-5}^{-1}
    \B_{n-4}^{-1}\delta_{n-2}^{-1}\\
  &=[\B_{n-4}]\B_{n-5}^{-1}
    \B_{n-4}^{-1}\delta_{n-2}^{-1}\\
 &=[\B_{n-5}]\delta_{n-2}^{-1}
\end{align*}

and therefore by the induction hypothesis (since
$\tilde\rho_{n-5,n-3,n-2,n-1}^n\in i_{(23)}(\CC_{(23)}^{n-1})$)
the lasso with head at $\tilde\rho_{n-2,n-1}$ labelled by
$\B_{n-2}$ is indeed in $G$.

\begin{description}
\item[Case b] Lassoes with heads at $i_{(14)}(\CC_{(14)}^{n-1})$.
\end{description}

 By symmetry and the previous case one needs to check only
lassoes with head at $\tilde\rho_{23}^n$ and when $n\geq10$ the
lasso with head at $\tilde\rho_{56}^n$ labelled by
$\B_{n-2}$ (corresponding to $\delta_6$).\\

For the lassoes based at $\tilde\rho_{23}^n$ and head labelled by
$\B_2$, $\B_4$, $\B_5,\ldots$, $\B_{n-2}$ notice that the path
$w=\B_{n-3}\B_{n-4}\dotsm\B_2\B_{n-2}\B_{n-3}\dotsm\B_3$ connects
$\tilde\rho_{n-2,n-1}^n$ to $\tilde\rho_{23}^n$ and

\begin{align*} w\B_2w^{-1}&=
             [\B_2]\B_3^{-1}\B_4^{-1}\dotsm\B_{n-3}^{-1}\B_{n-2}^{-1}\B_2^{-1}
              \B_2^{-1}\B_3^{-1}\dotsm\B_{n-4}^{-1}\B_{n-3}^{-1}\\
               &=[\B_2]\B_3^{-1}\B_2^{-1}\B_4^{-1}\dotsm\B_{n-3}^{-1}\B_{n-2}^{-1}
                \B_2^{-1}\B_3^{-1}\dotsm\B_{n-4}^{-1}\B_{n-3}^{-1}\\
                &=[\B_3]\B_4^{-1}\dotsm\B_{n-3}^{-1}\B_{n-2}^{-1}
                  \B_2^{-1}\B_3^{-1}\dotsm\B_{n-4}^{-1}\B_{n-3}^{-1}\\
                  &=\cdots\\
                  &=[\B_{n-3}]\B_{n-2}^{-1}\B_{n-3}^{-1}\\
                           &=\B_{n-2}.
\end{align*}

 and (via similar calculations), for $i=4,\ldots,n-2$
$$w\B_iw^{-1}=\B_{i-2}.$$
 Therefore by the previous case all of these lassoes belong to
 $G$.\\

 For the lasso with head at $\tilde\rho_{23}^n$ labelled by
 $\delta_4$ notice that the path $w=\B_{n-2}\B_{n-3}\dotsm\B_7\delta_4\B_5\B_4\B_3$
   connects $\tilde\rho_{3,4,5,n-1}^n$ to $\tilde\rho_{23}^n$ and
 \begin{align*} w\delta_4w^{-1}&=[\delta_4]\B_3^{-1}\B_4^{-1}\B_5^{-1}
                    \delta_4^{-1}\B_7^{-1}\dotsm\B_{n-3}^{-1}\B_{n-2}^{-1}\\
                      &=[\delta_4]\B_5^{-1}\delta_4^{-1}\\
                      &=\B_5
  \end{align*}
and therefore since $\tilde\rho_{3,4,5,n-1}^n\in
i_{(23)}(\CC_{(23)}^{n-1}$ by case a) this lasso is in
$G$.\\

For the lasso with head at $\tilde\rho_{56}^n$ labelled by
$\B_{n-2}$ notice that the path $w=\B_{n-2}\B_{n-3}\dotsm\B_7$
connects $\tilde\rho_{5,n-1}^n$ to $\tilde\rho_{56}^n$ and
\begin{align*} w\B_{n-2}w^{-1}&=
        [\B_{n-2}]\B_7^{-1}\dotsm\B_{n-3}^{-1}\B_{n-2}^{-1}\\
             &=[\B_{n-2}]\B_{n-3}^{-1}\B_{n-2}^{-1}\\
                         &=\B_{n-3}.
\end{align*}
and therefore since $\tilde\rho_{5,n-1}^n\in
i_{(23)}(\CC_{(23)}^{n-1}$ by case a) this lasso is in
$G$.\\
 This concludes the proof for $\overline{\Au}(\rho_{23}^n)$ when $n$
is even.\\

Continuing with $n$ even:
\begin{align*} \CC(\rho_{2,3,n-1}^n)&=\CC_{(23)}^n\\
                  &=i_{(23)}(\CC_{(14)}^{n-1}\sqcup \{\rho_\emptyset^{n-1}\})
                  \bigsqcup i_{(14)}(\CC_{(23)}^{n-1}\sqcup \{\tilde\rho_\emptyset^{n-1}\})
\end{align*}

and the lassoes fall in to four cases:
\begin{description}
\item[Case a] Lassoes with heads at
$i_{(14)}(\CC_{(23)}^{n-1})$.
\end{description}

In this case by the induction hypothesis one only needs to check
lassoes with
 head labelled by $\B_{n-2}$. But those actually have heads lying inside
 $i_{(14)}i_{(14)}(\CC_{(14)(23)}^{n-2}$ which is connected and all of its edges
commute with $\B_{n-2}$. It therefore follows that all of these
lassoes
 are equal to the lasso
 with head at $\rho_{2,n-2,n-1}$ labelled by $\B_{n-2}$. Noticing that the path
 $w=\B_3\B_4\dotsm\B_{n-3}$ connects $\rho_{2,3,n-1}^n$ to $\rho_{2,n-2,n-1}^n$
 one sees that the later lasso equals $[\B_{n-2}]\B_{n-3}^{-1}\dotsm\B_3^{-1}$.\\

\begin{description}
\item[Case b] Lassoes with heads at $\rho_{n-1}^n$.
\end{description}

Using as tail the path
$w=\B_{n-2}\B_{n-3}\dotsm\B_5\delta_4\B_4\B_5\dotsm\B_{n-2}$ one
sees that:
\begin{align*}w\B_0w^{-1}&=[\B_0]\B_{n-2}^{-1}\B_{n-3}^{-1}\dotsm\B_5^{-1}
         \B_4^{-1}\delta_4^{-1}\B_5^{-1}\dotsm\B_{n-3}^{-1}\B_{n-2}^{-1}\\
         &=\B_0
 \end{align*}

\begin{align*}w\B_2w^{-1}&=[\B_2]\B_{n-2}^{-1}\B_{n-3}^{-1}\dotsm\B_5^{-1}
         \B_4^{-1}\delta_4^{-1}\B_5^{-1}\dotsm\B_{n-3}^{-1}\B_{n-2}^{-1}\\
         &=\B_2
 \end{align*}

\begin{align*}w\B_3w^{-1}&=[\B_3]\B_{n-2}^{-1}\B_{n-3}^{-1}\dotsm\B_5^{-1}
         \B_4^{-1}\delta_4^{-1}\B_5^{-1}\dotsm\B_{n-3}^{-1}\B_{n-2}^{-1}\\
         &=[\B_3]\B_4^{-1}\delta_4^{-1}\B_5^{-1}\dotsm\B_{n-3}^{-1}\B_{n-2}^{-1}\\
         &=[\B_3]\delta_4^{-1}\B_4^{-1}\B_5^{-1}\dotsm\B_{n-3}^{-1}\B_{n-2}^{-1}\\
         &=[\B_3]\B_4^{-1}\B_5^{-1}\dotsm\B_{n-3}^{-1}\B_{n-2}^{-1}
 \end{align*}
which as seen in Figure \ref{Nolabel} equals
$[\B_{n-2}]\B_{n-3}^{-1}\B_{n-4}^{-1}\dotsm\B_4^{-1}$.\\

\begin{figure}[htp]
\[
\hspace{-2cm}
\begin{pspicture}(1,3.2)(7,9.6)
 \psdots(1,9)(2,9)(3,9)(5,9)(6,9)(7,9)
        (1,4)(2,4)(3,4)(5,4)(6,4)(7,4)
        \psdots[dotsize=0.05](3.4,9)(3.8,9)(4.2,9)(4.6,9)
                            (3.4,4)(3.8,4)(4.2,4)(4.6,4)
 \rput(1,9.3){14}
 \rput(2,9.3){23}
 \rput(3,9.3){23}
                 \rput(4,9.3){23's}
 \rput(5,9.3){23}
 \rput(6,9.3){23}
 \rput(7,9.3){14}
 \rput(1,3.7){14}
 \rput(2,3.7){23}
 \rput(3,3.7){23}
                  \rput(4,3.7){23's}
 \rput(5,3.7){23}
 \rput(6,3.7){23}
 \rput(7,3.7){14}
 \psline(1,7)(1.9,6.1)
 \psline(1.9,6.9)(1.56,6.56)
 \psline(1.44,6.44)(1,6)
 \psframe(2.7,4.6)(5.3,8.4)
 \psline(1,6)(1,4)
 \psline(1,9)(1,7)
 \psline(3,9)(3,8.4)
 \psline(3,4.6)(3,4)
 \psline(5,9)(5,8.4)
 \psline(5,4.6)(5,4)
 \psline(2,9)(2,4)
 \psline(6,9)(6,4)
  \psline(2.1,7.04)(2.6,7.24)
   \psline(2.8,7.32)(5.2,8.28)
   \psline(5.4,8.36)(5.9,8.56)
    \psline(6.1,8.64)(7,9)
    \psline(2.1,5.96)(2.6,5.76)
    \psline(2.8,5.68)(5.2,4.72)
    \psline(5.4,4.64)(5.9,4.44)
    \psline(6.1,4.36)(7,4)
                          \psline[linewidth=.005](3.8,9)(3.8,8.4)
                          \psline[linewidth=.005](4.2,9)(4.2,8.4)
                          \psline[linewidth=.005](4.6,9)(4.6,8.4)
                          \psline[linewidth=.005](3.4,9)(3.4,8.4)
                        \psline[linewidth=.005](3.8,4.6)(3.8,4)
                        \psline[linewidth=.005](4.2,4.6)(4.2,4)
                        \psline[linewidth=.005](4.6,4.6)(4.6,4)
                        \psline[linewidth=.005](3.4,4.6)(3.4,4)
                 \rput(3.5,5){id}
\end{pspicture}
\hspace{2cm}
\begin{pspicture}(1,3.2)(7,9.6)
 \psdots(1,9)(2,9)(3,9)(5,9)(6,9)(7,9)
        (1,4)(2,4)(3,4)(5,4)(6,4)(7,4)
        \psdots[dotsize=0.05](3.4,9)(3.8,9)(4.2,9)(4.6,9)
                            (3.4,4)(3.8,4)(4.2,4)(4.6,4)
 \rput(1,9.3){14}
 \rput(2,9.3){23}
 \rput(3,9.3){23}
                 \rput(4,9.3){23's}
 \rput(5,9.3){23}
 \rput(6,9.3){23}
 \rput(7,9.3){14}
 \rput(1,3.7){14}
 \rput(2,3.7){23}
 \rput(3,3.7){23}
                  \rput(4,3.7){23's}
 \rput(5,3.7){23}
 \rput(6,3.7){23}
 \rput(7,3.7){14}
  \psline(6.1,6.9)(7,6)
  \psline(7,7)(6.56,6.56)
  \psline(6.44,6.44)(6.1,6.1)
 \psframe(2.7,4.6)(5.3,8.4)
 \psline(7,6)(7,4)
 \psline(7,9)(7,7)
 \psline(3,9)(3,8.4)
 \psline(3,4.6)(3,4)
 \psline(5,9)(5,8.4)
 \psline(5,4.6)(5,4)
 \psline(2,9)(2,4)
 \psline(6,9)(6,4)
    \psline(5.9,5.96)(5.4,5.76)
    \psline(5.2,5.68)(2.8,4.72)
    \psline(2.6,4.64)(2.1,4.44)
    \psline(1.9,4.36)(1,4)
       \psline(5.9,7.04)(5.4,7.24)
       \psline(5.2,7.32)(2.8,8.28)
       \psline(2.6,8.36)(2.1,8.56)
       \psline(1.9,8.64)(1,9)
                         \psline[linewidth=.005](3.8,9)(3.8,8.4)
                          \psline[linewidth=.005](4.2,9)(4.2,8.4)
                          \psline[linewidth=.005](4.6,9)(4.6,8.4)
                          \psline[linewidth=.005](3.4,9)(3.4,8.4)
                        \psline[linewidth=.005](3.8,4.6)(3.8,4)
                        \psline[linewidth=.005](4.2,4.6)(4.2,4)
                        \psline[linewidth=.005](4.6,4.6)(4.6,4)
                        \psline[linewidth=.005](3.4,4.6)(3.4,4)
                 \rput(5,5){id}
\end{pspicture}
\]
\[
\begin{pspicture}(1,3.2)(7,9.6)
 \psdots(1,9)(2,9)(3,9)(5,9)(6,9)(7,9)
        (1,4)(2,4)(3,4)(5,4)(6,4)(7,4)
        \psdots[dotsize=0.05](3.4,9)(3.8,9)(4.2,9)(4.6,9)
                            (3.4,4)(3.8,4)(4.2,4)(4.6,4)
 \rput(1,9.3){14}
 \rput(2,9.3){23}
 \rput(3,9.3){23}
                 \rput(4,9.3){23's}
 \rput(5,9.3){23}
 \rput(6,9.3){23}
 \rput(7,9.3){14}
 \rput(1,3.7){14}
 \rput(2,3.7){23}
 \rput(3,3.7){23}
                  \rput(4,3.7){23's}
 \rput(5,3.7){23}
 \rput(6,3.7){23}
 \rput(7,3.7){14}
  \psline(6,7)(7,6)
  \psline(7,7)(6.56,6.56)
  \psline(6.44,6.44)(6,6)
 \psline(2.7,4.6)(5.3,4.6)
 \psline(2.7,8.4)(5.3,8.4)
 \psline(2.7,4.6)(2.7,4.63)
 \psline(2.7,4.8)(2.7,8.25)
 \psline(5.3,4.6)(5.3,5.6)
 \psline(5.3,5.8)(5.3,7.2)
 \psline(5.3,7.4)(5.3,8.4)
 \psline(7,6)(7,4)
 \psline(7,9)(7,7)
 \psline(3,9)(3,8.4)
 \psline(3,4.6)(3,4)
 \psline(5,9)(5,8.4)
 \psline(5,4.6)(5,4)
  \psline(2,9)(2,8.7)
  \psline(2,8.5)(2,4.5)
  \psline(2,4.3)(2,4)
  \psline(6,9)(6,7.1)
  \psline(6,6.9)(6,6.1)
  \psline(6,5.9)(6,4)
    \psline(6,6)(1,4)
       \psline(6,7)(1,9)
                         \psline[linewidth=.005](3.8,9)(3.8,8.4)
                          \psline[linewidth=.005](4.2,9)(4.2,8.4)
                          \psline[linewidth=.005](4.6,9)(4.6,8.4)
                          \psline[linewidth=.005](3.4,9)(3.4,8.4)
                        \psline[linewidth=.005](3.8,4.6)(3.8,4)
                        \psline[linewidth=.005](4.2,4.6)(4.2,4)
                        \psline[linewidth=.005](4.6,4.6)(4.6,4)
                        \psline[linewidth=.005](3.4,4.6)(3.4,4)
                 \rput(5,5){id}
\end{pspicture}
\]
 \caption{First shift to the right place by isotopy
and then use move $\mathcal{P}$ to pass over the strands in the
box.}\label{Nolabel}
\end{figure}
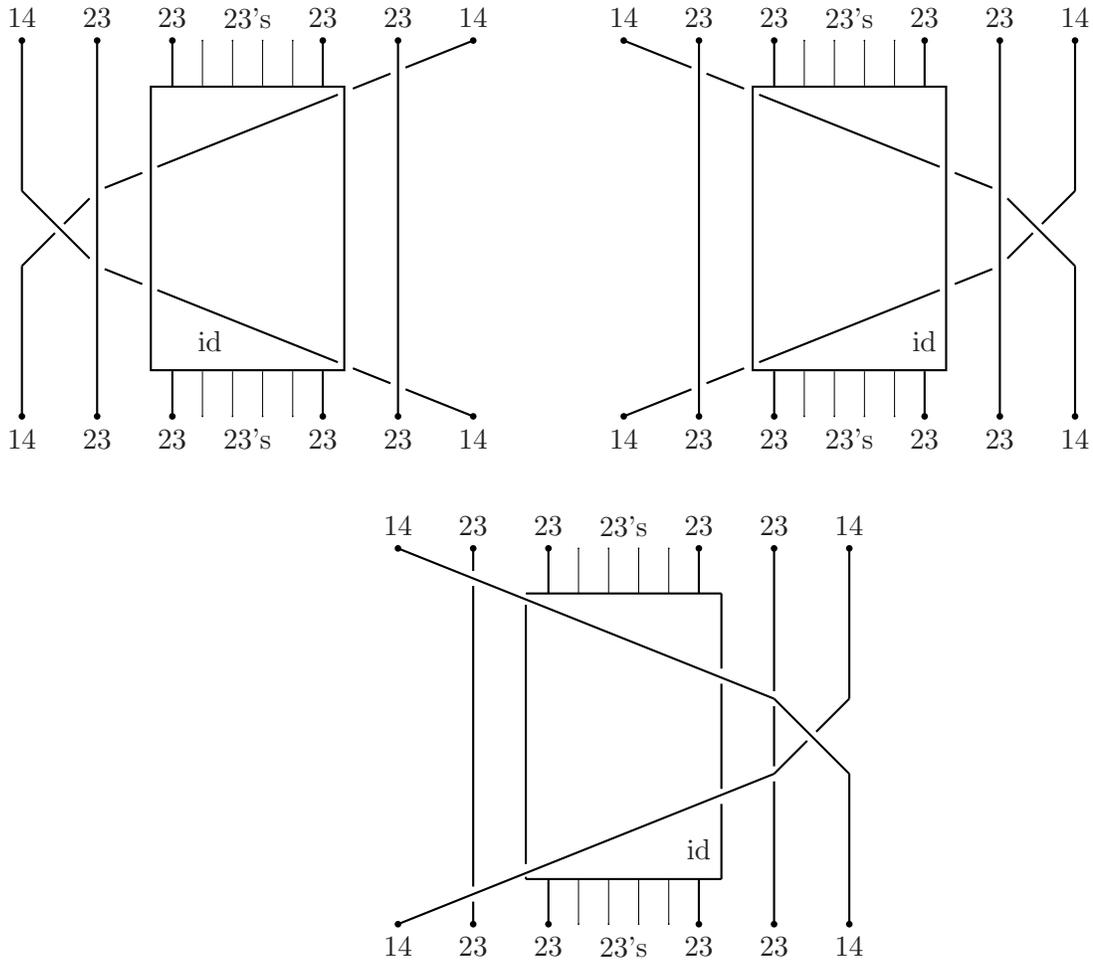

\begin{align*}w\B_4w^{-1}&=[\B_4]\B_{n-2}^{-1}\B_{n-3}^{-1}\dotsm\B_5^{-1}
         \B_4^{-1}\delta_4^{-1}\B_5^{-1}\dotsm\B_{n-3}^{-1}\B_{n-2}^{-1}\\
         &=[\B_4]\B_5^{-1}\B_4^{-1}\delta_4^{-1}\B_5^{-1}\dotsm\B_{n-3}^{-1}\B_{n-2}^{-1}\\
         &=[\B_5]\delta_4^{-1}\B_5^{-1}\dotsm\B_{n-3}^{-1}\B_{n-2}^{-1}\\
         &=[\delta_4]\B_6^{-1}\dotsm\B_{n-3}^{-1}\B_{n-2}^{-1}\\
         &=\delta_4
 \end{align*}

\begin{align*}w\delta_4w^{-1}&=[\delta_4]\B_{n-2}^{-1}\B_{n-3}^{-1}\dotsm\B_5^{-1}
         \B_4^{-1}\delta_4^{-1}\B_5^{-1}\dotsm\B_{n-3}^{-1}\B_{n-2}^{-1}\\
         &=[\delta_4]\B_5^{-1}\B_4^{-1}\delta_4^{-1}\B_5^{-1}\dotsm\B_{n-3}^{-1}\B_{n-2}^{-1}\\
         &=[\delta_4]\B_5^{-1}\delta_4^{-1}\B_4^{-1}\B_5^{-1}\dotsm\B_{n-3}^{-1}\B_{n-2}^{-1}\\
         &=[\B_5]\B_4^{-1}\B_5^{-1}\dotsm\B_{n-3}^{-1}\B_{n-2}^{-1}\\
         &=[\B_4]\B_6^{-1}\dotsm\B_{n-3}^{-1}\B_{n-2}^{-1}\\
         &=\B_4
 \end{align*}

\begin{align*}w\B_5w^{-1}&=[\B_5]\B_{n-2}^{-1}\B_{n-3}^{-1}\dotsm\B_5^{-1}
         \B_4^{-1}\delta_4^{-1}\B_5^{-1}\dotsm\B_{n-3}^{-1}\B_{n-2}^{-1}\\
         &=[\B_5]\B_6^{-1}\B_5^{-1}\B_4^{-1}\delta_4^{-1}\B_5^{-1}
\dotsm\B_{n-3}^{-1}\B_{n-2}^{-1}\\
         &=[\B_6]\delta_4^{-1}\B_4^{-1}\B_5^{-1}\dotsm\B_{n-3}^{-1}\B_{n-2}^{-1}\\
         &=[\B_6]\B_5^{-1}\B_6^{-1}\dotsm\B_{n-3}^{-1}\B_{n-2}^{-1}\\
         &=[\B_5]\B_7^{-1}\dotsm\B_{n-3}^{-1}\B_{n-2}^{-1}\\
         &=\B_5
 \end{align*}
 and similarly for $i=6,\ldots,\B_{n-3}$
$$w\B_iw^{-1}=\B_i$$

\begin{description}
\item[Case c] Lassoes with heads at
$i_{(23)}(\CC_{(14)}^{n-1})$.
\end{description}

By symmetry and case a) one needs only to check lassoes with head
at $\tilde\rho_{2,3,n-1}^n$, the lasso with head at
$\tilde\rho_{2,n-2,n-1}$ labelled by $\B_{n-2}$ (corresponding to
$[\B_{n-2}]\B_{n-3}^{-1}\dotsm\B_3^{-1}$),and for $n\geq 10$ the
lasso with head at $\tilde\rho_{2,3,n-1}^n$ labelled by
$\delta_4$(corresponding to $\delta_6$).\\
For the lassoes with head at $\tilde\rho_{2,3,n-1}^n$ labelled by
$\delta_4$ or $\B_i$ for $i=2,4,5,\ldots\B_{n-4}$ notice that
$\B_{n-2}$ connects $\tilde\rho_{2,3,n-1}^n$ to
$\tilde\rho_{2,3,n-2}^n\in i_{(14)}(\CC_{(12)}^{n-1})$ and so they
have been checked in case a). For the one with head labelled by
$\B_{n-3}$ ``pull further into $i_{(14)}(\CC_{(12)}^{n-1})$'' by
using the path $\B_{n-2}\B_{n-3}$.\\
For the lasso with head at $\tilde\rho_{2,n-2,n-1}$ labelled by
$\B_{n-2}$ use again the path $\B_{n-2}\B_{n-3}$ to ``pull it'' to
the lasso with head at $\tilde\rho_{2,n-3,n-2}\in
i_{(14)}(\CC_{(12)}^{n-1})$, labelled
by $\B_{n-3}$.\\
Finally for the lasso with head at $\tilde\rho_{2,3,n-1}^n$
labelled by $\delta_4$ use $\B_{n-2}$ to ``pull'' it to
$\tilde\rho_{2,3,n-1}^n \in i_{(14)}(\CC_{(23)}^{n-1})$ with label
$\delta_4$.

\begin{description}
\item[Case d] Lassoes with heads at $\tilde\rho_{n-1}^n$.
\end{description}
Analogously to case b) those are covered by case c).\\
This concludes the proof for $\overline{\Au}(\rho_{2,3,n-1}^n$ and
 the even $n$ case.\\

 For odd $n$:

\begin{align*} \CC(\rho_{23}^n) &=\CC_{(23)}^n\\
                                &=i_{(23)}(\CC_{\text{id}}^{n-1}\sqcup
                                 \{\tilde\rho_\emptyset^{n-1}\})
                              \bigsqcup i_{(14)}(\CC_{(14)(23)}^{n-1})
\end{align*}

and the lassoes fall into three cases:

\begin{description}
\item[Case a] Lassoes with heads at
$i_{(23)}(\CC_{\text{id}}^{n-1})$
\end{description}

 In this case by the induction hypothesis, one needs only to check lassoes with head
labelled by $\B_{n-2}$ . After writing
$$\CC_{\text{id}}^{n-1}=i_{(23)}(\CC_{(23)}^{n-2})
 \bigsqcup i_{(14)}(\CC_{(14)}^{n-2})$$
 and observing that for lassoes in
 $i_{(23)}i_{(23)}(\CC_{(23)}^{n-2})$ the tail can be chosen
 to commute with $\B_{n-2}$, one sees that all of these lassoes
 are equal to $\B_{n-2}$.

\begin{description}
\item[Case b] Lassoes with heads at
$i_{(14)}(\CC_{(14)(23)}^{n-1})$
\end{description}
By the induction hypothesis one needs only to check the lassoes
with heads at $\rho_{2,3,n-2,n-1}$, the lasso with head at
$\rho_{2,n-3,n-2,n-1}$ labelled by $\B_{n-3}$, and the lasso with
head at $\rho_{5,6,n-2,n-1}$ labelled by
$\delta_4$.\\
For the lassoes with head at $\rho_{2,3,n-2,n-1}$ labelled by
$\delta_4$ or by $\B_i$ with $i=2,4,5,\dotsc,\B_{n-5}$ use the
path $\B_{n-2}\B_{n-3}$ to pull them at $\rho_{2,3,n-3,n-2}\in
i_{(23)}(\CC_{\text{id}}^{n-1})$ with the same label and refer to
case a).\\
For the lasso with head at $\rho_{2,3,n-2,n-1}$ labelled by
 $\B_{n-2}$ use the path $\B_{n-2}\B_{n-3}$ to pull it at
 $\rho_{2,3,n-3,n-2}\in i_{(23)}(\CC_{\text{id}}^{n-1})$ with label
 $\B_{n-3}$ and refer to case a).\\
For the lasso with head at $\rho_{2,3,n-2,n-1}$ labelled by
 $\B_{n-4}$ use the path $\B_{n-3}\B_{n-2}\B_{n-4}\B_{n-3}$ to pull it at
 $\rho_{2,3,n-4,n-3}\in i_{(23)}(\CC_{\text{id}}^{n-1})$ with label
 $\B_{n-2}$ and refer to case a).\\
 For the lasso
with head at $\rho_{2,n-3,n-2,n-1}$ labelled by $\B_{n-3}$, use
the path $\B_{n-2}\B_{n-3}\B_{n-4}$ to pull it at
 $\rho_{2,n-4,n-3,n-2}\in i_{(23)}(\CC_{\text{id}}^{n-1})$ with label
 $\B_{n-4}$ and refer to case a).\\
 Finally for  the lasso with head at $\rho_{5,6,n-2,n-1}$ labelled by
$\delta_4$, use the path $\B_{n-2}\B_{n-3}$ to pull it at
 $\rho_{5,6,n-3,n-2}\in i_{(23)}(\CC_{\text{id}}^{n-1})$ with label
 $\delta_4$ and refer to case a).\\

\begin{description}
\item[Case c] Lassoes with heads at
$\tilde\rho_{n-1}^n$
\end{description}

For those with heads labelled by $\delta_4$ or $\B_i$ with
$i=1,2,\ldots,n-4$ use $\B_{n-2}$ to pull them at
$\tilde\rho_{n-2}^n\in i_{(14)}(\CC_{(14)(23)}^{n-1})$ with the
same
label, and refer to case b).\\
For the one with head labelled by $\B_{n-3}$ use
$\B_{n-3}\B_{n-2}$ to pull it at $\tilde\rho_{n-3}^n\in
i_{(14)}(\CC_{(14)(23)}^{n-1})$ with label $\B_{n-2}$, and refer
to case b).

\end{proof}

\subsection{The Kernel}\label{thekernel} In this subsection we
compute the kernel of the the (quotiented) lifting homomorphism in a way
completely analogous to \cite{BW1}, see
Section~\ref{2sphere3}.\

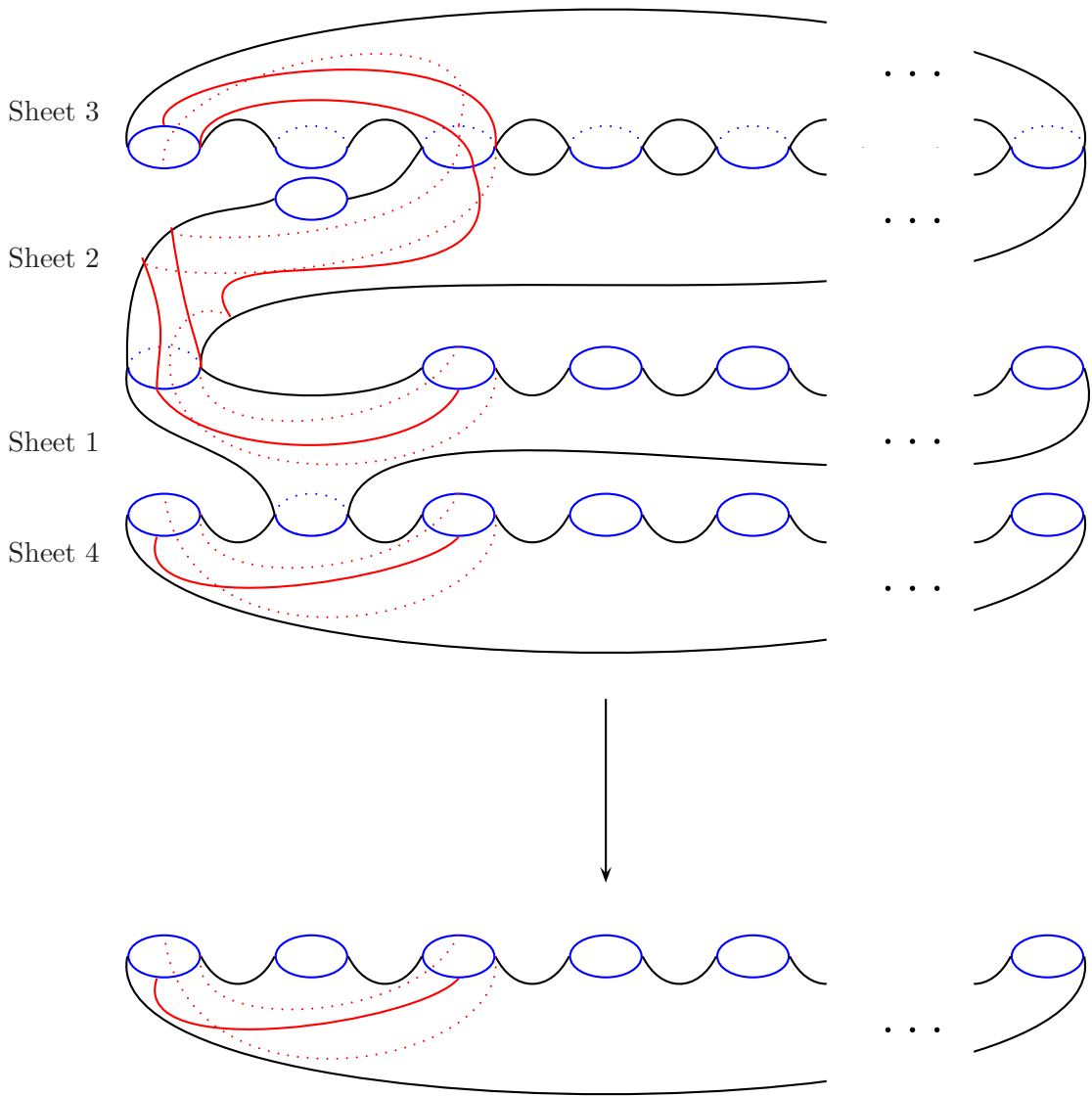
\begin{figure}[htp]
\[
\begin{pspicture}(-1,-1)(15,15)
 \psellipse[linecolor=blue](.5,1.5)(.5,.3)
 \psellipse[linecolor=blue](2.5,1.5)(.5,.3)
 \psellipse[linecolor=blue](4.5,1.5)(.5,.3)
 \psellipse[linecolor=blue](6.5,1.5)(.5,.3)
 \psellipse[linecolor=blue](8.5,1.5)(.5,.3)
 \psellipse[linecolor=blue](12.5,1.5)(.5,.3)
 \psbezier(1,1.5)(1.25,1)(1.75,1)(2,1.5)
  \psbezier(3,1.5)(3.25,1)(3.75,1)(4,1.5)
   \psbezier(5,1.5)(5.25,1)(5.75,1)(6,1.5)
    \psbezier(7,1.5)(7.25,1)(7.75,1)(8,1.5)
     \psbezier(9,1.5)(9.25,1)(9.75,1)(10,1.5)
          \psbezier(11,1.5)(11.25,1)(11.75,1)(12,1.5)
     \psbezier
     (0,1.5)(-.5,-1)(13.5,-1)(13,1.5)
  \psframe[linecolor=white,fillstyle=solid,fillcolor=white](9.5,-1)(11.5,2)
  \psdots[dotsize=.1](10.3333,.5)(10.6666,.5)(10.9999,.5)
  \psbezier[linestyle=dotted,linecolor=red](5,1.5)(5.2,0)(.8,-.8)(.5,1.8)
  \psbezier[linecolor=red](.4,1.2)(0,0)(4,.6)(4.5,1.2)
  \psbezier[linestyle=dotted,linecolor=red](4.5,1.8)(4,.6)(.8,.4)(1,1.5)

 \rput(0,6){\psellipse[linecolor=blue](2.5,1.5)(.5,.3)
 \psframe[linecolor=white,fillstyle=solid,fillcolor=white](0,1.5)(13,2)
\psellipse[linecolor=blue](.5,1.5)(.5,.3)
 \psellipse[linestyle=dotted,linecolor=blue](2.5,1.5)(.5,.3)
  \psellipse[linecolor=blue](4.5,1.5)(.5,.3)
 \psellipse[linecolor=blue](6.5,1.5)(.5,.3)
 \psellipse[linecolor=blue](8.5,1.5)(.5,.3)
 \psellipse[linecolor=blue](12.5,1.5)(.5,.3)
 \psbezier(1,1.5)(1.25,1)(1.75,1)(2,1.5)
  \psbezier(3,1.5)(3.25,1)(3.75,1)(4,1.5)
   \psbezier(5,1.5)(5.25,1)(5.75,1)(6,1.5)
    \psbezier(7,1.5)(7.25,1)(7.75,1)(8,1.5)
     \psbezier(9,1.5)(9.25,1)(9.75,1)(10,1.5)
          \psbezier(11,1.5)(11.25,1)(11.75,1)(12,1.5)
     \psbezier
     (0,1.5)(-.5,-1)(13.5,-1)(13,1.5)
  \psframe[linecolor=white,fillstyle=solid,fillcolor=white](9.5,-1)(11.5,2)
  \psdots[dotsize=.1](10.3333,.5)(10.6666,.5)(10.9999,.5)
  \psbezier[linestyle=dotted,linecolor=red](5,1.5)(5.2,0)(.8,-.8)(.5,1.8)
  \psbezier[linecolor=red](.4,1.2)(0,0)(4,.6)(4.5,1.2)
  \psbezier[linestyle=dotted,linecolor=red](4.5,1.8)(4,.6)(.8,.4)(1,1.5)}

  \rput(0,8){\psellipse[linecolor=blue](.5,1.5)(.5,.3)
  \psframe[linecolor=white,fillstyle=solid,fillcolor=white](0,1.5)(13,2)
  \psellipse[linestyle=dotted,linecolor=blue](.5,1.5)(.5,.3)
 \psellipse[linecolor=blue](2.5,1.5)(.5,.3)
 \psellipse[linecolor=blue](4.5,1.5)(.5,.3)
 \psellipse[linecolor=blue](6.5,1.5)(.5,.3)
 \psellipse[linecolor=blue](8.5,1.5)(.5,.3)
 \psellipse[linecolor=blue](12.5,1.5)(.5,.3)
 \psbezier(1,1.5)(1.25,1)(1.75,1)(2,1.5)
  \psbezier(3,1.5)(3.25,1)(3.75,1)(4,1.5)
   \psbezier(5,1.5)(5.25,1)(5.75,1)(6,1.5)
    \psbezier(7,1.5)(7.25,1)(7.75,1)(8,1.5)
     \psbezier(9,1.5)(9.25,1)(9.75,1)(10,1.5)
         \psbezier(11,1.5)(11.25,1)(11.75,1)(12,1.5)
       \psframe[linecolor=white,fillstyle=solid,fillcolor=white](9.5,-1)(11.5,2)
  \psdots[dotsize=.1](10.3333,.5)(10.6666,.5)(10.9999,.5)}

  \rput(0,11){\psellipse[linecolor=blue](.5,1.5)(.5,.3)
  \psellipse[linecolor=blue](2.5,1.5)(.5,.3)
 \psellipse[linecolor=blue](4.5,1.5)(.5,.3)
 \psellipse[linecolor=blue](6.5,1.5)(.5,.3)
 \psellipse[linecolor=blue](8.5,1.5)(.5,.3)
 \psellipse[linecolor=blue](12.5,1.5)(.5,.3)
 \psframe[linecolor=white,fillstyle=solid,fillcolor=white](0,1.5)(13,2)
 \psellipse[linestyle=dotted,linecolor=blue](2.5,1.5)(.5,.3)
 \psellipse[linestyle=dotted,linecolor=blue](4.5,1.5)(.5,.3)
 \psellipse[linestyle=dotted,linecolor=blue](6.5,1.5)(.5,.3)
 \psellipse[linestyle=dotted,linecolor=blue](8.5,1.5)(.5,.3)
 \psellipse[linestyle=dotted,linecolor=blue](12.5,1.5)(.5,.3)
 \psbezier(1,1.5)(1.25,1)(1.75,1)(2,1.5)
   \psbezier(5,1.5)(5.25,1)(5.75,1)(6,1.5)
    \psbezier(7,1.5)(7.25,1)(7.75,1)(8,1.5)
     \psbezier(9,1.5)(9.25,1)(9.75,1)(10,1.5)
          \psbezier(11,1.5)(11.25,1)(11.75,1)(12,1.5)

  \psframe[linecolor=white,fillstyle=solid,fillcolor=white](9.5,-1)(11.5,2)
  \psdots[dotsize=.1](10.3333,.5)(10.6666,.5)(10.9999,.5)}

 \rput(-1,7){Sheet 4}
 \rput(-1,8.5){Sheet 1}
 \rput(-1,11){Sheet 2}
 \rput(-1,13){Sheet 3}
 \psframe[linecolor=white,fillstyle=solid,fillcolor=white](1,8)(4,11)
\psframe[linecolor=white,fillstyle=solid,fillcolor=white](0,10)(2,13)

\psellipse[linecolor=blue](.5,12.5)(.5,.3)
 \psbezier(1,12.5)(1.25,13)(1.75,13)(2,12.5)
 \psbezier(5,12.5)(5.25,13)(5.75,13)(6,12.5)
 \psbezier(3,12.5)(3.25,13)(3.75,13)(4,12.5)
 \psbezier(7,12.5)(7.25,13)(7.75,13)(8,12.5)
 \psbezier(9,12.5)(9.25,13)(9.75,13)(10,12.5)
 \psbezier(11,12.5)(11.25,13)(11.75,13)(12,12.5)
 \psbezier(0,12.5)(-.5,15)(13.5,15)(13,12.5)
 \psframe[linecolor=white,fillstyle=solid,fillcolor=white](9.5,12.5)(11.5,15)
  \psdots[dotsize=.1](10.3333,13.5)(10.6666,13.5)(10.9999,13.5)

\psbezier
     (0,9.5)(-.2,8.4)(1.8,8.6)(2,7.5)
  \psbezier
     (3,7.5)(3.2,9.8)(14,6.5)(13,9.5)
  \psbezier
     (1,9.5)(1.5,9)(3.5,9)(4,9.5)
  \psbezier
     (0,9.5)(-.1,12)(1.5,11.5)(2,11.8)   \psbezier
     (1,9.5)(.8,12)(13.5,9)(13,12.5)

  \psbezier[linestyle=dotted,linecolor=red](5,9.5)(5.2,7.8)(.2,7.6)(.6,9.8)
  \psbezier[linecolor=red](.4,9.2)(1,8.2)(4,8.2)(4.5,9.2)
  \psbezier[linestyle=dotted,linecolor=red](4.5,9.8)(4,8.6)(.8,8.4)(1,9.5)

  \psbezier[linecolor=red](5,12.5)(5.2,14.2)(.2,13.5)(.5,12.8)
  \psbezier[linecolor=red](4.7,12.2)(4.7,13.4)(.8,13.4)(1,12.5)
  \psbezier[linestyle=dotted,linecolor=red](4.5,12.8)(5.2,14.6)(.2,13.6)(.5,12.2)

  \psbezier[linestyle=dotted,linecolor=red](5,12.5)(5.2,10.6)(0,10.6)(.2,11)
  \psbezier[linecolor=red](.2,11)(.6,10)(.4,9.8)(.4,9.2)
  \psbezier[linestyle=dotted,linecolor=red](4.5,12.8)(4.8,11)(0,11.2)(.6,11.4)
  \psbezier[linecolor=red](.6,11.4)(.8,10)(1,9.8)(1,9.5)
  \psbezier[linestyle=dotted,linecolor=red](.6,9.8)(.8,10.4)(1.2,10.3)(1.4,10.2)
  \psbezier[linecolor=red](1.4,10.2)(.5,11.5)(5.5,10)(4.7,12.2)

\psframe[linecolor=white,fillstyle=solid,fillcolor=white](9.5,8)(11.5,8.3)
\psframe[linecolor=white,fillstyle=solid,fillcolor=white](9.5,10)(11.5,11.3) \psline{->}(6.5,5)(6.5,2.5)

\psellipse[linecolor=blue](2.5,11.8)(.5,.3)
  \pscurve(3,11.8)(3.6,12)(4,12.5)
\end{pspicture}
\]
\caption{The covering $\rho_n$ and how $d_4$ lifts}\label{malako}
\end{figure}

 For this subsection covering
means a covering
over $S^2$, in particular the number of branch values of a
covering is even. Also notation is simplified by writing $\rho_n$
instead of $\rho_{23}^n$. Then using the Riemann-Hurwitz formula
one gets that the total space of $\rho_n$ has genus
$$g=\frac{n-6}{2}.$$
In fact Figure~\ref{malako} shows an explicit model of
$\rho_{2g+6}$ as a genus-$g$ surface. This model is constructed by
cutting and pasting. $S^2$ is cut open along the intervals $x_i$
for even $i$. The cuts are shown in blue and it is assumed that
the boundary circles in this picture are ``sewed'' to become
intervals.

The commutative diagram (\ref{lambda}) then becomes:

$$\xymatrix{
 {\Au(\rho_{2g+6})} \ar[d] \ar[r]^-{\lambda} &{\CM_g}\\
{ \overline{\Au}(\rho_{2g+6}) }\ar[ur]_-{\bar{\lambda}}   } $$

Notice that all the generators of $\overline{\Au}(\rho_n)$ given
in Theorem~\ref{main} are rotations around intervals and thus
their image under $\bar{\lambda}$ can be determined according to
Lemma~\ref{intlift} by lifting the intervals. Using the model for
$\rho_{2g+6}$ in Figure~\ref{malako} one gets the following:

\begin{thm}$\bar{\lambda}$ is determined by:

\begin{equation}\label{themap}
\bar\lambda(x)=
\begin{cases} \text{id} & \text{if} \quad x=\B_0,\B_2\\
                       a_{i-1} & \text{if $ \quad x=\B_{2i}$ for $i\geq2$}\\
                       b_{i-1} &\text{if  $\quad x=\B_{2i+1}$ with
                       $i\geq2$}\\
                       a_1 & \text{if $ \quad x=\delta_4$}\\
                       d & \text{if $ \quad x=\delta_6$}
                   \end{cases}
  \end{equation}
where $a_i$, $b_i$, and $d$ are Wajnryb's generators of $\CM_g$ as
described in Appendix~\ref{Pres}.
\end{thm}

\begin{proof} Notice that $\rho_n$ is obtained by $\rho(n)$ by
adding a fourth trivial sheet. After ``cutting off'' the fourth
 sheet as in the proof of Lemma~\ref{deltakey} $d_6$ becomes $d_4$
 which as shown in Theorem~\ref{bw2} in the $3$-sheeted case lifts
 to $d$.\\
 In Figure~\ref{malako}, the interval $d_4$ is
shown and how it lifts to the disjoint union
 of a curve isotopic to $a_1$ and two arcs.\\
 The checks for the remaining generators are rather trivial.
\end{proof}
One can now prove:
\begin{thm}\label{kernel} The kernel of $\bar{\lambda}$ is the
smallest normal subgroup of $\overline{\Au}(\rho_{2g+6})$
containing the elements $\B_0$, $\B_2$, $B$, $D$, and all words
obtained by $B$ or $D$ by replacing some appearances of $\B_4$
with $\delta_4$, where:
$$B=(\B_4\B_5\B_6)^4([\delta_6^{-1}]\B_6^{-1}\B_5^{-1}\B_4^{-2}
\B_5^{-1}\B_6^{-1}\B_7^{-1})\delta_6^{-1},$$
$$D=\B_{2g+4}\chi\B_{2g+4}^{-1}\chi^{-1},$$
Where:
$$\chi=\B_{2g+3}\B_{2g+2}\dotsm\B_5\B_4^2\B_5\dotsm\B_{2g+2}\B_{2g+3},$$
\end{thm}
\begin{proof} Given the Wajnryb's presentation of $\CM_g$  and
the fact that $\bar{\lambda}$ is given by \eqref{themap} the proof
reduces to checking that the words obtained by the relators in
\cite{W} after replacing in all possible ways the generators of
$\CM_g$ with their preimages to generators of
$\overline{\Au}(\rho_{2g+6})$, are in the normal closure of the
given elements.\

 The calculations are identical to those in Theorems 5.1 and 6.1 of
 \cite{BW1} after replacing of $x_i$ in \cite{BW1} with $\B_{i+2}$
 for $i\geq2$, and $\delta_4$ in \cite{BW1} with $\delta_6$,
 and then replacing $\B_4$ with $\delta_4$ in all possible ways.
\end{proof}

\section{Four sheets II: Coverings over $S^3$}\label{3sphere4}

\subsection{Bi-tricolored links}

Recall from Section~\ref{3sphere} how colored links represent $3$-manifolds.
 Also recall that $4$-colored links
are called bi-tricolored. That name is justified by the
following combinatorial description which is analogous to
(but of course more complicated than) the combinatorial description of
tricolored links given in
Section~\ref{3sphere3}.

\begin{defn} A bi-tricolored link  diagram is a link
diagram whose arcs are colored by three colors that come in
light and dark shades (see Section~\ref{dim}), in such a manner that:
\begin{itemize}
\item By disregarding the shades (i.e. dimming the lights) we get a tricolored diagram.
\item There is at least one dark shade.
\item At trichromatic crossings there is an even number of dark shades,
while at monochromatic crossings the two
pieces of the under strand have the same shade.
\end{itemize}
\end{defn}

Recall also that a normalized (bitricolored) diagram is a
link diagram in plat form whose top and bottom are
equal to:

\[
\rho_{\text{stand}}=\begin{array}{lr}
\begin{xy}/r.9cm/:
\vcap|{12}
\end{xy}
\hspace{5mm}
\begin{xy}/r.9cm/:
\vcap|{14}
\end{xy}
\hspace{5mm}
\begin{xy}/r.9cm/:
\vcap|{23}
\end{xy}
\hspace{2.5mm} \cdots \hspace{2.5mm}
\begin{xy}/r.9cm/:
\vcap|{23}
\end{xy}
&\\
& \\
\end{array}
\]

Quite often in the next subsection a slightly different normalization will be used,
namely the positions of
$(12)$'s and $(14)$'s will be interchanged. The proof of Lemma~\ref{Previous}
shows how this can be done.

\subsection{Moves} 
\begin{figure}[htp]
\[
\begin{pspicture}(-.5,0)(5,5)
 \psdots(0,4)(.5,4)(1,4)(1.5,4)(2,4)(2.5,4)(4,4)(4.5,4)
        (0,1)(.5,1)(1,1)(1.5,1)(2,1)(2.5,1)(4,1)(4.5,1)
                \psframe(-.3,0)(4.8,1)
           \psframe(-.3,4)(4.8,5)
\psellipse(1.25,2.5)(.5,.3)
 \psframe[linecolor=white,fillstyle=solid,fillcolor=white](.9,2.65)(1.1,2.8)
 \psframe[linecolor=white,fillstyle=solid,fillcolor=white](1.4,2.65)(1.6,2.8)
\psline(1,4)(1,2.35) \psline(1,2.15)(1,1)
   \psline(1.5,4)(1.5,2.35)\psline(1.5,2.15)(1.5,1)
   \psline(2,4)(2,1)
   \psline(2.5,4)(2.5,1)
   \psline(4,4)(4,1)
   \psline(4.5,4)(4.5,1)
 \psline(0,4)(0,1)
 \psline(0.5,4)(0.5,1)
    \rput(0,4.3){14}
 \rput(0.5,4.3){14}
 \rput(1,4.3){12}
 \rput(1.5,4.3){12}
                  \rput(3.2,4.3){light shades}
 \rput(0,.7){14}
 \rput(0.5,.7){14}
 \rput(1,.7){24}
 \rput(1.5,.7){24}
                   \rput(3.2,.7){light shades}
  \psdots[dotsize=0.1](3,4)(3.3,4)(3.6,4)(3,1)(3.3,1)(3.6,1)
                      (3,2.5)(3.3,2.5)(3.6,2.5)
     \rput(2.5,4.8){$\tau$}
     \rput(2.5,.2){$\tau'$}
\end{pspicture}
\hspace{2cm}
\begin{pspicture}\begin{pspicture}(-.5,0)(5,5)
 \psdots(0,4)(.5,4)(1,4)(1.5,4)(2,4)(2.5,4)(4,4)(4.5,4)
        (0,1)(.5,1)(1,1)(1.5,1)(2,1)(2.5,1)(4,1)(4.5,1)
               \psframe(-.3,0)(4.8,1)
           \psframe(-.3,4)(4.8,5)
\psellipse(.75,2.5)(.9,.3)
 \psframe[linecolor=white,fillstyle=solid,fillcolor=white](.9,2.65)(1.1,2.8)
 \psframe[linecolor=white,fillstyle=solid,fillcolor=white](1.4,2.55)(1.6,2.7)
  \psframe[linecolor=white,fillstyle=solid,fillcolor=white](.4,2.65)(.6,2.8)
 \psframe[linecolor=white,fillstyle=solid,fillcolor=white](-.1,2.55)(.1,2.7)
\psline(1,4)(1,2.35) \psline(1,2.15)(1,1)
   \psline(1.5,4)(1.5,2.45)\psline(1.5,2.25)(1.5,1)
   \psline(2,4)(2,1)
   \psline(2.5,4)(2.5,1)
   \psline(4,4)(4,1)
   \psline(4.5,4)(4.5,1)
 \psline(0,4)(0,2.45) \psline(0,2.25)(0,1)
 \psline(0.5,4)(.5,2.35) \psline(.5,2.15)(0.5,1)
 \rput(-.45,2.5){14}
    \rput(0,4.3){14}
 \rput(0.5,4.3){14}
 \rput(1,4.3){12}
 \rput(1.5,4.3){12}
                   \rput(3.2,4.3){light shades}
 \rput(0,.7){14}
 \rput(0.5,.7){14}
 \rput(1,.7){24}
 \rput(1.5,.7){24}
                    \rput(3.2,.7){light shades}
  \psdots[dotsize=0.1](3,4)(3.3,4)(3.6,4)(3,1)(3.3,1)(3.6,1)
                      (3,2.5)(3.3,2.5)(3.6,2.5)
     \rput(2.5,4.8){$\tau$}
     \rput(2.5,.2){$\tau'$}
\end{pspicture}
\]
 \caption{Move $\mathcal{X}$}\label{X}
\end{figure}
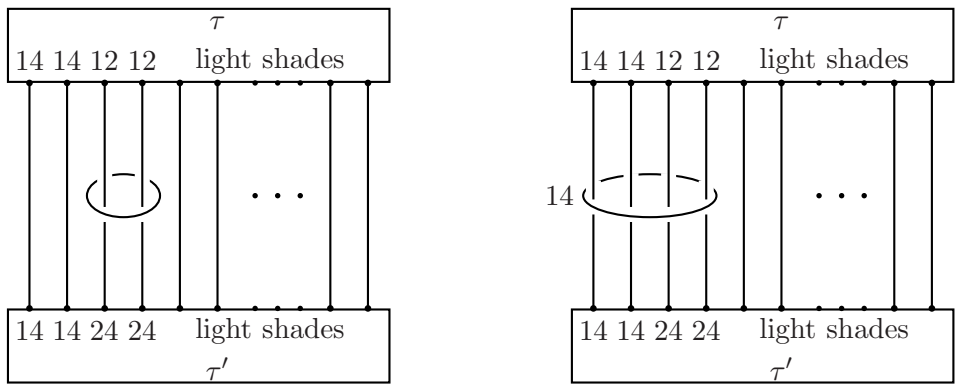

The goal of this subsection is to prove the
following two theorems:
\begin{thm}\label{theoremX}  Two bi-tricolored links represent the same
 $3$-manifold iff they can be related by a finite sequence of
 moves $\mathcal{M}$, $\mathcal{P}$ and the non-local moves
  $\mathcal{X}$, II and V shown in Figures \ref{X}, \ref{II} and \ref{V} respectively.
 \end{thm}

\begin{thm}\label{stab} Moves $\mathcal{M}$, $\mathcal{P}$ plus allowing the
addition/deletion of a trivial fifth sheet are enough to relate
any two bi-tricolored links representing the same $3$-manifold.
\end{thm}

The following lemma will be used through out the proofs.
\begin{lem}The move in Figure \ref{circumcision} is a consequence of move
$\mathcal{M}$.
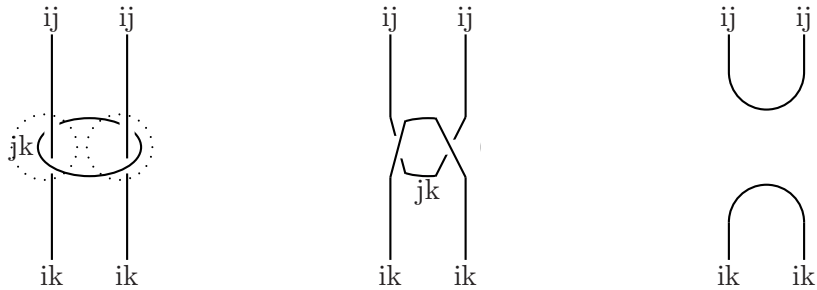
\begin{figure}[h]
\[
\begin{pspicture}(0,-.7)(3.5,2.7)
 \rput(1,2.7){ij} \rput(2,2.7){ij} \rput(1,-.7){ik} \rput(2,-.7){ik}
   \psellipse(1.5,1)(.7,.4)
    \psframe[linecolor=white,fillstyle=solid,fillcolor=white](.9,1.1)(1.1,1.4)
    \psframe[linecolor=white,fillstyle=solid,fillcolor=white](1.9,1.1)(2.1,1.4)
 \psline(1,2.5)(1,.85) \psline(1,.65)(1,-.5)
   \psline(2,2.5)(2,.85) \psline(2,.65)(2,-.5)
  \rput(.6,1){jk}
\end{pspicture}
\hspace{3cm}
\begin{pspicture}(0,-.7)(3.5,2.7)
 \rput(1,2.7){ij} \rput(2,2.7){ij} \rput(1,-.7){ik}
\rput(2,-.7){ik}
 \psline(1,2.5)(1,2) \psline(2,2.5)(2,2)
      \psarc(1.5,2){.5}{180}{360}
      \psarc(1.5,0){.5}{0}{180}
 \psline(1,-.5)(1,0) \psline(2,-.5)(2,0)
\end{pspicture}
\]
\caption{The circumcision move}\label{circumcision}
\end{figure}
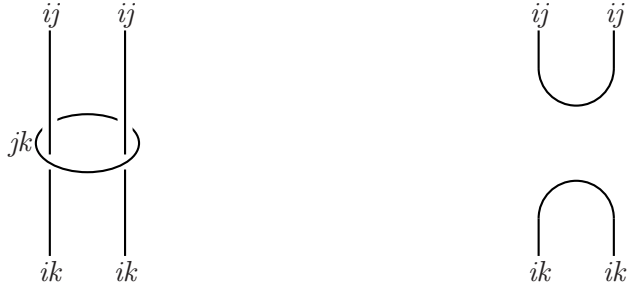
\end{lem}
\begin{proof}
\[
\begin{pspicture}(0,-1.5)(3.5,2.7)
 \rput(1,2.7){ij} \rput(2,2.7){ij} \rput(1,-.7){ik} \rput(2,-.7){ik}
   \psellipse(1.5,1)(.7,.4)
    \psframe[linecolor=white,fillstyle=solid,fillcolor=white](.9,1.1)(1.1,1.4)
    \psframe[linecolor=white,fillstyle=solid,fillcolor=white](1.9,1.1)(2.1,1.4)
 \psline(1,2.5)(1,.85) \psline(1,.65)(1,-.5)
   \psline(2,2.5)(2,.85) \psline(2,.65)(2,-.5)
  \rput(.6,1){jk}
 \pscircle[linestyle=dotted](.9,1){.45}
 \pscircle[linestyle=dotted](1.9,1){.45}
\end{pspicture}
\hspace{1cm}
\begin{pspicture}(0,-1.5)(3.5,2.7)
 \rput(1,2.7){ij} \rput(2,2.7){ij} \rput(1,-.7){ik} \rput(2,-.7){ik}
   \psellipse(1.5,1)(.7,.4)
 \psline(1,2.5)(1,.85) \psline(1,.65)(1,-.5)
   \psline(2,2.5)(2,.85) \psline(2,.65)(2,-.5)
       \psframe[linecolor=white,fillstyle=solid,fillcolor=white](.6,.6)(1.2,1.4)
    \psframe[linecolor=white,fillstyle=solid,fillcolor=white](1.6,.6)(2.2,1.4)
   \psline(1,1.4)(1.2,.645) \psline(2,1.4)(1.6,.61)
 \psframe[linecolor=white,fillstyle=solid,fillcolor=white](1.7,.9)(1.9,1.1)
\psframe[linecolor=white,fillstyle=solid,fillcolor=white](1,.85)(1.2,1.15)
\psline(2,.6)(1.6,1.39) \psline(1,.6)(1.2,1.355)
  \rput(1.5,.4){jk}
\end{pspicture}
 \hspace{1cm}
\begin{pspicture}(0,-1.5)(3.5,2.7)
 \rput(1,2.7){ij} \rput(2,2.7){ij} \rput(1,-.7){ik}
\rput(2,-.7){ik}
 \psline(1,2.5)(1,2) \psline(2,2.5)(2,2)
      \psarc(1.5,2){.5}{180}{360}
      \psarc(1.5,0){.5}{0}{180}
 \psline(1,-.5)(1,0) \psline(2,-.5)(2,0)
\end{pspicture}
\]

 Perform $\mathcal{M}$ moves inside the dotted circles and then
 isotope.
\end{proof}
 The proof proceeds in the manner of
 \cite{P1}and \cite{P2}. In particular we first prove the following lemma
 which is analogous to the main theorem of \cite{P1}, see Theorem~\ref{pier1}.

 \begin{lem} Two bi-tricolored links represent the same
 $3$-manifold iff they can be related by a finite sequence of
 moves $\mathcal{M}$, $\mathcal{P}$ and
 the non-local moves I-V shown in figures \ref{I} through \ref{V}, where $\tau$,
 $\tau'$ are arbitrary tangles with the shown endpoints.
\end{lem}

\begin{figure}[htp]
\[
\begin{pspicture}(-.5,0)(5,5)
 \psdots(0,4)(.5,4)(1,4)(1.5,4)(2,4)(2.5,4)(4,4)(4.5,4)
        (0,1)(.5,1)(1,1)(1.5,1)(2,1)(2.5,1)(4,1)(4.5,1)
   \psline(1,4)(1,1)
   \psline(1.5,4)(1.5,1)
   \psline(2,4)(2,1)
   \psline(2.5,4)(2.5,1)
   \psline(4,4)(4,1)
   \psline(4.5,4)(4.5,1)
 \psline(0,4)(0,2.75)
 \psline(0,2.25)(0,1)
 \psline(0.5,4)(0.5,2.75)
 \psline(0.5,2.25)(0.5,1)
  \psline(0,2.75)(.5,2.25)
  \psline(.5,2.75)(.3,2.55)
  \psline(.2,2.44)(0,2.25)
           \psframe(-.3,0)(4.8,1)
           \psframe(-.3,4)(4.8,5)
    \rput(0,4.3){12}
 \rput(0.5,4.3){12}
 \rput(1,4.3){14}
 \rput(1.5,4.3){14}
 \rput(2,4.3){23}
 \rput(2.5,4.3){23}
 \rput(4,4.3){23}
 \rput(4.5,4.3){23}
                   \rput(3.2,4.3){23's}
 \rput(0,.7){12}
 \rput(0.5,.7){12}
 \rput(1,.7){14}
 \rput(1.5,.7){14}
 \rput(2,.7){23}
 \rput(2.5,.7){23}
 \rput(4,.7){23}
 \rput(4.5,.7){23}
                   \rput(3.2,.7){23's}
  \psdots[dotsize=0.1](3,4)(3.3,4)(3.6,4)(3,1)(3.3,1)(3.6,1)
                      (3,2.5)(3.3,2.5)(3.6,2.5)
     \rput(2.5,4.8){$\tau$}
     \rput(2.5,.2){$\tau'$}
\end{pspicture}
\hspace{2cm}
\begin{pspicture}(-.5,0)(5,5)
 \psdots(0,4)(.5,4)(1,4)(1.5,4)(2,4)(2.5,4)(4,4)(4.5,4)
        (0,1)(.5,1)(1,1)(1.5,1)(2,1)(2.5,1)(4,1)(4.5,1)
   \psline(1,4)(1,1)
   \psline(1.5,4)(1.5,1)
   \psline(2,4)(2,1)
   \psline(2.5,4)(2.5,1)
   \psline(4,4)(4,1)
   \psline(4.5,4)(4.5,1)
 \psline(0,4)(0,1)
 \psline(0.5,4)(0.5,1)
             \psframe(-.3,0)(4.8,1)
             \psframe(-.3,4)(4.8,5)
 \rput(0,4.3){12}
 \rput(0.5,4.3){12}
 \rput(1,4.3){14}
 \rput(1.5,4.3){14}
 \rput(2,4.3){23}
 \rput(2.5,4.3){23}
 \rput(4,4.3){23}
 \rput(4.5,4.3){23}
                   \rput(3.2,4.3){23's}
 \rput(0,.7){12}
 \rput(0.5,.7){12}
 \rput(1,.7){14}
 \rput(1.5,.7){14}
 \rput(2,.7){23}
 \rput(2.5,.7){23}
 \rput(4,.7){23}
 \rput(4.5,.7){23}
                   \rput(3.2,.7){23's}
  \psdots[dotsize=0.1](3,4)(3.3,4)(3.6,4)(3,1)(3.3,1)(3.6,1)
                      (3,2.5)(3.3,2.5)(3.6,2.5)
     \rput(2.5,4.8){$\tau$}
     \rput(2.5,.2){$\tau'$}
\end{pspicture}
\]
 \caption{Move I}\label{I}
\end{figure}

\begin{figure}[htp]
\[
\begin{pspicture}(-.5,0)(5,5)
 \psdots(0,4)(.5,4)(1,4)(1.5,4)(2,4)(2.5,4)(4,4)(4.5,4)
        (0,1)(.5,1)(1,1)(1.5,1)(2,1)(2.5,1)(4,1)(4.5,1)
   \psline(0,4)(0,1)
   \psline(0.5,4)(0.5,1)
   \psline(2,4)(2,1)
   \psline(2.5,4)(2.5,1)
   \psline(4,4)(4,1)
   \psline(4.5,4)(4.5,1)
 \psline(1,4)(1,2.75)
 \psline(1,2.25)(1,1)
 \psline(1.5,4)(1.5,2.75)
 \psline(1.5,2.25)(1.5,1)
  \psline(1,2.75)(1.5,2.25)
  \psline(1.5,2.75)(1.3,2.55)
  \psline(1.2,2.44)(1,2.25)
           \psframe(-.3,0)(4.8,1)
           \psframe(-.3,4)(4.8,5)
    \rput(0,4.3){12}
 \rput(0.5,4.3){12}
 \rput(1,4.3){14}
 \rput(1.5,4.3){14}
 \rput(2,4.3){23}
 \rput(2.5,4.3){23}
 \rput(4,4.3){23}
 \rput(4.5,4.3){23}
                   \rput(3.2,4.3){23's}
 \rput(0,.7){12}
 \rput(0.5,.7){12}
 \rput(1,.7){14}
 \rput(1.5,.7){14}
 \rput(2,.7){23}
 \rput(2.5,.7){23}
 \rput(4,.7){23}
 \rput(4.5,.7){23}
                   \rput(3.2,.7){23's}
  \psdots[dotsize=0.1](3,4)(3.3,4)(3.6,4)(3,1)(3.3,1)(3.6,1)
                      (3,2.5)(3.3,2.5)(3.6,2.5)
     \rput(2.5,4.8){$\tau$}
     \rput(2.5,.2){$\tau'$}
\end{pspicture}
\hspace{2cm}
\begin{pspicture}(-.5,0)(5,5)
 \psdots(0,4)(.5,4)(1,4)(1.5,4)(2,4)(2.5,4)(4,4)(4.5,4)
        (0,1)(.5,1)(1,1)(1.5,1)(2,1)(2.5,1)(4,1)(4.5,1)
   \psline(1,4)(1,1)
   \psline(1.5,4)(1.5,1)
   \psline(2,4)(2,1)
   \psline(2.5,4)(2.5,1)
   \psline(4,4)(4,1)
   \psline(4.5,4)(4.5,1)
 \psline(0,4)(0,1)
 \psline(0.5,4)(0.5,1)
             \psframe(-.3,0)(4.8,1)
             \psframe(-.3,4)(4.8,5)
 \rput(0,4.3){12}
 \rput(0.5,4.3){12}
 \rput(1,4.3){14}
 \rput(1.5,4.3){14}
 \rput(2,4.3){23}
 \rput(2.5,4.3){23}
 \rput(4,4.3){23}
 \rput(4.5,4.3){23}
                   \rput(3.2,4.3){23's}
 \rput(0,.7){12}
 \rput(0.5,.7){12}
 \rput(1,.7){14}
 \rput(1.5,.7){14}
 \rput(2,.7){23}
 \rput(2.5,.7){23}
 \rput(4,.7){23}
 \rput(4.5,.7){23}
                                    \rput(3.2,.7){23's}
  \psdots[dotsize=0.1](3,4)(3.3,4)(3.6,4)(3,1)(3.3,1)(3.6,1)
                      (3,2.5)(3.3,2.5)(3.6,2.5)
     \rput(2.5,4.8){$\tau$}
     \rput(2.5,.2){$\tau'$}
\end{pspicture}
\]
 \caption{Move II}\label{II}
\end{figure}

\begin{figure}[htp]
\[
\begin{pspicture}(-.5,0)(6.5,5)
 \psdots(0,4)(.5,4)(1,4)(1.5,4)(2,4)(2.5,4)(3,4)(3.5,4)(5.5,4)(6,4)
        (0,1)(.5,1)(1,1)(1.5,1)(2,1)(2.5,1)(3,1)(3.5,1)(5.5,1)(6,1)
   \psline(0,4)(0,1)
   \psline(0.5,4)(0.5,1)
   \psline(3,4)(3,1)
   \psline(3.5,4)(3.5,1)
   \psline(5.5,4)(5.5,1)
   \psline(6,4)(6,1)
  \psline(1,4)(2,2.5)
  \psline(1.5,4)(2.5,2.5)
  \psline(1,2.5)(2,1)
  \psline(1.5,2.5)(2.5,1)
   \psline(1,2.5)(1.4,3.1)
   \psline(1.6,3.4)(1.7,3.55)
   \psline(1.85,3.755)(2,4)
    \psline(1.5,2.5)(1.7,2.8)
    \psline(1.85,3.025)(1.95,3.175)
    \psline(2.05,3.325)(2.5,4)
    \psline(2,2.5)(1.8,2.2)
    \psline(1.65,1.975)(1.56,1.84)
    \psline(1.4,1.6)(1,1)
    \psline(2.5,2.5)(2.1,1.9)
    \psline(1.9,1.6)(1.8,1.45)
    \psline(1.7,1.3)(1.5,1)
           \psframe(-.3,0)(6.3,1)
           \psframe(-.3,4)(6.3,5)
 \rput(0,4.3){14}
 \rput(0.5,4.3){14}
 \rput(1,4.3){12}
 \rput(1.5,4.3){12}
 \rput(2,4.3){23}
 \rput(2.5,4.3){23}
 \rput(3,4.3){23}
 \rput(3.5,4.3){23}
 \rput(5.5,4.3){23}
 \rput(6,4.3){23}
                   \rput(4.35,4.3){23's}
 \rput(0,.7){14}
 \rput(0.5,.7){14}
 \rput(1,.7){12}
 \rput(1.5,.7){12}
 \rput(2,.7){23}
 \rput(2.5,.7){23}
 \rput(3,.7){23}
 \rput(3.5,.7){23}
 \rput(5.5,.7){23}
 \rput(6,.7){23}
                   \rput(4.35,.7){23's}
  \psdots[dotsize=0.1](4.2,4)(4.5,4)(4.8,4)(4.2,1)(4.5,1)(4.8,1)
                      (4.2,2.5)(4.5,2.5)(4.8,2.5)
     \rput(3.2,4.8){$\tau$}
     \rput(3.2,.2){$\tau'$}
\end{pspicture}
\hspace{2cm}
\begin{pspicture}(-.5,0)(5,5)
 \psdots(0,4)(.5,4)(1,4)(1.5,4)(2,4)(2.5,4)(4,4)(4.5,4)
        (0,1)(.5,1)(1,1)(1.5,1)(2,1)(2.5,1)(4,1)(4.5,1)
   \psline(1,4)(1,1)
   \psline(1.5,4)(1.5,1)
   \psline(2,4)(2,1)
   \psline(2.5,4)(2.5,1)
   \psline(4,4)(4,1)
   \psline(4.5,4)(4.5,1)
 \psline(0,4)(0,1)
 \psline(0.5,4)(0.5,1)
             \psframe(-.3,0)(4.8,1)
             \psframe(-.3,4)(4.8,5)
 \rput(0,4.3){14}
 \rput(0.5,4.3){14}
 \rput(1,4.3){12}
 \rput(1.5,4.3){12}
 \rput(2,4.3){23}
 \rput(2.5,4.3){23}
 \rput(4,4.3){23}
 \rput(4.5,4.3){23}
                   \rput(3.2,4.3){23's}
 \rput(0,.7){14}
 \rput(0.5,.7){14}
 \rput(1,.7){12}
 \rput(1.5,.7){12}
 \rput(2,.7){23}
 \rput(2.5,.7){23}
 \rput(4,.7){23}
 \rput(4.5,.7){23}
                   \rput(3.2,.7){23's}
  \psdots[dotsize=0.1](3,4)(3.3,4)(3.6,4)(3,1)(3.3,1)(3.6,1)
                      (3,2.5)(3.3,2.5)(3.6,2.5)
     \rput(2.5,4.8){$\tau$}
     \rput(2.5,.2){$\tau'$}
\end{pspicture}
\]
 \caption{Move III}\label{III}
\end{figure}
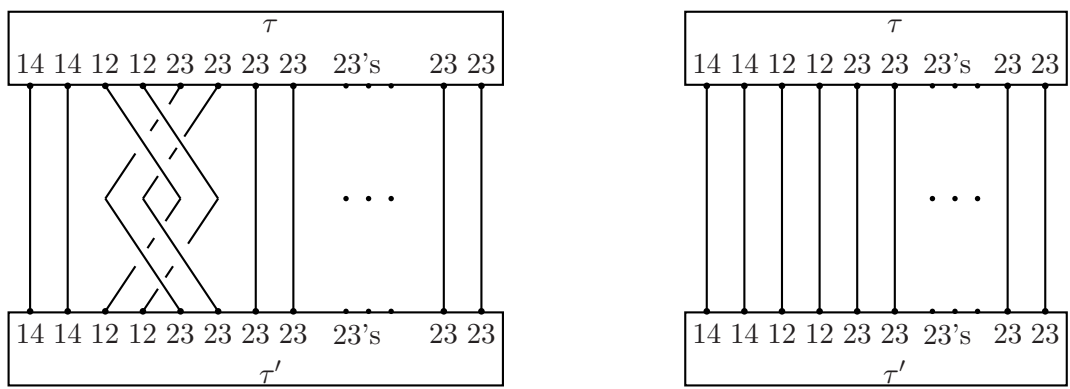

\begin{figure}[htp]
\[
\begin{pspicture}(-2,-4)(6.5,5)
 \psdots(-1,4)(-.5,4)(0,4)(.5,4)(1,4)(1.5,4)(2,4)(2.5,4)(3,4)(3.5,4)(5.5,4)(6,4)
        (-1,-3)(-.5,-3)(0,-3)(.5,-3)(1,-3)(1.5,-3)(2,-3)(2.5,-3)(3,-3)(3.5,-3)(5.5,-3)(6,-3)
   \psline(0,4)(0,1)
   \psline(0.5,4)(0.5,1)
   \psline(3,4)(3,-3)
   \psline(3.5,4)(3.5,-3)
   \psline(5.5,4)(5.5,-3)
   \psline(6,4)(6,-3)
  \psline(1,4)(2,2.5)
  \psline(1.5,4)(2.5,2.5)
  \psline(1,2.5)(2,1)
  \psline(1.5,2.5)(2.5,1)
   \psline(1,2.5)(1.4,3.1)
   \psline(1.6,3.4)(1.7,3.55)
   \psline(1.85,3.755)(2,4)
    \psline(1.5,2.5)(1.7,2.8)
    \psline(1.85,3.025)(1.95,3.175)
    \psline(2.05,3.325)(2.5,4)
    \psline(2,2.5)(1.8,2.2)
    \psline(1.65,1.975)(1.56,1.84)
    \psline(1.4,1.6)(1,1)
    \psline(2.5,2.5)(2.1,1.9)
    \psline(1.9,1.6)(1.8,1.45)
    \psline(1.7,1.3)(1.5,1)
           \psframe(-1.3,-4)(6.3,-3)
           \psframe(-1.3,4)(6.3,5)
 \rput(-1,4.3){14}
 \rput(-.5,4.3){14}
 \rput(0,4.3){12}
 \rput(.5,4.3){12}
 \rput(1,4.3){13}
 \rput(1.5,4.3){13}
 \rput(2,4.3){12}
 \rput(2.5,4.3){12}
 \rput(3,4.3){23}
 \rput(3.5,4.3){23}
 \rput(5.5,4.3){23}
 \rput(6,4.3){23}
                   \rput(4.35,4.3){23's}
 \rput(-1,-3.3){14}
 \rput(-.5,-3.3){14}
 \rput(0,-3.3){12}
 \rput(.5,-3.3){12}
 \rput(1,-3.3){13}
 \rput(1.5,-3.3){13}
 \rput(2,-3.3){12}
 \rput(2.5,-3.3){12}
 \rput(3,-3.3){23}
 \rput(3.5,-3.3){23}
 \rput(5.5,-3.3){23}
 \rput(6,-3.3){23}
                   \rput(4.35,-3.3){23's}
  \psdots[dotsize=0.1](4.2,4)(4.5,4)(4.8,4)(4.2,-3)(4.5,-3)(4.8,-3)
                      (4.2,.5)(4.5,.5)(4.8,.5)
     \rput(2.2,4.8){$\tau$}
     \rput(2.2,-3.8){$\tau'$}
   \psline(0,1)(1,-.5)
  \psline(.5,1)(1.5,-.5)
  \psline(0,-.5)(1,-2)
  \psline(.5,-.5)(1.5,-2)
   \psline(0,-.5)(.4,.1)
   \psline(.6,.4)(.7,.55)
   \psline(.85,.755)(1,1)
    \psline(.5,-.5)(.7,-.2)
    \psline(.85,.025)(.95,.175)
    \psline(1.05,.325)(1.5,1)
    \psline(1,-.5)(.8,-.8)
    \psline(.65,-1.025)(.56,-1.16)
    \psline(.4,-1.4)(0,-2)
    \psline(1.5,-.5)(1.1,-1.1)
    \psline(.9,-1.4)(.8,-1.55)
    \psline(.7,-1.7)(.5,-2)
 \psline(1,-2)(1.5,-2.5)
 \psline(1.5,-2)(1.3,-2.2)
 \psline(1.2,-2.3)(1,-2.5)
   \psline(1,-2.5)(1.5,-3)
   \psline(1.5,-2.5)(1.3,-2.7)
   \psline(1.2,-2.8)(1,-3)
  \psline(0,-2)(0,-3)
  \psline(.5,-2)(.5,-3)
  \psline(2,1)(2,-3)
  \psline(2.5,1)(2.5,-3)
  \psline(-.5,4)(-.5,-3)
  \psline(-1,4)(-1,-3)
\end{pspicture}
\hspace{2cm}
\begin{pspicture}(-.5,0)(5,5)
 \psdots(0,4)(.5,4)(1,4)(1.5,4)(2,4)(2.5,4)(4,4)(4.5,4)
        (0,1)(.5,1)(1,1)(1.5,1)(2,1)(2.5,1)(4,1)(4.5,1)
   \psline(1,4)(1,1)
   \psline(1.5,4)(1.5,1)
   \psline(2,4)(2,1)
   \psline(2.5,4)(2.5,1)
   \psline(4,4)(4,1)
   \psline(4.5,4)(4.5,1)
 \psline(0,4)(0,1)
 \psline(0.5,4)(0.5,1)
             \psframe(-.3,0)(4.8,1)
             \psframe(-.3,4)(4.8,5)
 \rput(0,4.3){14}
 \rput(0.5,4.3){14}
 \rput(1,4.3){12}
 \rput(1.5,4.3){12}
 \rput(2,4.3){13}
 \rput(2.5,4.3){13}
 \rput(4,4.3){23}
 \rput(4.5,4.3){23}
                   \rput(3.2,4.3){23's}
 \rput(0,.7){14}
 \rput(0.5,.7){14}
 \rput(1,.7){12}
 \rput(1.5,.7){12}
 \rput(2,.7){13}
 \rput(2.5,.7){13}
 \rput(4,.7){23}
 \rput(4.5,.7){23}
                   \rput(3.2,.7){23's}
  \psdots[dotsize=0.1](3,4)(3.3,4)(3.6,4)(3,1)(3.3,1)(3.6,1)
                      (3,2.5)(3.3,2.5)(3.6,2.5)
     \rput(2.5,4.8){$\tau$}
     \rput(2.5,.2){$\tau'$}
\end{pspicture}
\]
\caption{Move IV}\label{IV}
\end{figure}

\begin{figure}[htp]
\[
\begin{pspicture}(8,10)
 \psline(4.5,8)(1.5,7)
 \psline(1.5,4)(4.5,3)
 \psline(2,9)(2,7.25)
 \psline(2,7.05)(2,3.95)
 \psline(2,3.7)(2,1.833)
 \psline(3,9)(3,7.6)
 \psline(3,7.4)(3,3.6)
 \psline(3,3.4)(3,1.833)
 \psline(4,9)(4,7.95)
 \psline(4,7.7)(4,3.3)
 \psline(4,3.05)(4,1.833)
 \psline(1,9)(1.9,8.746)
 \psline(2.1,8.68)(2.9,8.45)
 \psline(3.1,8.39)(3.9,8.16)
 \psline(4.1,8.1)(4.5,8)
 \psline(1.5,7)(1.9,6.867)
 \psline(2.1,6.8)(2.9,6.533)
 \psline(3.1,6.467)(3.9,6.2)
 \psline(4.1,6.133)(4.5,6)
 \psline(4.5,6)(5,5.5)
 \psline(5,5.5)(5,1.833)
 \psline(5,9)(5,6.2)
 \psline(5,6.2)(4.75,5.9)
 \psline(4.62,5.744)(4.1,5.12)
 \psline(3.9,4.96)(3.1,4.64)
 \psline(2.9,4.56)(2.1,4.24)
 \psline(1.9,4.16)(1.5,4)
 \psline(4.5,3)(4.1,2.867)
 \psline(3.9,2.8)(3.1,2.533)
 \psline(2.9,2.467)(2.1,2.2)
 \psline(1.9,2.133)(1.5,2)
 \psline(1.5,2)(1,1.833)
 \psdots(0,9)(1,9)(2,9)(3,9)(4,9)(5,9)(6,9)(7.5,9)
        (0,1)(1,1)(2,1)(3,1)(4,1)(5,1)(6,1)(7.5,1)
      \psline(1,1.833)(1,1)
      \psline(2,1.833)(2,1)
      \psline(3,1.833)(3,1)
   \psline(0,9)(0,1)
   \psline(6,9)(6,1)
   \psline(7.5,9)(7.5,1)
     \psline(5,1.833)(4,1)
     \psline(4,1.833)(4.4,1.4998)
     \psline(4.6,1.3332)(5,1)
 \rput(0,9.3){12}
 \rput(6,9.3){23}
 \rput(7.5,9.3){23}
 \rput(1,9.3){12}
 \rput(2,9.3){14}
 \rput(3,9.3){14}
 \rput(4,9.3){23}
 \rput(5,9.3){23}
                 \psdots[dotsize=.1](6.5,9)(6.8,9)(7.1,9)
                 \rput(6.8,9.3){23's}
\rput(0,.7){12}
 \rput(6,.7){23}
 \rput(7.5,.7){23}
 \rput(1,.7){12}
 \rput(2,.7){14}
 \rput(3,.7){14}
 \rput(4,.7){23}
 \rput(5,.7){23}
                 \psdots[dotsize=.1](6.5,1)(6.8,1)(7.1,1)
                 \rput(6.8,.7){23's}
       \psframe(-.3,9)(7.8,10)
       \psframe(-.3,0)(7.8,1)
      \rput(3.5,9.8){$\tau$}
      \rput(3.5,.2){$\tau'$}
                 \psdots[dotsize=.1](6.5,4.5)(6.8,4.5)(7.1,4.5)
\end{pspicture}
\hspace{2cm}
\begin{pspicture}(-.5,0)(5,5)
 \psdots(0,4)(.5,4)(1,4)(1.5,4)(2,4)(2.5,4)(4,4)(4.5,4)
        (0,1)(.5,1)(1,1)(1.5,1)(2,1)(2.5,1)(4,1)(4.5,1)
   \psline(1,4)(1,1)
   \psline(1.5,4)(1.5,1)
   \psline(2,4)(2,1)
   \psline(2.5,4)(2.5,1)
   \psline(4,4)(4,1)
   \psline(4.5,4)(4.5,1)
 \psline(0,4)(0,1)
 \psline(0.5,4)(0.5,1)
             \psframe(-.3,0)(4.8,1)
             \psframe(-.3,4)(4.8,5)
 \rput(0,4.3){12}
 \rput(0.5,4.3){12}
 \rput(1,4.3){14}
 \rput(1.5,4.3){14}
 \rput(2,4.3){23}
 \rput(2.5,4.3){23}
 \rput(4,4.3){23}
 \rput(4.5,4.3){23}
                   \rput(3.2,4.3){23's}
 \rput(0,.7){12}
 \rput(0.5,.7){12}
 \rput(1,.7){14}
 \rput(1.5,.7){14}
 \rput(2,.7){23}
 \rput(2.5,.7){23}
 \rput(4,.7){23}
 \rput(4.5,.7){23}
                                    \rput(3.2,.7){23's}
  \psdots[dotsize=0.1](3,4)(3.3,4)(3.6,4)(3,1)(3.3,1)(3.6,1)
                      (3,2.5)(3.3,2.5)(3.6,2.5)
     \rput(2.5,4.8){$\tau$}
     \rput(2.5,.2){$\tau'$}
\end{pspicture}
\]
\caption{Move V}\label{V}
\end{figure}

\begin{proof} That moves I, II and V   do not change $M(L)$ follows from the
fact that the braid in their right side is in the kernel of $\bar{\lambda}$.
That the moves III and IV do not
change $M(L)$ follows from the fact that they can be realized by $\mathcal{M}$, $\mathcal{P}$
and $\mathcal{X}$
as shown in the proof of the theorem \ref{theoremX}.
(There is no vicious circle as this part of the argument is
not used in that proof).\

The proof of the reverse direction is completely analogous to Piergallini's proof
of Theorem~\ref{pier1}.
Referring to the sketch of that proof given in Section~\ref{3sphere3} we just
comment on how each step of that
proof goes through in the present situation:

\begin{itemize}
\item Each bi-tircolored link has a normalized diagram.
\item Step $1$ (Heegard stabilization) can be realized in exactly the same way.
\item For step $2$, in order to get braids that lift to the Suzuki generators
one just has to add two trivial strands colored by $(14)$ to the braids Piergallini uses.
\item For step $3$, we have to show how to add the normal generators of $Ker\bar{\lambda}$
given in
Theorem~\ref{kernel} to the top (and bottom) of normalized diagrams using the moves.\
$\B_0$ and $\B_2$ can be obviously added using moves I and II respectively. Also by
slight modification of
Piegallini's proof (by just adding two trivial strands colored by $(14)$) $B_1$
and $D_1$ can be added using
moves $\mathcal{M}$, I-IV.\
Now move V transforms $\delta_4$ to $\B_4$
 and therefore using move V  we can transform the remaining
  normal generators of $Ker\bar{\lambda}$ to $B$ or $D$
(it easily checked that the colors are right).
     Thus any element of
$Ker\bar{\lambda}$ can be added on the top and the bottom of a normalized
diagram using the given moves.

\end{itemize}

\end{proof}

\begin{proof}[Proof of \ref{theoremX}]Moves II and V do not change $M(L)$
since the braids on their left side belong in the kernel of the
lifting homomorphism $\bar\lambda$.
 That the move $\mathcal{X}$
does not change $M(L)$ follows from the fact that it can be
implemented using moves $\mathcal{M}$, $\mathcal{P}$ and addition
of a trivial fifth sheet as shown in the proof of theorem
\ref{stab}. (Again there is no vicious circle.)\
 To show the reverse it suffices to show that the
moves I through V of the previous theorem can be realized using those moves.
 This is done in the next few pages.

\newpage

For move I:

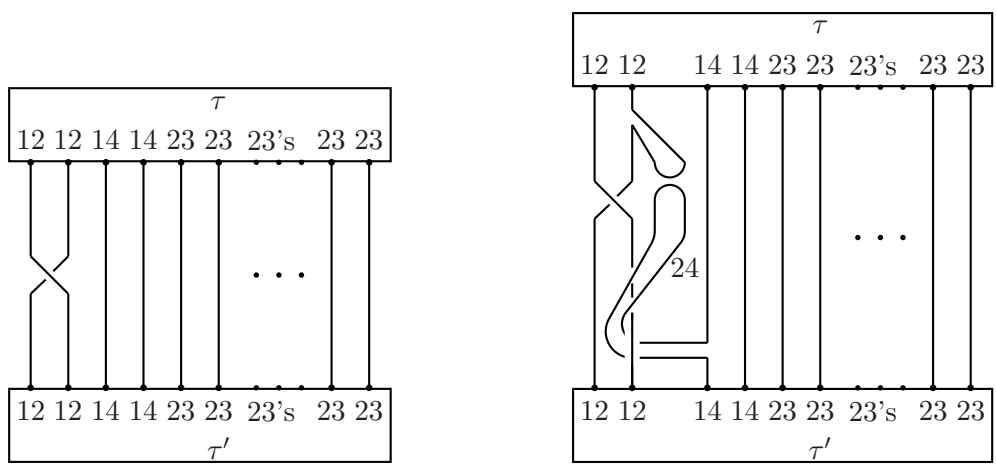
\[

\]

by isotopy. Then by reverse circumcision one gets:

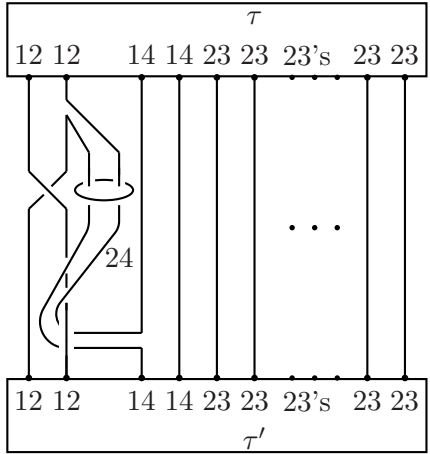
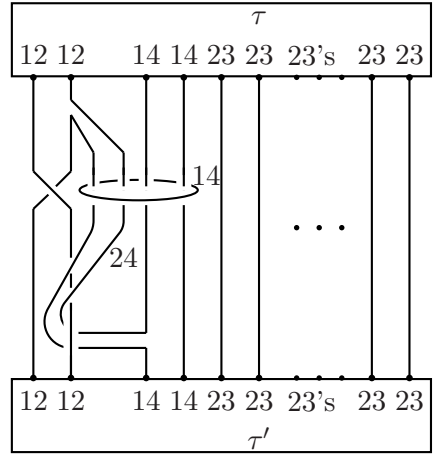
\[
%
\]

by using  move $\mathcal{X}$. Then by repeatedly using move $\mathcal{P}$ one gets:

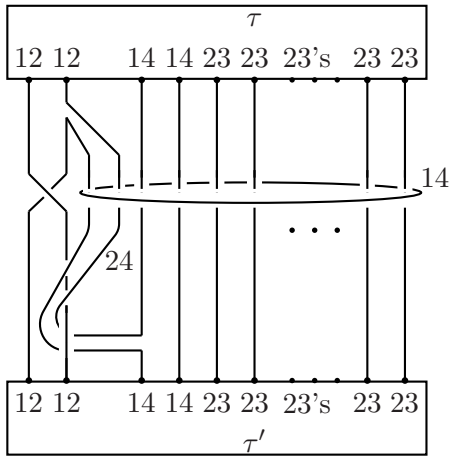
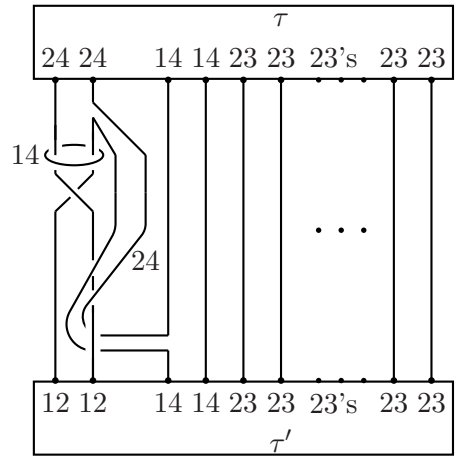
\[
%
\]

 by isotopically flipping the $(14)$-circle over the whole diagram.
Now circumcise to get:

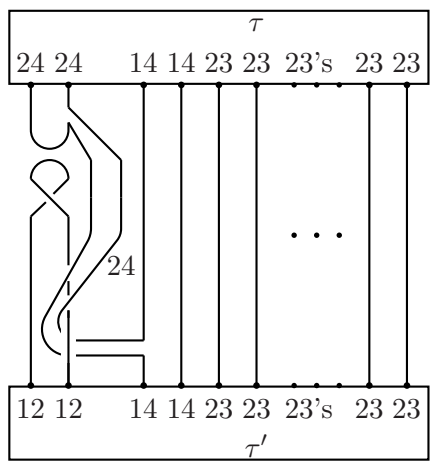
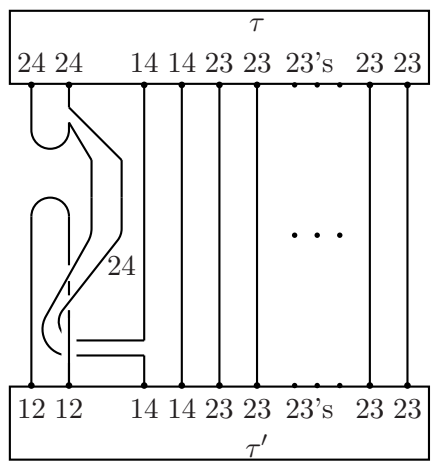
\[
%
\]

 by isotopy. Now that there is no crossing reverse the process
starting from the previous picture to complete move I.\\

 \newpage

 For move III:  Use isotopy and reverse circumcision to get:
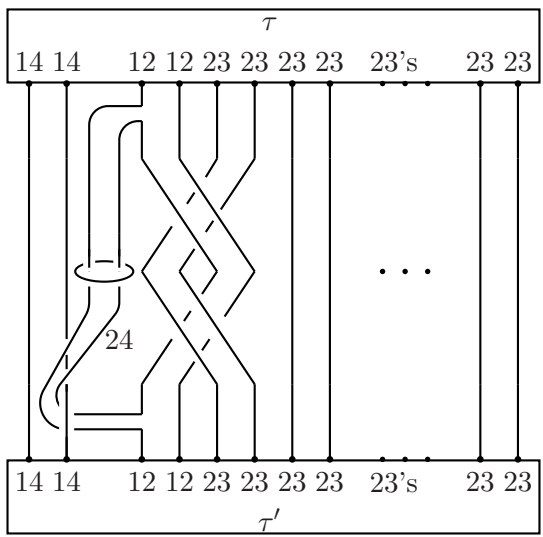
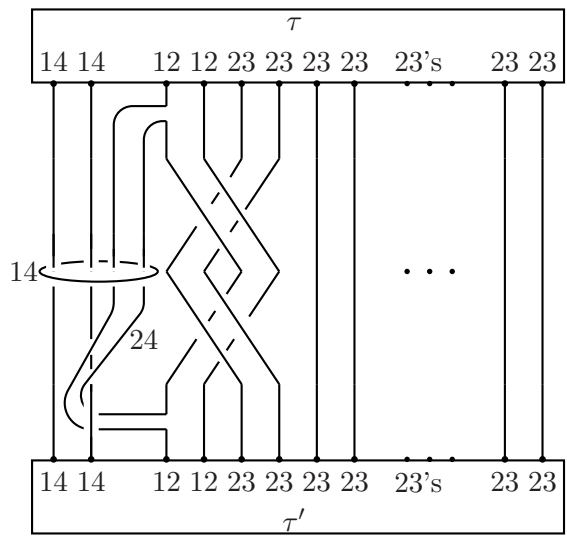
\[
\hspace{-.5cm}
%
\]

by using move $\mathcal{X}$. Now isotopically flip the $(14)$-circle over
 the whole diagram and apply
$\mathcal{P}$ moves to get it around the $(12)$ strip and then circumcise to get:

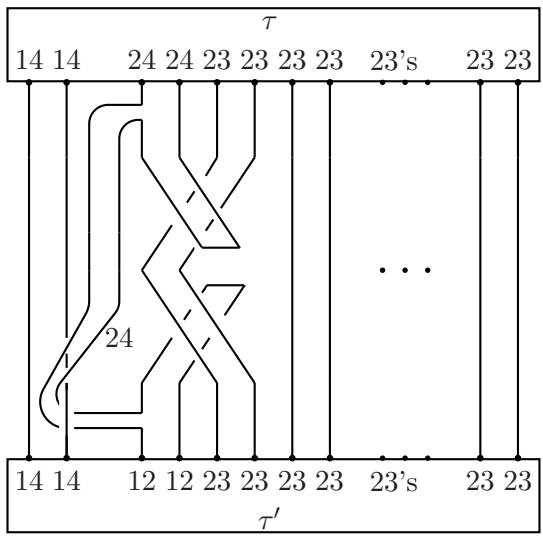
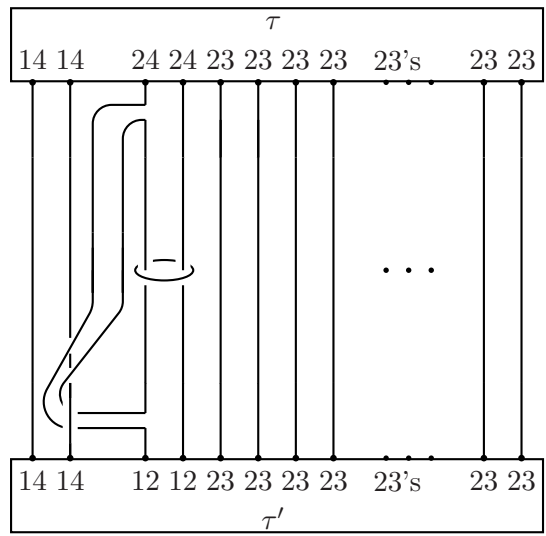
\[
\hspace{-.5cm}
%
\]

by isotopy and reverse circumcision. Now after moving the $(14)$ circle
beyond the $(23)$ strands using move
$\mathcal{P}$, flip it over the whole diagram,
 apply move $\mathcal{X}$, circumcise and isotope to complete move
III.

\newpage

For move IV: Start again with isotopy and reverse circumcision to get:

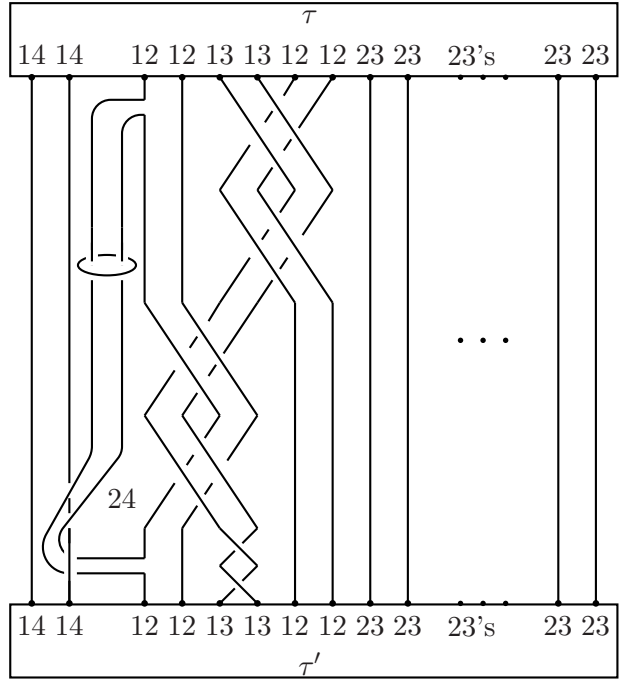
\[
%
\]

Then use move $\mathcal{X}$, isotopy and $\mathcal{P}$ moves to get:

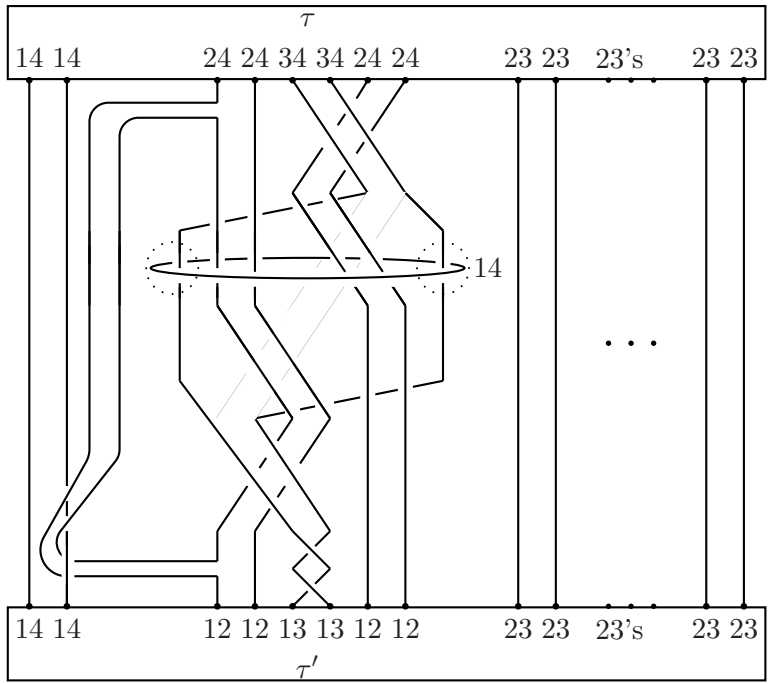
\[
%
\]

 \newpage

Then perform $\mathcal{P}$ moves inside the dotted circles to get:

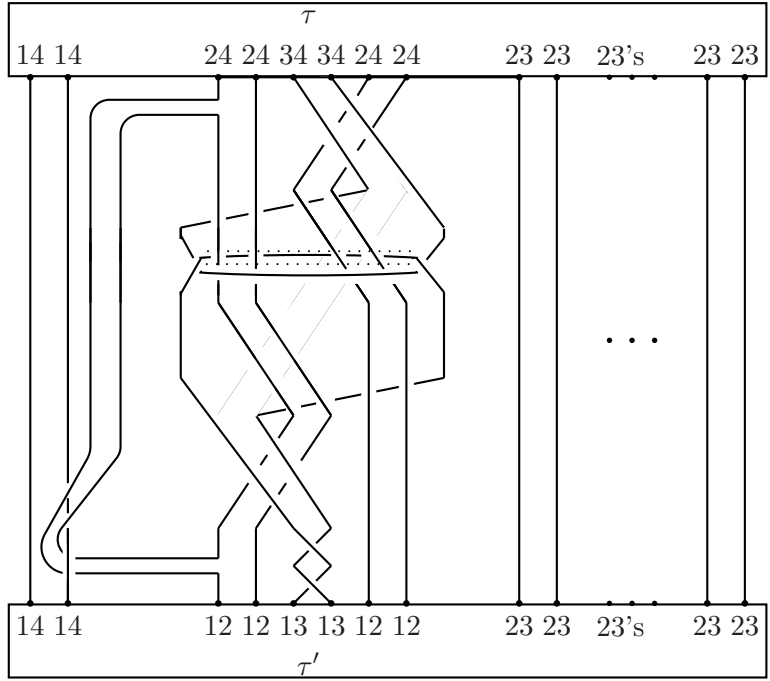
\[
%
\]
 Then after changing all crossings inside the dotted rectangle
 using move $\mathcal{P}$ isotope:
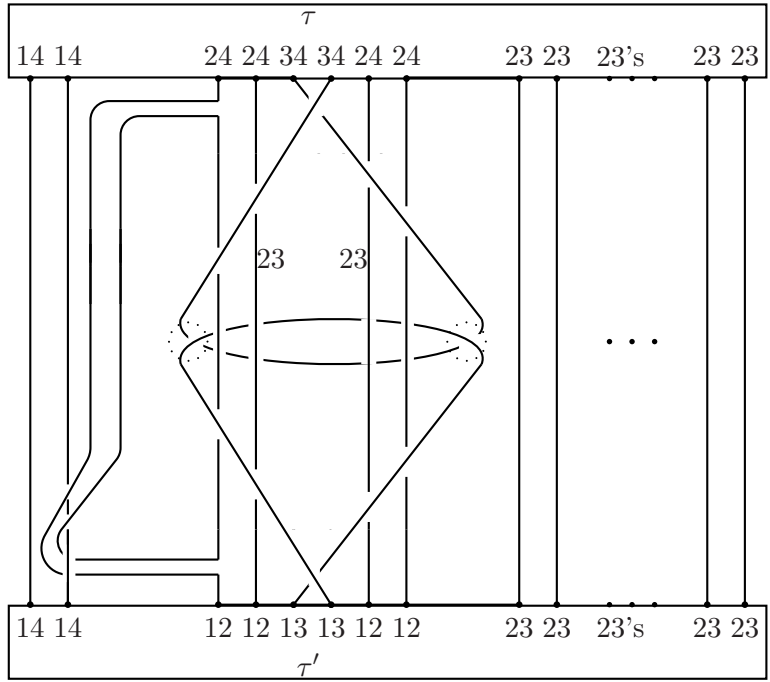
\[
%
\]

\newpage

Perform $\mathcal{P}$ moves inside the dotted circles:
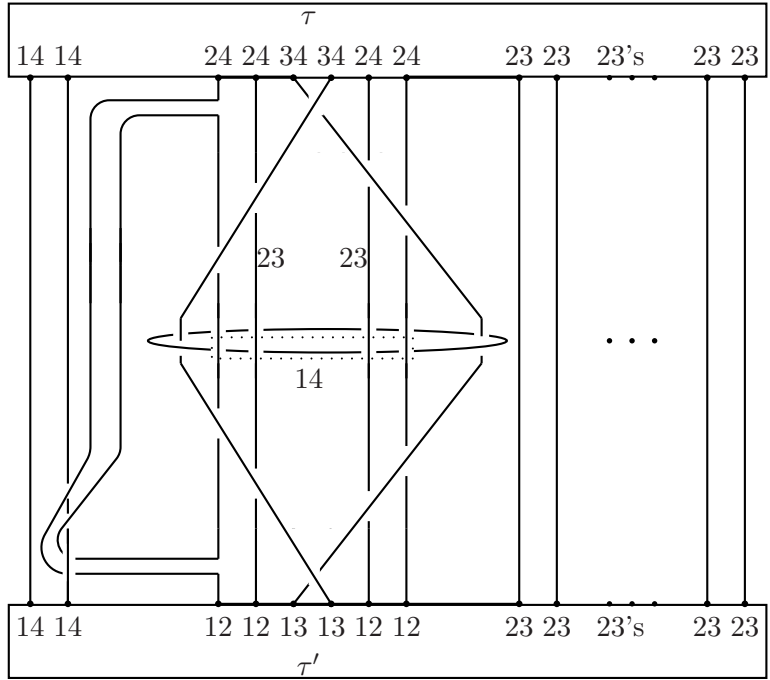
\[
%
\]

Perform $\mathcal{P}$ moves inside the dotted rectangle:

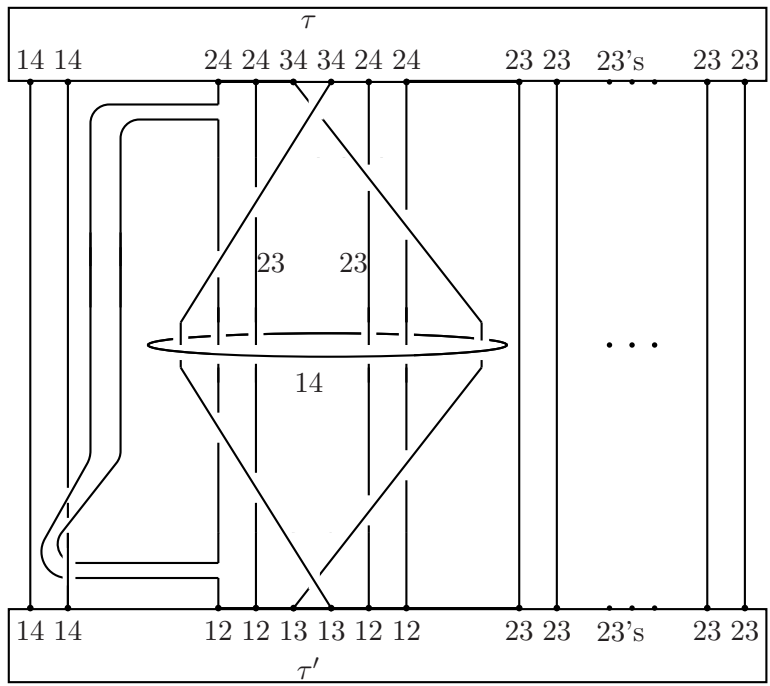
\[
%
\]
 Now reverse the process starting from the second picture to
 complete move V.\end{proof}

\newpage

\begin{proof}[Proof of \ref{stab}] It suffices to show that
moves $\mathcal{X}$, II and IV can be realized:\\

For move $\mathcal{X}$, first add a trivial fifth sheet with monodromy $(45)$ to get:\\

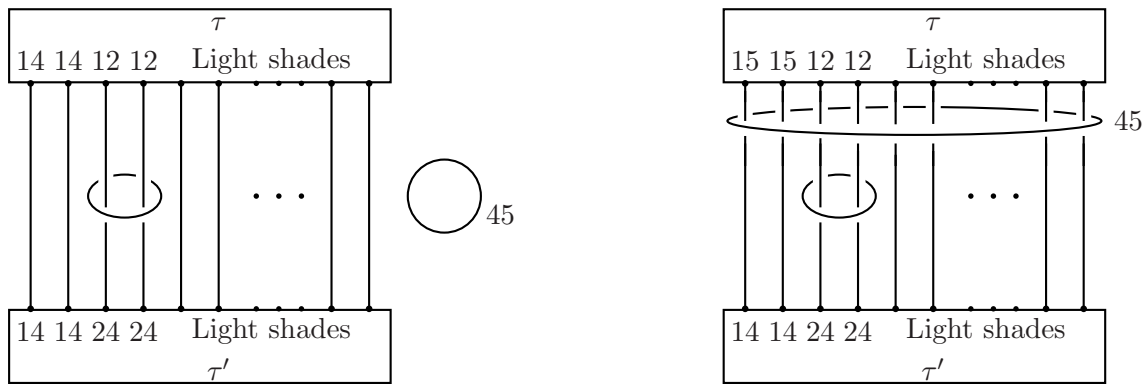
\[
%
\]

by isotopically flipping the trivial sheet over the whole diagram.
Then apply a sequence of $\mathcal{P}$ moves to get:\\

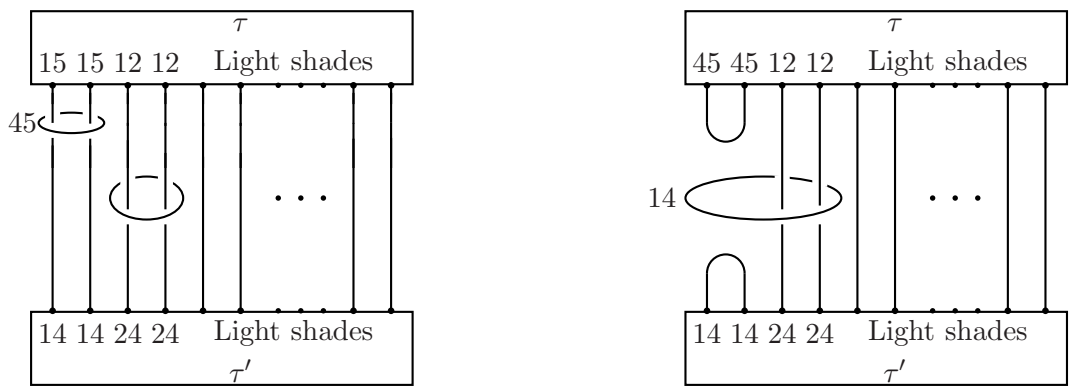
\[
%
\]

 by circumcision and isotopy. Now isotopically push the $(14)$ cup
 through the $(14)$ circle and reverse the process to complete
 move $\mathcal{X}$.

\newpage

 For move II: Add a trivial fifth sheet with monodromy $(45)$
and after flipping it over the whole diagram, use $\mathcal{P}$
moves and circumcision to get:\\

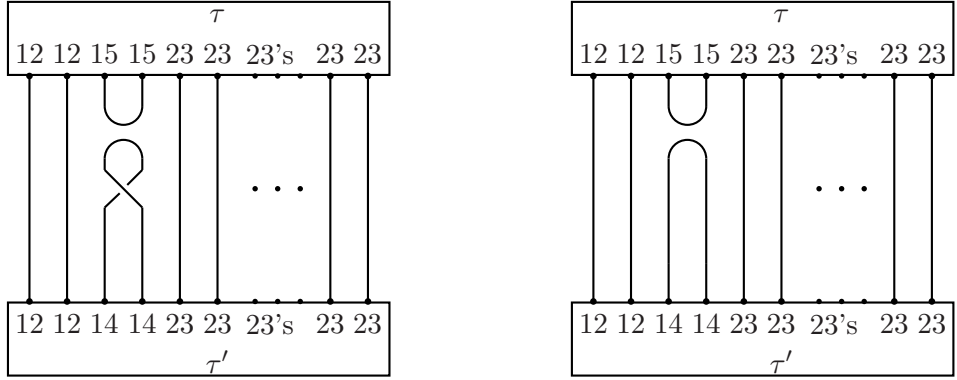
\[
%
\]

by isotopy. Now reverse the process to complete move II.

 \newpage

For move V:\\

 First add a trivial fifth sheet with monodromy $(45)$, and then use isotopy and several
 $\mathcal{P}$ moves to get:

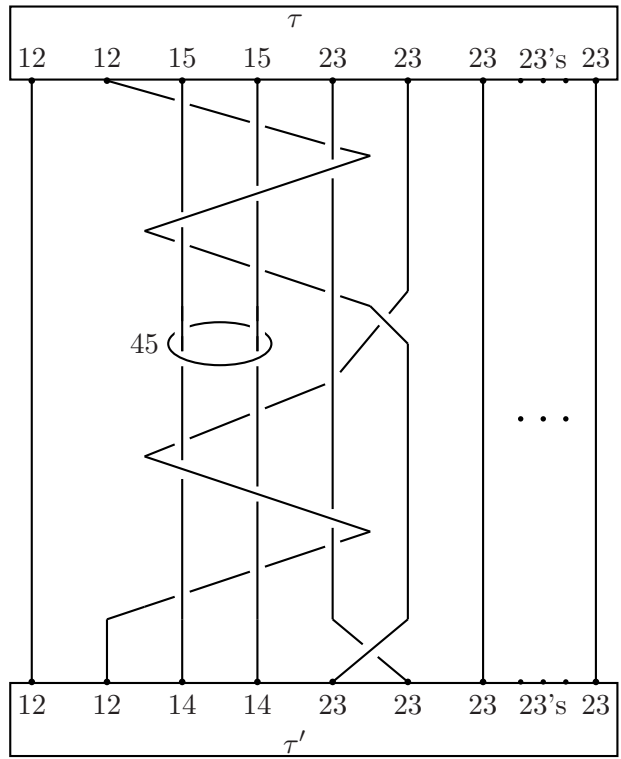
\[
%
\]

Then circumcise and isotope to get:

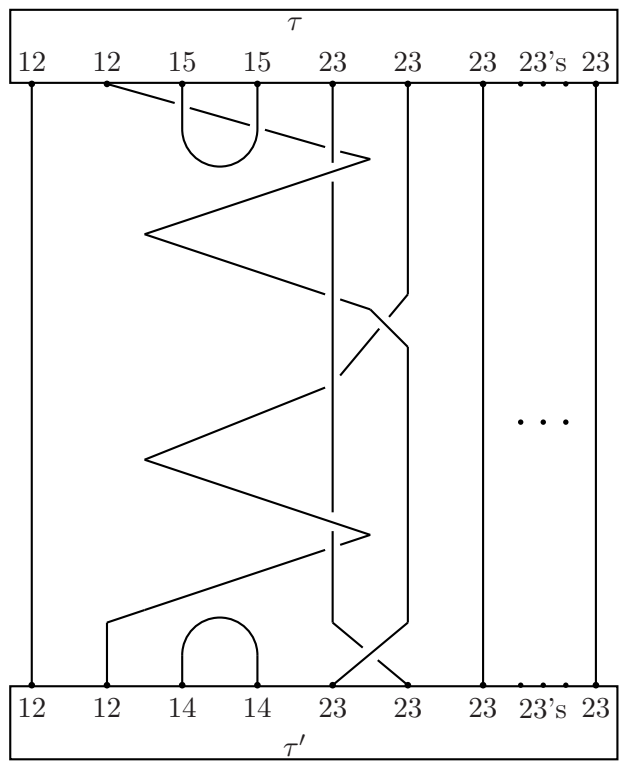
\[
%
\]

Then apply two $\mathcal{M}$ moves and isotopy to get:

 \[
%
\]

 From the last diagram apply reverse circumcision, $\mathcal{P}$
 moves, isotopy and erase the trivial sheet to complete move V.\\

  This completes the proof of the theorem.
\end{proof}

\newpage
\section{Some Open Questions}\label{questions}

This work naturally suggests the following questions:\

\begin{itemize}
\item[{\bf 1.}] Is the fifth trivial sheet really necessary or are moves
 $\mathcal{M}$ and $\mathcal{P}$ enough by
themselves? Equivalently can any bi-tricolored link be related to
 a tricolored link plus a trivial sheet, via
moves $\mathcal{M}$ and $\mathcal{P}$?
\item[{\bf 2.}] More generally is there a degree $m$ such that  moves
$\mathcal{M}$ and $\mathcal{P}$ are enough to related any two manifestations
of the same $3$-manifold as a simple
$m$-sheeted branched covering of the three sphere?
\item[{\bf 3.}] \cite{P3} Are moves $\mathcal{M}$ and $\mathcal{P}$
plus addition/deletion of trivial sheets enough
to relate any two manifestations of the same $3$-manifold as a
simple branched covering of the three sphere? Or,
referring to the move of adding/deleting a trivial sheet
as stabilization, are any two colored link presentations
of the same $3$-manifold stably equivalent?

\end{itemize}

\newpage
\appendix
\section{Wajnryb's presentation of $\CM_g$}\label{Pres}
 Let $\Sigma$ be an orientable surface. By the {\em mapping class
 group} $\CM(\Sigma)$ of $\Sigma$ we mean the group of isotopy
 classes of orientation preserving homeomorphisms of $\Sigma$.
$\CM(\Sigma)$ is indeed a group since it is the quotient of the
group of orientation preserving homeomorphisms of $\Sigma$ by the
normal subgroup of homeomorphisms isotopic to identity.

\begin{defn}\label{Dehn} Let $\alpha$ be a simple closed curve in
$\Sigma$. A Dehn twist around $\alpha$ is a homeomorphism of $\Sigma$ that
is identity outside a tubular
neighborhood of $\alpha$ and it is a full twist inside that neighborhood.
More precisely if the tubular
neighborhood of $\alpha$ is parameterized as
$\{(r,\theta)|1\leq r\leq 2,\quad \theta\in [0,2\pi] \}$ with
$\alpha$ being $\{r=1\}$ then the map inside the tubular neighborhood is given by
$(r,\theta)\mapsto(r,\theta+2\pi r)$ (see Figure~\ref{Dehntwist}).
\begin{figure}[htp]
\[
\begin{pspicture}(0,-1)(3,3)
 \pscircle(1.5,1.5){.5}
  \pscircle(1.5,1.5){1}
 \pscircle(1.5,1.5){1.5}
 \psline[linestyle=dashed](0,1.5)(1,1.5)
 \rput(2.65,1.5){$\alpha$}
 \rput(1.5,-.5){BEFORE}
 \end{pspicture}
 \hspace{2cm}
 \begin{pspicture}(0,-1)(3,3)
 \pscircle(1.5,1.5){.5}
  \pscircle(1.5,1.5){1}
 \pscircle(1.5,1.5){1.5}
 \pscurve[linestyle=dashed](0,1.5)(1.5,.2)(2.8,1.5)(1.5,2.4)(1,1.5)
 \rput(1.5,-.5){AFTER}
 \end{pspicture}
 \]
\caption{A Dehn twist around $\alpha$}\label{Dehntwist}
\end{figure}
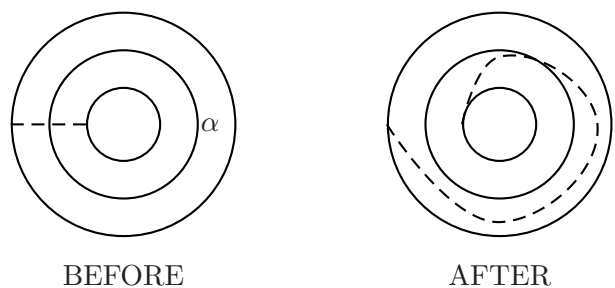
The isotopy class of a Dehn twist around $\alpha$ depends only on
the isotopy class of $\alpha$ and is called
{\em the} Dehn twist around $\alpha$. \\
Curves in a surface are denoted by Greek letters and Dehn twists
around them by the corresponding Latin letters.
\end{defn}

The mapping class group of the standard genus $g$ surface is
denoted by $\CM_g$. Wajnryb gave the following presentation of
$\CM_g$ in~\cite{W}.

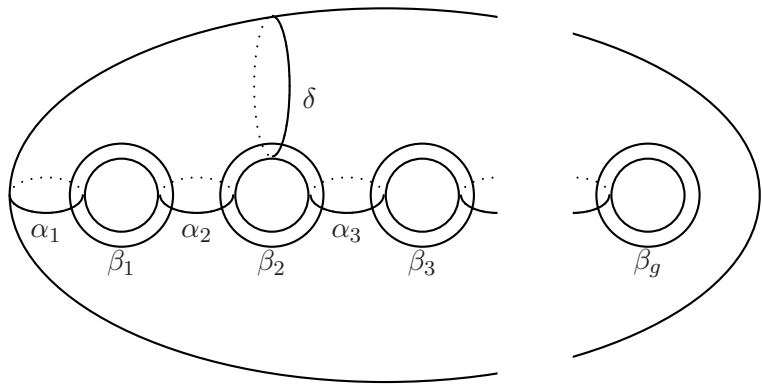
\begin{figure}[htp]
\begin{pspicture}(10,5)
 \psellipse(.5,2.5)(.5,.25)
 \psellipse(2.5,2.5)(.5,.25)
 \psellipse(4.5,2.5)(.5,.25)
 \psellipse(6.5,2.5)(.5,.25)
 \psellipse(7.5,2.5)(.5,.25)
 \psellipse(3.5,3.95)(.25,.95)
 \psframe[linecolor=white,fillstyle=solid,fillcolor=white](0,2.5)(8,3)
 \psframe[linecolor=white,fillstyle=solid,fillcolor=white](3.5,3)(3,5)
  \psellipse[linestyle=dotted](.5,2.5)(.5,.25)
 \psellipse[linestyle=dotted](2.5,2.5)(.5,.25)
 \psellipse[linestyle=dotted](4.5,2.5)(.5,.25)
 \psellipse[linestyle=dotted](7.5,2.5)(.5,.25)
 \psellipse[linestyle=dotted](3.5,3.95)(.25,.95)
 \psellipse[linestyle=dotted](6.5,2.5)(.5,.25)
 \psellipse(5,2.5)(5,2.5)
 \pscircle(1.5,2.5){.5}
 \pscircle(3.5,2.5){.5}
 \pscircle(5.5,2.5){.5}
 \pscircle(8.5,2.5){.5}
 \pscircle(1.5,2.5){.7}
 \pscircle(3.5,2.5){.7}
 \pscircle(5.5,2.5){.7}
 \pscircle(8.5,2.5){.7}
 \psframe[linecolor=white,fillstyle=solid,fillcolor=white](6.5,-.5)(7.5,5.5)
 \psdots[dotsize=.03](6.8,2.5)(7,2.5)(7.2,2.5)
 \rput(.5,2){$\alpha_1$}
 \rput(2.5,2){$\alpha_2$}
 \rput(4.5,2){$\alpha_3$}
 \rput(1.5,1.6){$\B_1$}
 \rput(3.5,1.6){$\B_2$}
 \rput(5.5,1.6){$\B_3$}
 \rput(8.5,1.6){$\B_g$}
 \rput(4,3.8){$\delta$}
\end{pspicture}
\caption{The generators of $\CM_g$}\label{Mg}
\end{figure}

\begin{thm}[\bf Wajnryb] The mapping class group of $\CM_g$ admits a
presentation with generators $a_1$, $b_1$,$\dotsc,a_n$, $b_n$, $d$
(see Figure~\ref{Mg}) and relations
\begin{itemize}
\item[(A)] $a_ib_ia_i=b_ia_ib_i,\quad
a_{i+1}b_ia_{i+1}=b_ia_{i+1}b_i,\quad b_2db_2=db_2d,$ and every
other pair of generators commute.
\item[(B)]
$(a_1b_1a_2)^4=d(b_2a_2b_1a_1a_1b_1a_2b_2)^{-1}db_2a_2b_1a_1a_1b_1a_2b_2$.
\item[(C)]
$dt_2dt_2^{-1}t_1t_2d(a_1a_2a_3t_1t_2)^{-1}=(ub_1a_2b_2a_3b_3)^{-1}vub_1a_2b_2a_3b_3\quad$
where $\quad
t_1=b_1a_1a_2b_1,\quad t_2=b_2a_2a_3b_2,\quad
u=a_3b_3t_2d(a_3b_3t_2)^{-1},$\\
$ v=a_1b_1a_2b_2d(a_1b_1a_2b_2)^{-1}.$
\item[(D)] $d_n$ commutes with $b_na_n\dotsm b_1\dotsm a_nb_n$
where
$$d_n=(u_1u_2\dotsm u_{n-1})^{-1}a_1u_1u_2\dotsm u_{n-1},$$
$$u_i=b_ia_{i+1}b_{i+1}v_i(b_{i+1}a_{i+1}b_ia_i)^{-1}\qquad
\text{for }i=1,\dotsc,n-1,$$
$$v_1=(b_2a_2b_1a_1a_1b_1a_2b_2)^{-1}d(b_2a_2b_1a_1a_1b_1a_2b_2),$$
$$v_i=t_{i-1}t_iv_{i-1}(t_{i-1}t_i)^{-1}\qquad
\text{for }i=2,\dotsc,n-1,$$
$$t_i=b_ia_ia_{i+1}b_i\qquad
\text{for }i=1,\dotsc,n-1.$$
\end{itemize}

\end{thm}

\section{Braids}\label{braids}
References for the material in this section are \cite{B} and
\cite{BZ}.

\begin{defn} Let $\text{Conf}_n(D^2)$ denote the configuration space of
 unordered $n$-tuples of distinct points of the $2$-dimensional disc and chose
 a basepoint $L=\{A_0,\dotsc,A_{n-1}\}$. The fundamental group of  $\text{Conf}_n(D^2)$
 is called  the  braid group on $n$ strands and is denoted by
 $B_n$, that is
 $$B_n:=\pi_1(\text{Conf}_n(D^2),L)  .$$
\end{defn}

 So a braid is (the homotopy class of) a closed path in $\text{Conf}_n(D^2)$.
 Such a path is equivalent to  $n$ nonintersecting paths in the disc
 with the property that the set of initial points and the
 set of endpoints are equal to $L$.
 We represent such an $n$-tuple of paths  by $n$
 arcs in the solid cylinder $I\times D^2$, which are the graphs of
 the paths thought of as functions $\xymatrix@1{
 {I}\ar[r]&{D^2}}$. We draw generic planar projections of braids
 with conventions analogous to the conventions for drawing links
  (see Appendix~\ref{diagrams}).

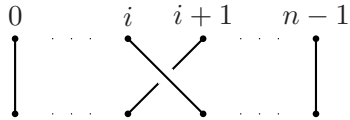
\begin{figure}[htp]
\[
\begin{pspicture}(5,2)
\psdots(.5,0)(2,0)(3,0)(4.5,0)
 \psdots(.5,1)(2,1)(3,1)(4.5,1)
 \psdots[dotsize=.01](1,1)(1.25,1)(1.5,1)(1,0)(1.25,0)(1.5,0)
 (3.5,1)(3.75,1)(4,1)(3.5,0)(3.75,0)(4,0)
 \psline(.5,1)(.5,0) \psline(4.5,1)(4.5,0)  \psline(3,1)(2,0)
\psframe[linecolor=white,fillstyle=solid,fillcolor=white](2.4,.4)(2.6,.6)
\psline(2,1)(3,0)
 \rput(.5,1.3){$0$}
 \rput(2,1.3){$i$}
 \rput(3,1.3){$i+1$}
 \rput(4.5,1.3){$n-1$}
\end{pspicture}
\]
\caption{The generator $\B_i$ of $B_n$}\label{braidgen}
\end{figure}

\begin{thm}[\bf Artin] For each $n\geq1$ the  braid group on $n$ strands has the following
 presentation :
$$B_n=\langle\B_0,\B_1,\dotsc,\B_{n-2}|\B_i\B_{i+1}\B_i=\B_{i+1}\B_i\B_{i+1}
\text{ and } \B_i\B_j=\B_j\B_i\text{ if }|i-j|\geq 2\rangle,$$
where $\B_i$ interchanges the $i^{\text{th}}$ and $(i+1)^{\text{th}}$
strands as shown in Figure~\ref{braidgen}.
\end{thm}

\section{Link diagrams}\label{diagrams}

\begin{defn} A {\em link} $L$ is a $1$-dimensional submanifold of the three
sphere $S^3$. A {\em knot} is a connected link.\\
 Two links $L$, $L'$ are called isotopic if there is
an ambient isotopy carrying one to the other, that is there is a
homotopy $\xymatrix@1{ {h_t:S^3}\ar[r]&{S^3}}$, $0\leq t\leq 1$
with the properties:
\begin{itemize}
\item $h_0=$id
\item $\forall t$, $h_t$ is a homeomorphism.
\item $L'=h_1(L)$
\end{itemize}
The fundamental group of the complement of $L$ is called the group
of the link $L$.
\end{defn}

 Isotopic links have isomorphic groups, and furthermore an isotopy
 between links
 induces an isomorphism between link groups.

\begin{figure}[htp]
\[
\begin{xy}/r9mm/:
 \vover
 \end{xy}
 \]
\caption{A crossing}\label{cross}
\end{figure}
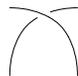

 A link is usually represented by
drawing its image under a generic (i.e. with only double points)
projection onto a plane. Furthermore at each double point of such
a diagram the relative distance of the two points from the plane
of projection is recorded by drawing the nearest piece broken as
in Figure~\ref{cross}. (See any book in knot theory, for example
\cite{R}.) The double points of such a diagram (called {\em link
diagram}) are called {\em crossings}. A link diagram is thus
broken into connected pieces called the {\em arcs} of the diagram.

\begin{figure}[htp]
\[
\begin{pspicture}(0,-1)(8,8)
 \rput(1,7.5){\begin{xy}/r1cm/: \vcap \vcross \end{xy}}
 \psline{<->}(2.5,7)(4.5,7)
 \rput(3.5,7.2){R1}
 \rput(6,7){\begin{xy}/r.9cm/:  \vuncross \end{xy}}
 \psframe[linecolor=white,fillstyle=solid,fillcolor=white](5,7)(7,8)
 \rput(3.5,6){a) First}
 \psarc[linewidth=.01](1,4.4){.45}{0}{180}
 \psframe[linecolor=white,fillstyle=solid,fillcolor=white](.6,4.6)(.7,4.7)
 \psframe[linecolor=white,fillstyle=solid,fillcolor=white](1.3,4.6)(1.4,4.7)
 \psarc[linewidth=.01](1,4.9){.45}{180}{0}
 \psline{<->}(2.5,4.6)(4.5,4.6)
  \rput(3.5,4.8){R2}
 \rput(6,4.6){\begin{xy}/r.9cm/:  \vuncross \end{xy}}
 \rput(3.5,3.6){b) Second}
 \psline[linewidth=.01](1,3)(0,2)
 \psline[linewidth=.01](2,2)(1,1)
 \psline[linewidth=.01](0,0)(1,1)
 \psframe[linecolor=white,fillstyle=solid,fillcolor=white](.4,.4)(.6,.6)
 \psframe[linecolor=white,fillstyle=solid,fillcolor=white](.4,2.4)(.6,2.6)
 \psframe[linecolor=white,fillstyle=solid,fillcolor=white](1.4,1.4)(1.6,1.6)
\psline[linewidth=.01](0,3)(2,1)
 \psline[linewidth=.01](2,1)(2,0)
 \psline[linewidth=.01](0,2)(0,1)
 \psline[linewidth=.01](0,1)(1,0)
 \psline[linewidth=.01](2,3)(2,2)
 \psline{<->}(2.5,1.5)(4.5,1.5)
  \psline[linewidth=.01](7,3)(5,1)
 \psline[linewidth=.01](7,1)(6,0)
  \psframe[linecolor=white,fillstyle=solid,fillcolor=white](6.4,2.4)(6.6,2.6)
 \psframe[linecolor=white,fillstyle=solid,fillcolor=white](5.4,1.4)(5.6,1.6)
 \psframe[linecolor=white,fillstyle=solid,fillcolor=white](6.4,.4)(6.6,.6)
\psline[linewidth=.01](5,3)(5,2)
 \psline[linewidth=.01](5,2)(7,0)
 \psline[linewidth=.01](6,3)(7,2)
 \psline[linewidth=.01](5,1)(5,0)
 \psline[linewidth=.01](7,1)(7,2)
 \rput(3.5,1.7){R3}
 \rput(3.5,-.6){c) Third}
\end{pspicture}
\]
\caption{The Reidemeister moves}\label{REID}
\end{figure}
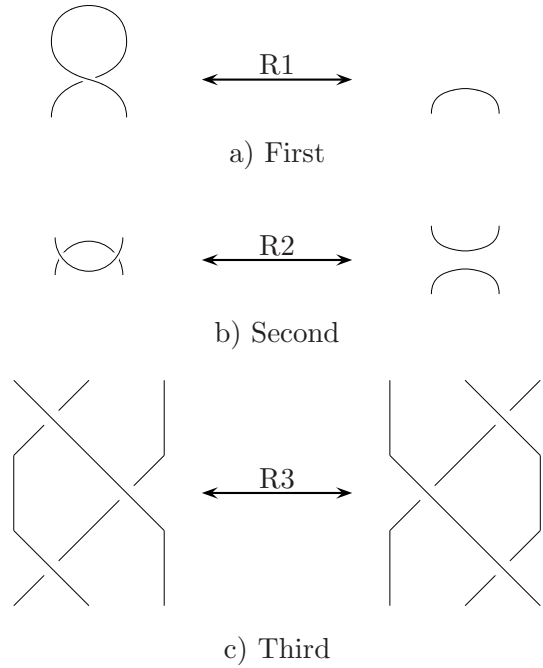

\begin{thm} Two diagrams represent isotopic links iff they can be
related via planar isotopy and a finite number of Reidemeister
moves R1, R2 and R3 shown in Figure~\ref{REID}.
\end{thm}

Given a diagram of a link there is an associated presentation
(called the Wirtinger presentation) of the link group obtained as
follows:
\begin{itemize}
\item[{\bf Generators}] There is one generator for each arc of the
diagram representing a lasso that goes around that arc.
\item[{\bf Relations}] There is one relation at each crossing
asserting that the (generator corresponding to) the over arc
conjugates one of the under arcs to the other.
\end{itemize}

\begin{rem} See~\cite{R} for a precise statement of the conjugation relations. For
the purposes of this work the above vague description is enough.
\end{rem}

\begin{rem}\label{precoliso} Notice that one can get a precise description of
the isomorphism between link groups induced by an isotopy by
tracing how the Wirtinger presentation changes under Reidemeister
moves.
\end{rem}

\section{Heegaard splittings}\label{Heegaard}
\begin{defn} Let $M$ be a $3$-manifold with boundary $\partial M$.
We say that a manifold $M'$ is obtained by $M$ by attaching a
$1$-handle iff $M'$ is homeomorphic to the space obtained by
attaching a ``solid cylinder'', $I\times D^2$ to $M$ via a map
$$\xymatrix@1{ {\{0,1\}\times D^2}\ar[r]& {\partial M}}$$
where $\{0,1\}=\partial I$.
\end{defn}

\begin{defn}A ($3$-dimensional) handlebody $H$ is a $3$-manifold
obtained by attaching some $1$-handles to the $3$-ball $D^3$. The
genus of $H$ is the genus of the surface $\partial H$.
\end{defn}

\begin{defn} A genus $g$ Heegaard splitting of a closed orientable
 $3$-manifold
$M$ is a decomposition
       $$M=H_1\bigcup_{\phi} H_2$$
 where $H_1$ and $H_2$ are genus $g$ handlebodies
 and $\phi$, called the splitting
 homeomorphism, is a homeomorphism $\xymatrix@1{ {\phi:\partial
 H_1}\ar[r]&{\partial H_2}}$.\\
 Two Heegaard splittings $H_1\bigcup_{\phi} H_2$ and
 $H_1'\bigcup_{\phi'} H_2'$ are called equivalent if there exist
 homeomorphisms $f_1$, $\bar{f_1}$ and $f_2$, $\bar{f_2}$ so that
 the following diagram commutes:
$$
\xymatrix{ {H_1}\ar[d]_{\bar{f_1}}& {\supset}&{\partial
H_1}\ar[d]_{f_1}\ar[r]^{\phi}&{\partial H_2}
\ar[d]^{f_2}& {\subset}&{H_2}\ar[d]^{\bar{f_2}}\\
{H_1'}& {\supset}&{\partial H_1'}\ar[r]^{\phi'}&{\partial H_2'}&
{\subset} &{H_2'} }
$$
 \end{defn}\

\begin{thm} Every closed orientable $3$-manifold $M$ admits a Heegaard
splitting.
\end{thm}
\begin{proof}[Sketch of proof]  $M$  admits a triangulation. Take
a triangulation of $M$ and observe that $M$ is the union of
thickenings of the $1$-skeleton of this triangulation and the
$1$-skeleton of the dual cellular decomposition. Both thickenings
are handlebodies. See~\cite{R} for details.
\end{proof}

Of course the same manifold admits many Heegaard splittings. For
example given Heegaard splittings of two manifolds $M_1$ and $M_2$
there is a connected sum Heegaard splitting of the connected sum
 $M_1\# M_2$. Therefore taking connected sums with splittings of
 the $3$-sphere produces new splittings of the same manifold.
 Theorem~\ref{HS} below asserts that essentially this is all that
 can happen.

 \begin{defn} The {\em standard} genus $1$ splitting of the $3$-sphere
 is the following:

 $$S^3=S^1\times D^2\bigcup_{S^1\times S^1} D^2\times S^1.$$

That is the splitting of $S^3$ as two solid tori (genus $1$
handlebodies) glued via the map of the torus $S^1\times S^1$ to
itself that interchanges the coordinates.\\
Taking connected sum with the standard genus $1$ splitting of
$S^3$ is referred to as {\em stabilization}.

\end{defn}

\begin{thm}\label{HS} Any two Heegaard splittings of the same $3$-manifold
are stably equivalent, that is they become equivalent after enough stabilizations. We emphasize that usually {\em
both} splittings need to be stabilized before they become equivalent.
\end{thm}
\begin{proof} See~\cite{Wal}.\end{proof}

\newpage

\end{document}